\documentstyle[amssymb,amscd,epsfig]{amsart}

\newtheorem{theorem}{Theorem}[part]
\newtheorem{corollary}[theorem]{Corollary}

\newtheorem{lemma}[theorem]{Lemma}
\newtheorem{proposition}[theorem]{Proposition}
\newtheorem{remark}[theorem]{Remark}
\newtheorem{example}[theorem]{Example}
\newtheorem{definition}[theorem]{Definition}
\theoremstyle{definition}
\theoremstyle{remark}
\numberwithin{section}{part}
\numberwithin{equation}{part}

\begin{document}
\pagenumbering{roman}

\title{On the work of Givental relative to mirror symmetry}
\author{G. Bini}
\author{C. De Concini}
\author{M. Polito}
\author{C. Procesi}

\address{Scuola Normale Superiore\\
P.zza dei Cavalieri,7\\
56126 Pisa, Italia}
\address{\noindent Dipartimento di Matematica ``Istituto Guido Castelnuovo''\\
Universit\`a di Roma ``La Sapienza''\\
P.zzale A. Moro, 2\\
00185 Roma, Italia
\bigskip}

\email[Gilberto Bini]{bini@@cibs.sns.it}
\email[Corrado De Concini]{deconcin@@mat.uniroma1.it}
\email[Marzia Polito]{polito@@cibs.sns.it}
\email[Claudio Procesi]{claudio@@mat.uniroma1.it}

\begin{abstract}
These are the informal notes of two seminars
held at the Universit\`a di Roma ``La Sapienza'', and at the Scuola
Normale Superiore in Pisa in Spring and Autumn 1997.

We discuss in detail the content of the parts of   Givental's paper \cite{G1}
dealing with
mirror symmetry for projective complete intersections.
\end{abstract}

\maketitle
\tableofcontents

\newpage
\pagenumbering{arabic}
\part{Introduction}

Let $X$ be a generic quintic threefold in ${\Bbb P}^4$. In 1991 Candelas, de
la Ossa, Green and Parkes \cite{COGP} ``predicted'' the numbers $n_d$ of
degree $d$ rational curves in $X$ conjecturing that the generating
function 
\begin{equation}
K(q)=5+\sum_{d=1}^\infty n_dd^3\frac{q^d}{1-q^d}  \label{kappa}
\end{equation}
could be recovered via elementary transformations from the hypergeometric
series 
\begin{equation*}
\sum_{d=0}^\infty \frac{(5d)!}{(d!)^5}q^d\text{,} \label{hyper} 
\end{equation*}
which is annihilated by the Picard-Fuchs linear differential operator 
\begin{equation}
D=\left( \frac d{dt}\right) ^4-5\exp \left( t\right) \prod_{m=1}^4\left(
5\frac d{dt}+m\right) \text{.}  \label{pieffe}
\end{equation}

\noindent Indeed, pick a basis $\left\{ f_i(t)\right\} _{i=0,...,3}$ of the space of 
solutions of $D(f)=0$, introduce the new variable $T(t)=\frac{f_1(t)}{f_0(t)}$, and
consider functions $g_i(T):=\frac{f_i(t(T))}{f_0(t(T))}$, $i=0,...,3$ .
These functions form a basis of solutions of the following differential
equation: 
\begin{equation*}
\left( \frac d{dt}\right) ^2\frac 1{K(\exp t)}\left( \frac d{dt}\right)
^2g(t)=0\text{,} 
\end{equation*}
where $K(q)$ is the function in \ref{kappa}.

This conjecture was motivated by a fascinating phenomenon which is known by
physicists as ``mirror symmetry''. By this they mean that given $X$ as above,
it is possible to construct a one parameter family $Y_t$ of Calabi-Yau
threefolds which are {\it mirrors} of the quintic $X$, i.e. the Hodge
numbers of each $Y_t$ enjoy the property $h^{r,s}(Y_t)=h^{3-r,s}(X)$. In
this case, $Y_t$ may be described as follows. Take the family of quintics 
\begin{equation*}
Y_\lambda ^{\prime }=\left\{ \left[ y_0\cdot ...\cdot y_4\right] \in {\Bbb P}%
^5:y_0\cdot ...\cdot y_4=\left( \lambda \right) ^{\frac 15}\left(
y_0^5+...+y_4^5\right) \right\} \text{,} 
\end{equation*}
act by multiplications of variables by fifth roots of unity, and resolve the
singularities to obtain 
\begin{equation*}
Y_t=\left\{ \left( u_0,...,u_4\right) \in {\Bbb C}^5:u_0\cdot ...\cdot
u_4=\exp (t),u_0+...+u_4=1\right\} \text{.} 
\end{equation*}

\noindent So in this context, there seems to arise a relation between
enumerative geometry in $X$ and solutions of differential equations for
periods of the family $\{Y_t\}$.

\medskip\ \ 

It is also possible to extend observations made for a generic quintic to a
broad class of complete intersections varieties in projective spaces and to
toric varieties. In the sequel, we shall assume that $X$ is a smooth
projective complete intersection of codimension $r$ and degree $%
l=(l_1,...,l_r)$ in ${\Bbb P}^n$, with $\sum_{i=1}^rl_i\leq n+1$. As for the
quintic threefold, we consider a linear differential operator $D_\hbar $ of
order $n+1-r$ which depends on a parameter $\hbar $. We still refer to $%
D_\hbar $ as the Picard-Fuchs operator associated with $X$. The explicit
form of this operator is: 
\begin{equation*}
D_\hbar =\left( \hbar \frac d{dt}\right) ^{n+1-r}-\exp \left( t\right)
\prod_{j=1}^rl_j\prod_{m=1}^{l_j-1}\left( l_j\hbar \frac d{dt}+m\hbar
\right) \text{.} 
\end{equation*}
Note that for $n=4$ and $r=1$, $D_\hbar $ coincides with the operator \ref
{pieffe} up to the factor $\hbar ^4$. Eventually, in \cite{G2}, Givental
explains how to construct a family $Y_t$ of {\it mirrors } from solutions of
the differential equation $D_\hbar (f(t,\hbar ))=0$.

\noindent Even in this more general setting, if we pick solutions for the differential
equation $D_\hbar (f(t,\hbar ))=0$, and manipulate them in a suitable
way (we shall give details later), then it is conjectured that we obtain
functions whose coefficients contain numbers $n_d$ which {\it virtually}
count rational curves on $X$. In fact, under the assumptions that the number
of degree $d$ curves in $X$ is finite, the numbers $n_d$ coincides exactly
with the number of these curves (cfr.\cite{Ma}).

\medskip\ 

Givental 's goal is to {\em formulate in a proper way and prove the link
between solutions of Picard-Fuchs equation }$D_\hbar (f(t,\hbar ))=0$ {\em %
and numbers of rational curves on }$X${\em \ }.

\medskip\ 

Let us briefly sketch Givental 's strategy. The numbers counting rational
curves on $X$ (i.e. the Gromov Witten invariants of $X$) are
all contained in a function called {\em potential}, of the following form: 
\begin{equation*}
\phi (t_0,...,t_m)=\sum_{n_0+...+n_m\geq 3}\frac{t_0^{n_0}}{n_0!}...\frac{%
t_m^{n_m}}{n_m!}\sum_{\beta \in H_2(X,{\Bbb {Z}})}I_\beta (n_0,...,n_m). 
\end{equation*}

\noindent We construct a family $\left\{ \nabla _\hbar \left( s\right) \right\} $ of
connections on the tangent bundle $TH^{*}(X)$ to the cohomology of $X$: if we
fix a basis $\left\{ T_0,...,T_m\right\} $ of $H^{*}(X)$, and coordinates $%
\left\{ t_0,...,t_m\right\} $, then we can take $\left\{ T_0,...,T_m\right\} 
$ as a global frame for the trivial bundle $TH^{*}(X)$. The connections are
defined as follows: 
\begin{equation*}
\nabla _\hbar (T_i)=-\frac 1\hbar \sum_{k=0}^m\phi _{ijk}dt_jT_k\text{,} 
\end{equation*}
i.e., the structure constants of the connections are the third derivatives of
the potential.

\noindent Pick a basis of solutions of $\nabla _\hbar \left( s\right) =0$, and
consider their components along $T_0$; the claim is that these functions $%
\left\{ s_i\right\} $, suitably ``manipulated'' (as in the case of the
quintic) form a basis of solutions for the Picard-Fuchs equation.

Let us summarize: 
\begin{equation*}
\begin{array}{ccc}
X &  & X \\ 
\downarrow &  & \downarrow \\ 
\text{Mirror family }\left\{ Y_t\right\} &  & \text{Potential (numbers of
rational curves on }X\text{)} \\ 
\downarrow &  & \downarrow \\ 
\text{Picard Fuchs operator }D_\hbar &  & \text{Connections }\nabla _\hbar
\\ 
\downarrow &  & \downarrow \\ 
\text{Solutions of }D_\hbar f=0 & \overset{\text{manipulating}}{%
\longleftrightarrow } & \text{Solutions of }\nabla _\hbar g=0
\end{array}
\text{.} 
\end{equation*}

An important remark is that $\nabla _\hbar \left( s\right) =0$ is a system
of first order differential equations, and a basis of solutions gives us
complete information about the coefficients of the system, i.e. about the
potential; hence the second column of the previous diagram can be also
interpreted in the other way: by this we mean that from solutions of
Picard-Fuchs equation of the mirror family, we can recover the numbers of rational
curves on $X$. The diagram should better look like:

\begin{equation*}
\begin{array}{ccc}
X &  & X \\ 
\downarrow &  & \downarrow \\ 
\text{Mirror family }\left\{ Y_t\right\} &  & \text{Potential (numbers of
rational curves on }X\text{)} \\ 
\downarrow &  & \downarrow \uparrow \\ 
\text{Picard Fuchs operator }D_\hbar &  & \text{Connections }\nabla _\hbar
\\ 
\downarrow &  & \downarrow \uparrow \\ 
\text{Solutions of }D_\hbar f=0 & \overset{\text{manipulating}}{%
\longleftrightarrow } & \text{Solutions of }\nabla _\hbar g=0
\end{array}
\text{.} 
\end{equation*}

\bigskip
\bigskip
\bigskip
\bigskip

In writing these notes we followed very closely Givental's paper \cite{G1}.
For a different approach one can see \cite{LLY}; also the recent 
Sem. Bourbaki talk  \cite{Pa}
by R. Pandharipande contains a discusssion of these matters and further 
information on the literature. We are  very grateful to him
for reading a preliminary version of our notes and making many useful observations.

We heartily thank Enrico Arbarello for the extremely helpful conversations 
and suggestions during the development of our work.

\newpage

\part{Preliminaries}

\section{Quantum cohomology and small quantum cohomology\label{smallq}}

\subsection{The case of convex varieties}

Let $X$ be a complex projective algebraic variety which satisfies the
following properties:

\begin{itemize}
\item  $X$ is convex, i.e. for every morphism $\mu :{\Bbb P}^1\rightarrow X$%
, $H^1({\Bbb P}^1,\mu ^{*}(T_X))=0$, with $T_X$ the tangent bundle to $X$,

\item  $A^{*}(X)\cong H^{2*}(X,{\Bbb {Z})}$ , i.e. the even integer
cohomology is isomorphic to the Chow ring,

\item  the effective cone in $H_2(X,{\Bbb {Z})}$ is positively generated by
a finite number of effective classes.
\end{itemize}

For the purpose of what follows, we may assume that $X$ is a homogeneous
space $G/P$, where $G$ is a semisimple Lie group and $P$ is a parabolic
subgroup. Therefore the second hypothesis is automatically satisfied.
Moreover, we shall deal with ${\Bbb P}^n$ for explicit calculations.

\noindent Under these hypotheses, we can build the moduli space of stable maps of $n$
pointed genus $0$ curves to $X$, whose image lies in the homology class $%
\beta \in H_2(X,{\Bbb {Z})}$, denoted by $\overline{{\cal M}}_{0,n}(X,\beta
) $ , which is an orbifold, and a locally normal projective variety of pure
dimension $\dim X+\int_\beta c_1(TX)+n-3$. Together with the moduli space we
get the evaluation maps: 
\begin{equation*}
\begin{array}{cccc}
\rho _i: & \overline{{\cal M}}_{0,n}(X,\beta ) & \rightarrow & X \\ 
& \left[ C,x_1,...,x_n,\psi \right] & \rightarrow & \psi (x_i)
\end{array}
\text{.} 
\end{equation*}

\noindent When we consider vector bundles on $\overline{{\cal M}}_{0,n}(X,\beta )$, we
always assume them to be bundles in the {\sl orbifold sense}, if not
explicitely stated.

Let $\gamma _1,...,\gamma _n$ be classes in $A^{*}(X)$, and $\beta \in H_2(X,%
{\Bbb {Z})}$; we define the Gromov Witten invariant:

\begin{equation*}
I_\beta (\gamma _1,...,\gamma _n):=\int_{\overline{{\cal M}}_{0,n}(X,\beta
)}\rho _1^{*}(\gamma _1)\cup ...\cup \rho _n^{*}(\gamma _n)\text{.} 
\end{equation*}

Three properties of these invariants are needed:

\begin{enumerate}
\item  $I_0(\gamma _1,...,\gamma _n)=\left\{ 
\begin{array}{c}
\int_X\gamma _1\cup \gamma _2\cup \gamma _3\text{ if }n=3 \\ 
0\text{ if }n>3
\end{array}
\right. $

\item  $I_\beta (1,\gamma _2,...,\gamma _n)=\left\{ 
\begin{array}{c}
\int_X\gamma _2\cup \gamma _3\text{ if }n=3\text{ and }\beta =0 \\ 
0\text{ otherwise}
\end{array}
\right. \label{unity}$

\item  if $\gamma _1\in A^1(X)$, then $I_\beta (\gamma _1,...,\gamma
_n)=\int_\beta \gamma _1\cdot I_\beta (\gamma _2,...,\gamma _n).$
\end{enumerate}

\noindent Let us fix a basis $\left\{ T_0,T_1,...,T_p,T_{p+1},...,T_m\right\} $ of $%
A^{*}(X)$ as a ${\Bbb {Z}-}$module, such that $T_0=1$, $<T_1,...,T_p>=A^1(X)$
. Consider the vector space 

\noindent $QH^{*}(X)=A^{*}(X)\otimes {\Bbb {Q}}\left[
\left[ t_0,...,t_m\right] \right] $.

\noindent By means of the invariants, it is possible to define a product in this space
so that it becomes the quantum cohomology ring of $X$.

\noindent Let $\gamma =t_0T_0+...+t_mT_m$ be a generic element of $A^{*}(X)\otimes 
{\Bbb {Q}}\left[ \left[ t_0,...,t_m\right] \right] $; the potential is
defined as follows:

\begin{eqnarray*}
\phi (t_0,...,t_m) &:&=\sum_{n\geq 3}\frac 1{n!}\sum_{\beta \in H_2(X,{\Bbb {%
Z}})}I_\beta \left( \underset{n}{\underbrace{\gamma ,...,\gamma }}\right) =
\\
&=&\sum_{n_0+...+n_m\geq 3}\frac{t_0^{n_0}}{n_0!}...\frac{t_m^{n_m}}{n_m!}%
\sum_{\beta \in H_2(X,{\Bbb {Z}})}I_\beta \left( \underset{n_0}{\underbrace{%
T_0,...,T_0}},...,\underset{n_m}{\underbrace{T_m,...,T_m}}\right) .
\end{eqnarray*}
The third derivatives, which can be easily computed, have the following
expression:

\begin{eqnarray*}
\phi _{ijk}(t_0,...,t_m) &=&\frac{\partial ^3\phi }{\partial t_i\partial
t_j\partial t_k}=\\
&=&\sum_{n\geq 0}\frac 1{n!}\sum_{\beta \in H_2(X,{\Bbb {Z}}%
)}I_\beta \left( T_i,T_j,T_k,\underset{n}{\underbrace{\gamma ,...,\gamma }}%
\right) = \\
&=&\sum_{n_0+...+n_m\geq 0}\frac{t_0^{n_0}}{n_0!}...\frac{t_m^{n_m}}{n_m!}%
\cdot \\
& &\cdot \sum_{\beta \in H_2(X,{\Bbb {Z})}}I_\beta \left( T_i,T_j,T_k,\underset{n_o}{%
\underbrace{T_0,...,T_0}},...,\underset{n_m}{\underbrace{T_m,...,T_m}}%
\right) \text{.}
\end{eqnarray*}
If we denote by $g^{ij}$ the inverse of the intersection matrix $%
g_{ij}=\int_XT_i\ \cup T_j$, we finally define the quantum product:

\begin{equation*}
T_i*T_j:=\sum_{k,l=0}^m\phi _{ijk}g^{kl}T_l\text{.} 
\end{equation*}

\begin{proposition}
\label{ass}$(QH^{*}(X),*)$ is a commutative, associative algebra with unit $%
T_0$.
\end{proposition}

The next goal is to define the small quantum cohomology algebra, which,
roughly speaking, involves only the potential restricted to the cohomology
classes in $A^1(X)$; more precisely, we define the new product structure
constants:

\begin{equation*}
\psi _{ijk}:=\phi _{ijk}(0,t_1,...,t_p,0,...,0)\text{.} 
\end{equation*}
Using properties of GW invariants, we see that

\begin{eqnarray*}
\phi _{ijk}(0,t_1,...,t_p,0,...,0)&=& \\
&=& \sum_{n_1+...+n_p\geq 0}\frac{t_1^{n_1}}{%
n_1!}...\frac{t_p^{n_p}}{n_p!}\cdot \\
& &\cdot \sum_{\beta \in H_2(X,{\Bbb {Z}})}I_\beta
\left( T_i,T_j,T_k,\underset{n_1}{\underbrace{T_1,...,T_1}},...,\underset{n_p%
}{\underbrace{T_p,...,T_p}}\right) = 
\end{eqnarray*}

\begin{eqnarray*}
=\int_XT_i\cup T_j\cup T_k+ \\
+ \sum\begin{Sb} \beta \in H_2(X,{\Bbb {Z}}) \\ %
\beta \neq 0  \end{Sb} \sum_{n_1+...+n_p\geq 0}\frac{t_1^{n_1}}{n_1!}...%
\frac{t_p^{n_p}}{n_p!}\left( \int_\beta T_1\right) ^{n_1}...\left(
\int_\beta T_p\right) ^{n_p}I_\beta (T_i,T_j,T_k)= 
\end{eqnarray*}

\begin{eqnarray*}
=\int_X T_i\cup T_j\cup T_k+\\
+ \sum_{\beta \neq 0}\left( \sum_{n_1+...+n_p\geq
0}\frac{t_1^{n_1}}{n_1!}...\frac{t_p^{n_p}}{n_p!}\left( \int_\beta
T_1\right) ^{n_1}...\left( \int_\beta T_p\right) ^{n_p}\right) I_\beta
(T_i,T_j,T_k)= 
\end{eqnarray*}

\begin{equation*}
=\int_XT_i\cup T_j\cup T_k+\sum_{\beta \neq 0}(\exp (t_1\int_\beta
T_1+...+t_p\int_\beta T_p))I_\beta (T_i,T_j,T_k)= 
\end{equation*}

\begin{equation*}
=\int_XT_i\cup T_j\cup T_k+\sum_{\beta \neq 0}q_1^{\int_\beta
T_1}...q_p^{\int_\beta T_p}I_\beta (T_i,T_j,T_k). 
\end{equation*}
where we have set $q_i=\exp (t_i)$ , for $i=1...p$.

\noindent Once more, setting

\begin{equation*}
SQH^{*}(X):=A^{*}(X)\otimes {\Bbb {Q}}\left[ \left[ t_1,...,t_p\right]
\right] 
\end{equation*}
and

\begin{equation*}
T_i*T_j:=\sum_{k,l=0}^m\psi _{ijk}g^{kl}T_l 
\end{equation*}
one gets immediatly from Proposition \ref{ass}

\begin{corollary}
$(SQH^{*}(X),*)$ is a commutative, associative algebra with unit $T_0$.
\end{corollary}

Although easier to deal with, the small quantum cohomology ring does not
provide relevant information on the enumerative geometry of $X$.

\begin{example} 
\begin{equation*}
SQH^{*}({\Bbb {P}}^r)\cong \frac{{\Bbb Q}\left[ p,q\right] }{\left(
p^{r+1}-q\right) }. 
\end{equation*}
\end{example}

\noindent Let us work out explicitely this example: first of all, we choose the basis
of $H^{*}({\Bbb P}^r)=\left\langle T_0,...,T_r\right\rangle $, where $%
T_0=1,T_1=p$ is the dual of the hyperplane class, and $T_i=p^i$.

\noindent The product formula reduces to $T_i*T_j=\sum_{k=0}^m\psi _{ijk}T_{r-k}$ .

\noindent About the computation of GW invariants, we observe that, taking the Poincar\'e
dual $\alpha $ of $p$ as a generator of $H_2(X,{\Bbb {Z}})$ , the formula
for $\psi _{ijk}$ becomes: 
\begin{equation*}
\psi _{ijk}=\int_XT_i\cup T_j\cup T_k+\sum_{d\in {\Bbb {N}}%
}q^dI_d(T_i,T_j,T_k) 
\end{equation*}
where the positive integer $d$ means that we are considering the class $%
\alpha d$.

\noindent Observe that $\dim \overline{{\cal M}}_{0,3}({\Bbb P}^r,d\alpha
)=r+\int_{d\alpha }c_i(T{\Bbb P}^r)+3-3=r+d(r+1)$ ; a necessary condition
for $I_d(T_i,T_j,T_k)\neq 0$ is: $\deg (\rho _1^{*}(T_i)\cup \rho
_2^{*}(T_j)\cup \rho _3^{*}(T_k))=\dim \overline{{\cal M}}_{0,3}({\Bbb {P}}%
^r,d\alpha )$, i.e. $i+j+k=r(d+1)+d$. This is possible only for $d=1$, i.e. $%
i+j+k=2r+1$.

\noindent We have to distinguish two cases:

\begin{enumerate}
\item  $i+j\leq r$ ; in this case $\psi _{ijk}$ has only the ``classical''
part: 
\begin{equation*}
\psi _{ijk}=\int_XT_i\cup T_j\cup T_k=\left\{ 
\begin{array}{c}
1\text{ if }k=r-i-j \\ 
0\text{ otherwise}
\end{array}
\right. \text{.}
\end{equation*}

\item  $i+j>r$ ; in this case $\psi _{ijk}$ has only the ``quantum '' part: 
\begin{equation*}
\psi _{ijk}=qI_1(T_i,T_j,T_k)=\left\{ 
\begin{array}{c}
qI_1(T_i,T_j,T_{2r+1-i-j})\text{ if }k=2r+1-i-j \\ 
0\text{ otherwise}
\end{array}
\right. \text{.}
\end{equation*}
\end{enumerate}

\noindent It is sufficient to evaluate $I_1(T_i,T_j,T_{2r+1-i-j})$ for $i=1,j=r$ ,
since the associativity law brings all other information. With the
enumerative interpretation, the invariant $I_1(T_1,T_r,T_r)=\int_\alpha
T_1\cdot I_1(T_r,T_r)$ $=I_1(T_r,T_r)$ counts the lines on ${\Bbb {P}}^r$
passing through two generic point, and equals $1$.

\noindent Finally, we can establish that $p^{r+1}=p*p^r=\psi _{1rr}=q$, and observe that
this is the only relation in the small quantum cohomology algebra.

\bigskip\ 

\subsection{The case of complete intersections}

Let $X$ be a smooth projective complete intersection, of codimension $r$ and
degree $l=\left( l_1,...,l_r\right) $, $X\overset{i}{\hookrightarrow }$ $%
{\Bbb P}^n$; generically, $X$ does not satisfy our hypothesis of convexity,
and up to now we cannot define its quantum cohomology ring; the next goal is
to solve this problem.

\noindent Let $V:=i^{*}\left( H^{*}\left( {\Bbb P}^n\right) \right) \subseteq A^{*}(X)$%
; $V$ is a subspace of dimension $m+1=n+1-r$ of the even cohomology of $X$;
fix a basis $\left\{ T_0,...,T_m\right\} $ of $V$, with coordinates $\left\{
t_0,...,t_m\right\} $.

\noindent We will define Gromov-Witten invariants on $X$ only for classes in $V$, and
then construct the quantum cohomology ring of $X$ defining a product in $%
QH^{*}\left( X\right) :=V\otimes {\Bbb Q}\left[ \left[ t_0,...,t_m\right]
\right] $ $.$

Consider the vector bundle $W_{k,d\text{ }}$on $\overline{{\cal M}}_{0,k}(%
{\Bbb P}^n,d)$ whose fibers are 
\begin{equation*}
W_{k,d}\left( \left[ C,x_1,...,x_k,\psi \right] \right) :=H^0\left( C,\psi
^{*}\left( \oplus_{i=1}^r{\cal O}\left( l_i\right) \right) \right) \text{,%
} 
\end{equation*}
and let ${\cal E}_{k,d}:=Euler\left( W_{k,d}\right) $.

\noindent Now we are ready to define the GW invariants for $X$; let $\left( \gamma
_1,...,\gamma _k\right) \in V^k$, and choose $\left( \delta _1,...,\delta
_k\right) \in \left( H^{*}\left( {\Bbb P}^n\right) \right) ^k$ such that $%
\gamma _1=i^{*}\left( \delta _1\right) ,...,\gamma _k=i^{*}\left( \delta
_k\right) $ ; fix $\beta \in H_2(X)$ such that $i_{*}\left( \beta
\right) =d\alpha $, then:

\begin{definition}
\label{gwforx} 
\[
I_d\left( \gamma _1,...,\gamma _k\right) :=\int_{\overline{{\cal M}}_{0,k}(%
{\Bbb P}^n,d)}\rho _1^{*}\left( \delta _1\right) \cup ...\cup \rho
_k^{*}\left( \delta _k\right) \cup {\cal E}_{k,d}\text{.}
\]
\end{definition}

\noindent We now prove that this is well defined, i.e. that $I_d\left( \gamma
_1,...,\gamma _k\right) $ does not depend on the choice of $\left( \delta
_1,...,\delta _k\right) $; for this we need the following lemma.

\begin{lemma}
$Ker$ $i^{*}\subset Ann(E)$, where $E:=Euler\left( \oplus_{i=1}^r{\cal O}%
\left( l_i\right) \right) $.
\end{lemma}

{\bf Proof. }For every $\alpha \in H^{*}\left( {\Bbb P}^n\right) $, we have
that $i_{*}\left( i^{*}\left( \alpha \right) \right) =E\cdot \alpha .$
\begin{flushright}
$\Box$
\end{flushright}

\noindent Take $\delta _1^{\prime }=\delta _1+\epsilon _1$, with $\epsilon _1\in Ker$ $%
i^{*}$. Then 
\begin{equation*}
\int_{\overline{{\cal M}}_{0,k}({\Bbb P}^n,d)}\rho _1^{*}\left( \delta
_1^{\prime }\right) \cup ...\cup \rho _k^{*}\left( \delta _k\right) \cup 
{\cal E}_{k,d}= 
\end{equation*}
\begin{equation*}
=\int_{\overline{{\cal M}}_{0,k}({\Bbb P}^n,d)}\rho _1^{*}\left( \delta
_1\right) \cup ...\cup \rho _k^{*}\left( \delta _k\right) \cup {\cal E}%
_{k,d}+\int_{\overline{{\cal M}}_{0,k}({\Bbb P}^n,d)}\rho _1^{*}\left(
\epsilon _1\right) \cup ...\cup \rho _k^{*}\left( \delta _k\right) \cup 
{\cal E}_{k,d}\text{.} 
\end{equation*}
Look at the following diagram:

\begin{equation*}
\begin{array}{ccc}
\overline{{\cal M}}_{0,k+1}({\Bbb P}^n,d) & \overset{\tau _1}{\rightarrow }
& {\Bbb P}^n \\ 
\pi \downarrow & \overset{\rho _1}{\nearrow } &  \\ 
\overline{{\cal M}}_{0,k}({\Bbb P}^n,d) &  & 
\end{array}
\text{} 
\end{equation*}
Since ${\cal E}_{k,d}=\pi _{*}\tau _1^{*}\left( E\right) $, we get 
\begin{equation*}
\int_{\overline{{\cal M}}_{0,k}({\Bbb P}^n,d)}\rho _1^{*}\left( \epsilon
_1\right) \cup ...\cup \rho _k^{*}\left( \delta _k\right) \cup {\cal E}%
_{k,d}= 
\end{equation*}
\begin{eqnarray*}
&=&\int_{\overline{{\cal M}}_{0,k}({\Bbb P}^n,d)}\rho _1^{*}\left( \epsilon
_1\right) \cup ...\cup \rho _k^{*}\left( \delta _k\right) \cup \pi _{*}\tau
_1^{*}\left( E\right) = \\
&=&\int_{\overline{{\cal M}}_{0,k+1}({\Bbb P}^n,d)}\pi ^{*}\left( \rho
_1^{*}\left( \epsilon _1\right) \cup ...\cup \rho _k^{*}\left( \delta
_k\right) \right) \cup \tau _1^{*}\left( E\right) = \\
&=&\int_{\overline{{\cal M}}_{0,k+1}({\Bbb P}^n,d)}\pi ^{*}\left( \rho
_2^{*}\left( \delta _2\right) \cup ...\cup \rho _k^{*}\left( \delta
_k\right) \right) \cup \tau _1^{*}\left( \epsilon _1E\right) =0\text{.}
\end{eqnarray*}

\smallskip\ 

We will now justify definition \ref{gwforx}. Consider the subset of $%
\overline{{\cal M}}_{0,k}({\Bbb P}^n,d)$, which consists of equivalence
classes of curves mapping to $X$: 
\begin{equation*}
S_{0,k}(X,d)=\left\{ \left[ C,x_1,...,x_k,\psi \right] \in \overline{{\cal M}%
}_{0,k}({\Bbb P}^n,d):\psi (C)\subset X\right\} \text{.} 
\end{equation*}
If $S_{0,k}(X,d)$ were a closed projective subvariety, with at most finite
quotient singularities, and if $i_{*}:H_2\left( X\right) \rightarrow
H_2\left( {\Bbb P}^n\right) $ were injective, then we could take $%
S_{0,k}(X,d)$ as the moduli space $\overline{{\cal M}}_{0,k}(X,\beta )$, and
we could claim that ${\cal E}_{k,d}$ is its fundamental class in $\overline{%
{\cal M}}_{0,k}({\Bbb P}^n,d)$. In fact, consider the section $\sigma \in
\oplus_{i=1}^r{\cal O}\left( l_i\right) $ whose locus of zeroes is $X$.
From this, build a section $\Gamma \,$ of $W_{k,d}$: $\Gamma \left( \left[
C,x_1,...,x_k,\psi \right] \right) :=\psi ^{*}\left( \sigma \right) $ ; note
that $S_{0,k}(X,d)$ is exactly the zero locus of $\Gamma $, and from this
the claim would follow.

\noindent Unfortunately, we do not know anything about $S_{0,k}(X,d)$, and we have to
take \ref{gwforx} just as a definition.

Now, we can extend all the previous results: let $\gamma =t_0T_0+...+t_mT_m$
be a generic element of $V\otimes {\Bbb {Q}}\left[ \left[
t_0,...,t_m\right] \right] $ ; we define the potential

\begin{equation*}
\phi (t_0,...,t_m)=\sum_{n_0+...+n_m\geq 3}\frac{t_0^{n_0}}{n_0!}...\frac{%
t_m^{n_m}}{n_m!}\sum_{d\in \i_{*}\left( H_2(X,{\Bbb {Z}})\right)
}I_d\left( \underset{n_0}{\underbrace{T_0,...,T_0}},...,\underset{n_m}{%
\underbrace{T_m,...,T_m}}\right) \text{,} 
\end{equation*}
consider the third derivatives:

\begin{eqnarray*}
\phi _{ijk}(t_0,...,t_m)\ =\sum_{n_0+...+n_m\geq 0}\frac{t_0^{n_0}}{n_0!}...%
\frac{t_m^{n_m}}{n_m!}\cdot \\
\cdot \sum_{d\in \i_{*}\left( H_2(X,{\Bbb {Z}})\right)
}I_d\left( T_i,T_j,T_k,\underset{n_o}{\underbrace{T_0,...,T_0}},...,%
\underset{n_m}{\underbrace{T_m,...,T_m}}\right) 
\end{eqnarray*}
and the inverse $g^{ij}$ of the intersection matrix $g_{ij}=\int_XT_i\ \cup
T_j$ , and then finally define the quantum product:

\begin{equation*}
T_i*T_j:=\sum_{k,l=0}^m\phi _{ijk}g^{kl}T_l\text{.} 
\end{equation*}

\begin{proposition}
$(QH^{*}(X),*)$ is a commutative, associative algebra with unit $T_0$.
\end{proposition}

{\bf Proof. }Commutativity is trivial, and the fact that $T_0$ is the unity
follows from property \ref{unity}. For the associativity, one can repeat
exactly the proof in \cite{FP}, with a suitable modification of Lemma 15, p.
37 we now explain. With the same notation, $D(A_1,A_2,d_1,d_2)$ is the
boundary divisor of $\overline{{\cal M}}_{0,k}({\Bbb P}^n,d)$ isomorphic to $%
\overline{{\cal M}}_{0,A_1\cup \left\{ \cdot \right\} }({\Bbb P}%
^n,d_1)\times _{{\Bbb P}^n}\overline{{\cal M}}_{0,A_2\cup \left\{ \cdot
\right\} }({\Bbb P}^n,d_2)$, $i$ is the natural inclusion of $%
D(A_1,A_2,d_1,d_2)$ in $\overline{{\cal M}}_{0,A_1\cup \left\{ \cdot
\right\} }({\Bbb P}^n,d_1)\times \overline{{\cal M}}_{0,A_2\cup \left\{
\cdot \right\} }({\Bbb P}^n,d_2)$, and $j$ the embedding of $%
D(A_1,A_2,d_1,d_2)$ in $\overline{{\cal M}}_{0,k}({\Bbb P}^n,d)$.

\noindent Let us now compute the restriction of the Euler class ${\cal E}_{k,d}$ to
the image of $j$; if $\rho _{k_i+1}^i$ is the evaluation map on the last
point of $\overline{{\cal M}}_{0,A_i\cup \left\{ \cdot \right\} }({\Bbb P}%
^n,d_i)$, we define $W_{k_i+1,d_i}^{\prime }$ as the kernel of the map 
\begin{eqnarray*}
W_{k_i+1,d_i} &\rightarrow &\rho _{k_i+1}^{i*}\left( \oplus_{i=1}^r{\cal O%
}(l_i)\right)   \\
\sigma &\rightarrow &\sigma \left( \psi \left( x_{k_i+1}\right) \right) 
\text{,}
\end{eqnarray*}
therefore ${\cal E}_{k_i+1,d_i}^{\prime }:=Euler(W_{k_i+1,d_i}^{\prime })$
satisfies ${\cal E}_{k_i+1,d_i}={\cal E}_{k_i+1,d_i}^{\prime }\rho
_{k_i+1}^{i*}\left( E\right) $. If $\nu $ denotes the evaluation map on the
meeting point in 
$$
D(A_1,A_2,d_1,d_2)\simeq \overline{{\cal M}}_{0,A_1\cup
\left\{ \cdot \right\} }({\Bbb P}^n,d_1)\times _{{\Bbb P}^n}\overline{{\cal M%
}}_{0,A_2\cup \left\{ \cdot \right\} }({\Bbb P}^n,d_2)
$$
then we get the
exact sequence: 
\begin{equation*}
\begin{array}{ccccccccc}
0 & \rightarrow & W_{k_1+1,d_1}^{\prime }\oplus W_{k_2+1,d_2}^{\prime } & 
\rightarrow & W_{k,d} & \rightarrow & \nu ^{*}\left( \oplus_i{\cal O}%
\left( l_i\right) \right) & \rightarrow & 0 \\ 
&  &  &  & \left[ \psi ^{*}\left( \sigma \right) ,\psi ^{*}\left( \tau
\right) \right] & \rightarrow & \sigma \left( x_{k_1+1}\right) =\tau \left(
x_{k_2+1}\right) &  & 
\end{array}
\text{,}  
\end{equation*}
and ${\cal E}_{k,d}={\cal E}_{k_1+1,d_1}^{\prime }{\cal E}%
_{k_2+1,d_2}^{\prime }\nu ^{*}\left( E\right) $.

\noindent Fix a suitable basis for $H^{*}\left( {\Bbb P}^n\right) $; let $\left\{
T_0,...,T_s\right\} $ be a set of independent elements such that $\left\{
\left[ T_0\right] ,...,\left[ T_s\right] \right\} $ form an orthonormal
basis in $\frac{H^{*}\left( {\Bbb P}^n\right) }{Ann(E)}$, hence $%
\int_XT_i\cup T_jE=\delta _{ij}$. Let us denote $T_iE$ with $T_i^{\prime }$
and complete $\left\{ T_0,...,T_s \right\} $ and $\left\{ T_0^{\prime
},...,T_s^{\prime }\right\} $
to orthonormal dual bases of $H^{*}\left( {\Bbb P}^n\right) $  
$\left\{ T_0,...,T_s,T_{s+1},...,T_n\right\} $ and $\left\{ T_0^{\prime
},...,T_s^{\prime },T_{s+1}^{\prime },...,T_n^{\prime }\right\} $, such
that the class of the diagonal in ${\Bbb P}^n$ is 
\begin{equation*}
\Delta =\sum_{i=0}^nT_i\otimes T_i^{\prime }\text{.} 
\end{equation*}
Putting everything toghether, and denoting with $\mu =\left( \rho
_{k_1+1}^1,\rho _{k_2+1}^2\right) $ the product of the evaluation maps on
the last marked points, we get 
\begin{equation*}
\int_{D(A_1,A_2,d_1,d_2)}\rho _1^{*}\left( \delta _1\right) ...\rho
_k^{*}\left( \delta _k\right) {\cal E}_{k,d}= 
\end{equation*}
\begin{eqnarray*}
\ &=&\int_{D(A_1,A_2,d_1,d_2)}\rho _1^{*}\left( \delta _1\right) ...\rho _k^{*}\left( \delta
_k\right) {\cal E}_{k_1+1,d_1}^{\prime }{\cal E}_{k_2+1,d_2}^{\prime }\mu
^{*}\left( E\cdot \left( 1\otimes 1\right) \right) \mu ^{*}\left( \Delta
\right) = \\
\ &=&\int_{D(A_1,A_2,d_1,d_2)}\rho _1^{*}\left( \delta _1\right) ...\rho _k^{*}\left( \delta
_k\right) {\cal E}_{k_1+1,d_1}^{\prime }{\cal E}_{k_2+1,d_2}^{\prime }\mu
^{*}\left( E\sum_{i=0}^mT_i\otimes T_i^{\prime }\right) = \\
\ &=&\int_{D(A_1,A_2,d_1,d_2)}\rho _1^{*}\left( \delta _1\right) ...\rho _k^{*}\left( \delta
_k\right) {\cal E}_{k_1+1,d_1}^{\prime }{\cal E}_{k_2+1,d_2}^{\prime }\mu
^{*}\left( \left( E\otimes E\right) \sum_{i=0}^mT_i\otimes T_i\right) =
\\
&=&\sum_{i=1}^n\int_{\overline{{\cal M}}_{0,A_1\cup \left\{ \cdot \right\} }(%
{\Bbb P}^n,d_1)}\rho _1^{*}\left( \delta _1\right) ...\rho _{k_1}^{*}\left(
\delta _k\right) \rho _{k_1+1}^1\left( T_i\right) {\cal E}_{k_1+1,d_1}\ \cdot
\\
&&\cdot \int_{\overline{{\cal M}}_{0,A_2\cup \left\{ \cdot \right\} }({\Bbb P%
}^n,d_2)}\rho _{k_1+1}^{*}\left( \delta _{k_1+1}\right) ...\rho _k^{*}\left(
\delta _k\right) \rho _{k_2+1}^2\left( T_i\right) {\cal E}_{k_2+1,d_2}\ 
\text{.}
\end{eqnarray*}
\begin{flushright}
$\Box$
\end{flushright}

In the same way we can define the small quantum cohomology ring: $SQH^{*}\left( X\right)
:=V\otimes {\Bbb Q}\left[ q\right] $, 
\begin{equation*}
\psi _{ijk}\left( q\right) :=\sum_{d\in \iota _{*}\left( H_2(X,{\Bbb {Z}}%
)\right) }q^dI_d(T_i,T_j,T_k)\text{,} 
\end{equation*}
and $T_i*T_j:=\sum_{k,l=0}^m\psi _{ijk}g^{kl}T_l$.

\begin{proposition}
$(SQH^{*}(X),*)$ is a commutative, associative algebra with unit $T_0$.
\end{proposition}

\begin{example}
In the case of a quintic threefold $X$ in ${\Bbb P}^4$, $%
SQH^{*}(X)=QH^{*}(X).$
\end{example}

\noindent This follows from a dimension computation: in fact, $\dim _{{\Bbb R}}\overline{%
{\cal M}}_{0,k}({\Bbb P}^4,d)=10d+2k+2$, and $\deg {\cal E}_{k,d}=10d+2$,
hence the non zero GW involve only degree two classes.

\bigskip
\bigskip

\section{Equivariant cohomology for torus actions\label{coequiv}}

Let $X$ be a topological space and $G\cong (S^1)^k$ be a compact torus
acting on $X$. Recall the following definition:

\begin{definition}
The universal bundle for a compact Lie group $G$ is the $G$-fibering $%
EG\rightarrow BG=EG/G$, where $EG$ is a contractible topological space on which 
$G$ acts freely, and $BG$ is called classifying space.
\end{definition}

\begin{remark}
The universal bundle satisfies the universal property that every $G$%
-principal bundle is obtained from it by pullback.
\end{remark}

\noindent In the case of the torus, a model for the universal bundle is the Hopf
fibering of infinite dimensional spheres on infinite dimensional complex
projective spaces:

\begin{equation*}
\begin{array}{ccc}
E(S^1)^k & \cong & (S^\infty )^k \\ 
\downarrow &  &  \\ 
B(S^1)^k & \cong & ({\Bbb {P}}^\infty )^k
\end{array}
\text{.} 
\end{equation*}
The group $G$ acts on $X\times EG$, and we denote by $X_G=X\times _GEG$ the
quotient; it is a bundle over $BG\,$, with fiber $X$:

\begin{equation*}
\begin{array}{cccc}
\pi : & X_G & \rightarrow & BG \\ 
& \left[ x,e\right] &  & \left[ e\right]
\end{array}
\text{.} 
\end{equation*}

\begin{definition}
The $G$-equivariant cohomology of $X$ is the cohomology of $X_G$: 
\[
H_G^{*}(X):=H^{*}(X_G).
\]
\end{definition}

\noindent From now on we will consider cohomology with coefficients in ${\Bbb {Q}}$.

\noindent Note that $H_G^{*}(\left\{ pt\right\} )=H^{*}(BG)$ . Therefore the
equivariant cohomology ring is a module (via $\pi ^{*}$) over this ring.

\noindent In our case, $H_{S^1}^{*}(\left\{ pt\right\} )=H^{*}(B(S^1))=H^{*}({\Bbb {P%
}}^\infty )\cong {\Bbb {Q}\left[ {\lambda }\right] }$ , and $%
H_{(S^1)^k}^{*}(\left\{ pt\right\} )\cong {\Bbb Q}\left[ \lambda
_1,...,\lambda _k\right] $.

The following example, and others which will come later in this section, are
the most relevant from our point of view, since Givental uses the
equivariant theory almost exclusively in these cases.

\begin{example}
$X={\Bbb {P}}^1,$ $G=(S^1)$\label{p1s1}.
\end{example}

\noindent The action is: 
\begin{equation}
\begin{array}{ccccc}
G & \times & X & \rightarrow & X \\ 
t & , & \left[ x_0,x_1\right] &  & \left[ t^2x_0,x_1\right]
\end{array}
\text{.}  \label{esempio}
\end{equation}
The bundle ${\Bbb {P}}_{S^1}^1{\Bbb ={P}}^1\times _{S^1}S^\infty \rightarrow 
{\Bbb {P}}^\infty $ is nothing else but the projectivization of a two
dimensional complex vector bundle on ${\Bbb {P}}^\infty $, namely ${\Bbb {P}}%
({\cal O}\left( -2\right) \oplus {\cal O})\rightarrow {\Bbb {P}}^\infty $.

\noindent Assuming that the cohomology of ${\Bbb P}^\infty $ is generated by $%
c_1\left( {\cal O}\left( -2\right) \right) $, by the Whitney formula, we can
compute:

\begin{equation*}
H_G^{*}(X)=H^{*}\left( {\Bbb P}\left( {\cal O}\left( -2\right) \oplus 
{\cal O}\right) \right) \cong \frac{{\Bbb Q}\left[ \lambda ,p\right] }{%
\left( p^2-p\lambda \right) }\text{,} 
\end{equation*}
where $p$ is to be interpreted as the generator of the cohomology of the
fiber of the projective bundle ${\Bbb {P}}({\cal O}\left( -2\right)
\oplus {\cal O})\rightarrow {\Bbb {P}}^\infty $.

\begin{example}
\label{pnsk}A generalization to a diagonal action of a $k$ dimensional torus 
$G$ on ${\Bbb {P}}^n$, with characters $\chi _i\in {\Bbb {Z}}^k$ , $i=0,...,n
$.
\end{example}

\noindent Let $L_j$ be the pull backs of the hyperplane bundle on the different 
components of 
$({\Bbb {P}}^\infty )^k$, and take the Chern classes of these bundles as generators 
of the cohomology of ${\Bbb ({P}}^\infty )^k$; we have: 
\begin{equation*}
{\Bbb {P}}^n\times _G(S^\infty )^k\cong {\Bbb {P}}\left( \left(
\otimes_{j=1}^kL_j^{\chi _{oj}}\right) \oplus ...\oplus \left(
\otimes_{j=1}^kL_j^{\chi _{nj}}\right) \right) \rightarrow {\Bbb ({P}}%
^\infty )^k\text{.} 
\end{equation*}
Setting $\chi _i(\lambda ):=\prod_j\lambda _j^{\chi _{ij}}$, and applying
Whitney formula, 
\begin{equation*}
H_G^{*}({\Bbb {P}}^n)=\frac{{\Bbb Q}\left[ \lambda _1,...,\lambda
_k,p\right] }{\prod_{i=0,...,n}(p-\chi _i(\lambda ))}\text{.} 
\end{equation*}

\subsection{$G$ -vector bundles and equivariant characteristic classes}

Suppose we have a complex rank $r$ vector bundle $V$ over $X$ with a $G$
linear action, compatible with the projection onto $X$ and linear on the
fibers, then $V_G:=V\times _GEG$ is still a ${\Bbb {C}}^r$ bundle over $X_G$.
 We define the equivariant Chern and Euler classes of $V$ : 
\begin{eqnarray*}
c_i^G(V) &:&=c_i\left( V_G\right) \text{ } \\
{\cal E}^G(V) &:&={\cal E}(V_G)\text{ .}
\end{eqnarray*}

\noindent If we reconsider Examples \ref{p1s1} and \ref{pnsk} more closely, we realize
that the actions described there correspond to a distinguished choice of
linearization of the tautological bundles over projective spaces. If another
choice had been made, then equivariant cohomology groups would have been
isomorphic. We now give a further example which will be useful for the
purpose of what folllows.

\begin{example}
$X={\Bbb {P}}^n$, $V=\oplus_{1=i}^r{\cal O}\left( l_i\right) $, $%
G=T^{n+1}\times T^r$. We compute the equivariant Euler class.
\end{example}

\noindent $T^{n+1}$ acts on $X$ as in case \ref{pnsk}, with $k=n+1$, and $\chi
_i=(0,...,0,\underset{i}{1},0,...,0)$\label{pnsn+1}, and $T^r$ acts
trivially; therefore $H_T^{*}\left( {\Bbb {P}}^n\right) \cong \frac{{\Bbb {Q}%
}\left[ \lambda _0,...,\lambda _n,\mu _1,...,\mu _r,p\right] }{%
\prod_{i=0}^n(p-\lambda _i)}$. Conversely, $T^r$ acts on a fiber of the
bundle by diagonal action with characters $\mu _i=(0,...,0,\underset{i}{1}%
,0,...,0)$, and the $T^{n+1}$action is induced by the one on $X$.

\noindent Let's look at the following diagram:

\begin{equation*}
\begin{array}{ccc}
\oplus_{i=1}^r{\cal O}\left( l_i\right) \times _{\left( S^1\right)
^{n+1}}\left( S^\infty \right) ^{n+1}\times _{\left( S^1\right) ^r}\left(
S^\infty \right) ^r & \overset{\phi }{\rightarrow } & {\Bbb {P}}^n\times
_{\left( S^1\right) ^{n+1}}\left( S^\infty \right) ^{n+1}\times \left( {\Bbb %
P}^\infty \right) ^r \\ 
& \overset{\pi ^{\prime }}{\searrow } & \downarrow \pi \\ 
&  & \left( {\Bbb {P}}^\infty \right) ^{n+1}\times \left( {\Bbb {P}}^\infty
\right) ^r
\end{array}
\end{equation*}

Restricted to a fiber of $\pi $, which is isomorphic to ${\Bbb {P}}^n$, the
diagram looks like: 
\begin{equation*}
\begin{array}{ccc}
\oplus_{i=1}^r{\cal O}\left( l_i\right) \otimes M_i^{-1} & \overset{%
\phi }{\rightarrow } & {\Bbb {P}}^n{\Bbb \times }\left( {\Bbb P}^\infty
\right) ^r
\end{array}
\text{,} 
\end{equation*}
where $M_i$ is the pull-back of the hyperplane bundle on the $i$-th copy of $%
{\Bbb {P}}^\infty $ , and has Chern class $\mu _i$.

\noindent The Euler class of this restricted bundle is $\prod_{i=1}^r(l_ip-\mu _i)$;
using this, we can prove that ${\cal E}^G(V)=\prod_{i=1}^r(l_ip-\mu _i)$.

\subsection{Equivariant integral and pairing}

Let now $X$ be a compact manifold with at most finite quotient
singularities.The fibering $
\begin{array}{cccc}
\pi : & X_G & \rightarrow & BG
\end{array}
$ induces a push-forward map $
\begin{array}{cccc}
\pi _{*}: & H_G^{*}(X) & \rightarrow & H_G^{*}(\left\{ pt\right\} )
\end{array}
$, which we will call equivariant integral and often denote with $\int^G$.
All the same for the equivariant pairing 
\begin{equation*}
\begin{array}{cccc}
\left\langle ,\right\rangle _G: & H_G^{*}(X)\times H_G^{*}(X) & \rightarrow
& H_G^{*}(\left\{ pt\right\} ) \\ 
& (\omega ,\eta ) &  & \pi _{*}(\omega \cup \eta ) \text{.}
\end{array}
\end{equation*}

\begin{example}
Case \ref{pnsk}.
\end{example}

\noindent The equivariant integral is given by the following formula: 
\begin{equation*}
\int^Gf(p,\lambda _0,...,\lambda _n)=\frac 1{2\pi \sqrt{-1}}\int \frac{%
f(p,\lambda _0,...,\lambda _n)}{\prod_{j=0}^n(p-\chi _j\left( \lambda
\right) )}dp\text{.} 
\end{equation*}
In fact, $\pi _{*}$ is just the integration along the fiber, which is a $n$
dimensional projective space whose cohomology is generated by $p$.
Therefore, writing $f(p,\lambda _0,...,\lambda _n)=\sum_jf_j(\lambda )p^j$,
we see that $\int^Tf(p,\lambda _0,...,\lambda _n)=f_n(\lambda )$ . On the
other hand, 
\begin{equation*}
\frac 1{2\pi \sqrt{-1}}\int \frac{\sum_jf_j(\lambda )p^j}{\prod_j(p-\chi _j)}%
dp=\frac 1{2\pi \sqrt{-1}}\sum_jf_j\left( \lambda \right) \int \frac{p^j}{%
\prod_j(p-\chi _j)}dp\text{;} 
\end{equation*}
these integrals all vanish by degree computation except for $j=n$ , and the
formula follows.

\subsection{Localization at fixed points and integration formula}

One of the main tools of equivariant cohomology is the localization at fixed
points. Equivariant cohomology rings for torus actions are modules over $%
{\Bbb {Q}}\left[ \lambda _1,...,\lambda _k\right] $ , and therefore it is
possible to localize on their support. Suppose $D_1,...,D_s$ are the
connected components of fixed points for the $T$ action on $X$, then $%
H_T^{*}(D_i)=H^{*}(D_i\times _TET)=H^{*}(D_i\times BT)\cong
H^{*}(D_i)\otimes H_T^{*}\left( \left\{ pt\right\} \right) \cong
H^{*}(D_i)\otimes {\Bbb {Q}}\left[ \lambda _1,...,\lambda _k\right] $.

\begin{proposition}
There is a natural morphism 
\[
\begin{array}{ccc}
H_T^{*}(X) & \stackrel{\delta =(\delta _1,...,\delta _s)}{\rightarrow } & 
\oplus_{i=1...,s}H_T^{*}(D_i)
\end{array}
\text{,}
\]
whose kernel and cokernel are torsion modules.
\end{proposition}

\begin{example}
Case \ref{p1s1}.
\end{example}

\noindent There are two fixed points, $\left[ 0,1\right] $, and $\left[ 1,0\right] $
(say $0$ and $\infty $). The morphism is 
\begin{equation*}
\begin{array}{ccc}
H_T^{*}({\Bbb {P}}^1{\Bbb )} & \rightarrow & H_T^{*}\left( \left\{ 0\right\}
\right) \oplus H_T^{*}\left( \left\{ \infty \right\} \right) \\ 
\cong &  & \cong \\ 
\frac{{\Bbb {Q}}\left[ \lambda ,p\right] }{\left( p^2-p\lambda \right) } & 
\rightarrow & {\Bbb {Q}\left[ \lambda \right] \oplus Q\left[ \lambda
\right] } \\ 
f(p,\lambda ) &  & (f(\lambda ,\lambda ),f(0,\lambda ))
\end{array}
\text{.} 
\end{equation*}

\noindent Localizing both sides outside the ideal $\left\langle \lambda \right\rangle $,
 we get an isomorphism. In particular, we see that 
\begin{eqnarray*}
\frac p\lambda &\rightarrow &\left( 1,0\right) \\
\frac{\lambda -p}\lambda &\rightarrow &\left( 0,1\right)
\end{eqnarray*}
and that $1=\frac p\lambda +\frac{\lambda -p}\lambda $ .

\begin{example}
Case \ref{pnsk}, with $k=n+1$, and $\chi _i=(0,...,0,\underset{i}{1}%
,0,...,0)$, as in \ref{pnsn+1}.
\end{example}

\noindent There are $n+1$ fixed points, $x_i=\left[ 0,...,,0,\underset{i}{1}%
,0,...,0\right] $; the morphism is 
\begin{equation*}
\begin{array}{ccc}
\frac{{\Bbb {Q}}\left[ \lambda _0,...,\lambda _n,p\right] }{%
\prod_{i=0}^n(p-\lambda _i)} & \rightarrow & \oplus_{i=0,...,n}{\Bbb {Q}}%
\left[ \lambda _0,...,\lambda _n\right] \\ 
f(p,\lambda _0,...,\lambda _n) &  & f(\lambda _i,\lambda _0,...,\lambda _n)
\end{array}
\text{.} 
\end{equation*}

\noindent The isomorphism comes localizing both sides ouside the ideal $\left\langle
\lambda _0,...,\lambda _n\right\rangle $. Once more, it is useful to know
that, if we set $\phi _i=\frac{\prod_{j\neq i}(p-\lambda _j)}{\prod_{j\neq
i}(\lambda _i-\lambda _j)}$, then $\phi _i$ goes to $(0,...,0,\underset{i}{1}%
,0,...,0)$ .

\noindent Notice that after having localized, 
\begin{equation*}
H_G^{*}(X)_{\left\langle \lambda _0,...,\lambda _n\right\rangle }\cong \frac{%
{\Bbb {Q}}\left( \lambda _0,...,\lambda _n\right) \left[ p\right] }{%
\prod_{i=0}^n(p-\lambda _i)} 
\end{equation*}
becomes a vector space of dimension $n+1$ over the field ${\Bbb {Q}}\left(
\lambda _0,...,\lambda _n\right) $, and that the equivariant pairing gives a
scalar product. The $\phi _i$'s form an orthogonal basis with respect to
this product: if $i\neq k$, $\phi _i\cdot \phi _k$ vanishes, and on the other hand 
\begin{eqnarray*}
\int^G\phi _i^2 &=&\frac 1{2\pi \sqrt{-1}\prod_{j\neq i}(\lambda _i-\lambda
_j)^2}\int^G\frac{\prod_{j\neq i}(p-\lambda _j)^2}{\prod_j(p-\lambda _j)}dp=
\\
\ &=&\frac 1{2\pi \sqrt{-1}\prod_{j\neq i}(\lambda _i-\lambda _j)^2}\int^G%
\frac{\prod_{j\neq i}(p-\lambda _j)}{(p-\lambda _i)}dp=\frac 1{\prod_{j\neq
i}(\lambda _i-\lambda _j)}
\end{eqnarray*}
The integration formula gives us a method for reducing equivariant integrals
to integrals on fixed points components; since equivariant cohomology on
fixed points is just usual cohomology tensorized by a ring of polynomials,
the formula simplifies computations in the sense that it reduces the problem
of integrating along the fiber isomorphic to $X$ to the one of integrating
on some subvarieties of $X$; let ${\cal E}_i$ be the equivariant Euler class
of the normal bundle ${\cal N}_{D_i/X}$ of the connected component of fixed
points $D_i$ in $X$.

\begin{proposition}
The following formula holds: 
\[
\int_X^G\omega =\sum_{i=1}^s\int_{D_i}^G\frac{\delta _i(\omega )}{{\cal E}_i}%
\text{.}
\]
\end{proposition}

\begin{example}
Case \ref{p1s1}.
\end{example}

\noindent The normal bundle of $\left\{ 0\right\} $ in ${\Bbb {P}}^1$ can be
identified equivariantly with ${\Bbb {P}}^1\backslash \left\{ \infty
\right\} ={\Bbb {C}}$, and $S^1$ acts on it by multiplication by $t^2$; the
Euler class of this bundle is $\lambda $ . Similarly, the action of $S^1$ on
the normal bundle to $\left\{ \infty \right\} $ is by multiplication by $%
t^{-2}$ and the Euler class is $-\lambda $. The formula is then: 
\begin{equation*}
\int_{{\Bbb {P}}^1}^{S^1}f(p,\lambda )=\frac{f(\lambda ,\lambda )}\lambda -%
\frac{f(0,\lambda )}\lambda \text{.} 
\end{equation*}

\begin{example}
\label{xp1s1} $G=S^1$ acts on $X\times {\Bbb {P}}^1$, trivially on the first
factor, and as in \ref{p1s1} on the second.
\end{example}

\noindent This example does not bring anything new from the point of view of
equivariant cohomology, but is frequentely used by Givental.

\noindent Clearly, $H_G^{*}(X\times {\Bbb {P}}^1){\Bbb \cong }H^{*}(X){\Bbb \otimes 
}H_G^{*}({\Bbb {P}}^1)$, and the localization isomorphism tells us: 
\begin{equation*}
H^{*}(X)\otimes \frac{{\Bbb {Q}}\left( \lambda \right) \left[ p\right] }{%
\left( p^2-p\lambda \right) }\cong \left( H^{*}(X)\otimes {\Bbb {Q}}%
\left( \lambda \right) \right) {\Bbb \oplus }\left( H^{*}(X){\Bbb %
\otimes {Q}}\left( \lambda \right) \right) \text{.} 
\end{equation*}
If $x\in H_G^{*}(X\times {\Bbb {P}}^1{\Bbb )}$ and $\delta (x)=(t,\tau )$,
the same as giving the decomposition $x=t\frac p\lambda +\tau \frac{\lambda
-p}\lambda $, we will say that $x$ is of type $0$ (resp. $\infty $) if $\tau
=0$ (resp. $t=0$).

\noindent The equivariant integral is given by 
\begin{equation*}
\int_{X\times {\Bbb P}^1}^{S^1}x=\frac{\int_Xt-\int_X\tau }\lambda \text{, } 
\end{equation*}
and the equivariant pairing is written as follows: 
\begin{equation*}
\left\langle x,x^{\prime }\right\rangle _{X\times {\Bbb P}^1}^G=\frac{%
\left\langle t,t^{\prime }\right\rangle _X-\left\langle \tau ,\tau ^{\prime
}\right\rangle _X}\lambda \text{.} 
\end{equation*}

\bigskip
\bigskip

\section{Equivariant Gromow Witten theory}

Keeping the notations of the previous sections, a $T$ action on $X$ induces a $T$
action on the space of maps of curves in $X$, as it acts on the image of
every map. Since the action is compatible with the equivalence relation on
maps, we define an induced action on $\overline{{\cal M}}_{0,n}(X,\beta )$: 
\begin{equation*}
\begin{array}{ccccc}
T & \times & \overline{{\cal M}}_{0,n}(X,\beta ) & \rightarrow & \overline{%
{\cal M}}_{0,n}(X,\beta ) \\ 
(t & , & \left[ C,x_1,...,x_n,\psi \right] ) &  & \left[
C,x_1,...,x_n,t\circ \psi \right]
\end{array}
\text{.} 
\end{equation*}

\noindent The evaluation maps are equivariant, and therefore induce morphisms in
equivariant cohomolgy: 
\begin{equation*}
\rho _i^{*}:H_T^{*}(X)\rightarrow H_T^{*}(\overline{{\cal M}}_{0,n}(X,\beta
))\text{.} 
\end{equation*}
Fix a basis $\left\{ 1=T_0,...,T_m\right\} $ of $H_T^{*}(X)$ as $%
H_T^{*}(\left\{ pt\right\} )-$ module, and take $n\,$classes $\gamma
_1,...,\gamma _n$ in it, and a class $\beta \in H_2(X)$.

\begin{definition}
The equivariant Gromov Witten invariants are: 
\[
I_\beta ^T(\gamma _1,...,\gamma _n):=\int_{\overline{{\cal M}}_{0,n}}^T\rho
_1^{*}(\gamma _1)\cup ...\cup \rho _n^{*}(\gamma _n)\text{.}
\]
\end{definition}

\noindent They satisfy the three properties listed in \ref{smallq}, with the usual
integral replaced by the equivariant one.

Now it is possible to repeat all constructions of the non equivariant case,
i.e. define the potential 
\begin{equation*}
{\cal F}=\sum_{n\geq 3}\frac 1{n!}\sum_{\beta \in H_2(X,{\Bbb Z})}I_\beta
^T\left( \underset{n}{\underbrace{\gamma ,...,\gamma }}\right) \text{,} 
\end{equation*}
evaluate its third derivatives 
\begin{equation*}
{\cal F}_{ijk}=\frac{\partial ^3{\cal F}}{\partial t_i\partial t_j\partial
t_k}=\sum_{n\geq 0}\frac 1{n!}\sum_{\beta \in H_2(X,{\Bbb Z})}I_\beta
^T\left( T_i,T_j,T_k,\underset{n}{\underbrace{\gamma ,...,\gamma }}\right) 
\text{,} 
\end{equation*}
consider the vector space $QH_T^{*}(X)=H_T^{*}(X)\otimes {\Bbb Q}\left[
\left[ t_0,...,t_m\right] \right] $, and define on it the equivariant
quantum product 
\begin{equation*}
T_i*T_j=\sum_{k,l=0}^m{\cal F}_{ijk}g^{kl}T_l\text{,} 
\end{equation*}
where $g^{ij}$ is still the inverse of the intersection matrix with respect
to the chosen basis of $H_T^{*}(X)$.

\begin{proposition}
$\left( QH_T^{*}(X),*\right) $ is a commutative, associative algebra with
unit $T_0$.
\end{proposition}

\begin{example}
Compute $QH_{S^1}^{*}({\Bbb P}^2)$, with the action: 
\[
\begin{array}{ccccc}
S^1 & \times  & {\Bbb P}^2 & \rightarrow  & {\Bbb P}^2 \\ 
t & , & \left[ x_0:x_1:x_2\right]  &  & \left[ t^ax_0:t^bx_1:t^cx_2\right] 
\end{array}
\text{.}
\]
\end{example}

\noindent As a first step, we compute the equivariant cohomology $H_{S^1}^{*}({\Bbb P}%
^2)$; the bundle ${\Bbb P}^2\times _{s^1}{\Bbb P}^\infty \rightarrow {\Bbb P}%
^\infty $ is exactly the bundle ${\Bbb P}({\cal O}\left( -a\right) \oplus 
{\cal O}\left( -b\right) \oplus {\cal O}\left( -c\right) )\rightarrow 
{\Bbb P}^\infty $; following \ref{pnsk}, 
\begin{equation*}
H_{S^1}^{*}({\Bbb P}^2)\cong \frac{{\Bbb Q}\left[ \lambda ,p\right] }{%
(p-a\lambda )(p-b\lambda )(p-c\lambda )}=\frac{{\Bbb Q}\left[ \lambda
,p\right] }{\left( p^3-e_1\lambda p^2+e_2\lambda ^2p-e_3\lambda ^3\right) }, 
\end{equation*}
where we have set: $e_1=a+b+c$, $e_2=ab+ac+bc$, $e_3=abc$.

\noindent A basis for $H_{S^1}^{*}({\Bbb P}^2)$ as a ${\Bbb Q}\left[ \lambda \right] $
- module is $\left\{ T_0,T_1,T_2\right\} $, with $T_i=p^i$. Let us compute $%
\int_{{\Bbb P}^2}^{S^1}T_i\cup T_j=\int_{{\Bbb P}^2}^{S^1}p^{i+j}$ .

\noindent Since $i+j\leq 4$, it is sufficient to compute: 
\begin{eqnarray*}
\int_{{\Bbb P}^2}^{S^1}p^0 &=&\int_{{\Bbb P}^2}^{S^1}p^1=0, \\
\int_{{\Bbb P}^2}^{S^1}p^2 &=&1, \\
\int_{{\Bbb P}^2}^{S^1}p^3 &=&\int_{{\Bbb P}^2}^{S^1}e_1\lambda
p^2-e_2\lambda ^2p+e_3\lambda ^3=e_1\lambda , \\
\int_{{\Bbb P}^2}^{S^1}p^4 &=&\int_{{\Bbb P}^2}^{S^1}\left( e_1^2-e_2\right)
\lambda ^2p^2-\left( e_1e_2-e_3\right) \lambda ^3p+e_1^2\lambda ^4=\left(
e_1^2-e_2\right) \lambda ^2,
\end{eqnarray*}
in order to write down the intersection matrix and its inverse: 
\begin{equation*}
g_{ij}=\left( 
\begin{array}{ccc}
0 & 0 & 1 \\ 
0 & 1 & e_1\lambda \\ 
1 & e_1\lambda & \left( e_1^2-e_2\right) \lambda ^2
\end{array}
\right) \text{, }g^{ij}=\left( 
\begin{array}{ccc}
e_2\lambda ^2 & -e_1\lambda & 1 \\ 
-e_1\lambda & 1 & 0 \\ 
1 & 0 & 0
\end{array}
\right) \text{. } 
\end{equation*}
Note that $\dim \overline{{\cal M}}_{0,k}\left( {\Bbb P}^2,d\right)
=2+\int_dc_1(T{\Bbb P}^2)+k-3=3d+k-1$, and, by degree computation and
properties of GW invariants, using the notation $%
I_d(T_0^{n_0},T_1^{n_1},T_2^{n_2})$ for equivariant invariants, we see that:

\begin{enumerate}
\item  if $n_0>0$, the GW invariants that do not vanish are just $%
I_0(T_0\cup T_1\cup T_1)$ and $I_0(T_0\cup T_0\cup T_2)$, and they are both equal
to $1$;

\item  if $n_0=0$, $I_d(T_1^{n_1},T_2^{n_2})=\left( \int_d^TT_1\right)
^{n_1}I_d(T_2^{n_2})=d^{n_1}I_d(T_2^{n_2})$; since equivariant integral on
moduli space is a push forward along a fiber of complex dimension $\geq 3d-1$%
, $I_d(T_2^{n_2})=0$ for $n_2<3d-1$.
\end{enumerate}

\noindent The potential is: 
\begin{eqnarray*}
{\cal F} &=&\sum_{n_0+n_1+n_2\geq 3}\frac{y_0^{n_0}}{n_0!}\frac{y_1^{n_1}}{%
n_1!}\frac{y_2^{n_2}}{n_2!}\sum_{d\geq 0}I_d\left(
T_0^{n_0},T_1^{n_1},T_2^{n_2}\right) = \\
&=&{\cal F}_{cl}+\sum_{d>0}\sum_{n_1\geq 0}\sum_{k\geq 0}\frac{%
d^{n_1}y_1^{n_1}}{n_1!}\frac{y_2^{3d-1+k}}{\left( 3d-1+k\right) !}I_d\left(
T_2^{3d-1+k}\right) = \\
&=&{\cal F}_{cl}+\sum_{d>0}\exp (dy_1)\sum_{k\geq 0}\frac{y_2^{3d-1+k}}{%
\left( 3d-1+k\right) !}I_d\left( T_2^{3d-1+k}\right) \text{,}
\end{eqnarray*}

\noindent Notice that, since $\dim _{{\Bbb R}}\overline{{\cal M}}_{0,3d-k+1}\left( 
{\Bbb P}^2,d\right) =2(3d-1+3d-1+k)=2(6d-2+k)$, and $T_2^{3d-1+k}$ is a form
of degree $4(3d-1+k)$, the push forward $I_d\left( T_2^{3d-1+k}\right) $ has
degree $2k=12d-4+4k-12d+4-2k$. This implies, looking at $H_{s^1}^{*}\left(
\left\{ pt\right\} \right) $, that $I_d\left( T_2^{3d-1+k}\right)
=N_{d,k}\lambda ^k$. 

Working out calculations, as in \cite{Al} , the associativity equation yields: 
\begin{eqnarray*}
{\cal F}_{222}+{\cal F}_{111}{\cal F}_{122}-{\cal F}_{112}^2 &=& - \left(
e_2+e_1^2\right) \lambda ^2{\cal F}_{112}+2e_1\lambda {\cal F}_{122}+ \left( e_{1} e_{2} - e_{3} \right)
 \lambda ^{3} {\cal F}_{111}
\text{.}
\end{eqnarray*}
From this we get recursive relations for $N_{d,k}^{\prime s}$: 
\begin{eqnarray*}
N_{d,k} &=&\sum\begin{Sb} d_1+d_2=d\geq 1  \\ k_1+k_2=k\geq 0  \end{Sb} %
N_{d_1,k_1}N_{d_1,k_1}d_1d_2\left[ d_1d_2\binom{3d-4+k}{3d_1-2+k_1}-d_1^2%
\binom{3d-4+k}{3d_1-1+k_1}\right] + \\
&&+2e_1dN_{d,k-1}-(e_1^2+e_2)d^2N_{d,k-2}+ \left( e_{1} e_{2} - e_{3} \right)
 d^{3} N_{d,k-3} \text{.}
\end{eqnarray*}

This formula is a deformation of Kontsevich's one in the
non-equivariant setting (see \cite{FP}).

\bigskip\ 

\bigskip\ 

\newpage

\part{Numbers counting curves and differential equations}

\section{Connections}

Let $V$ be the vector subspace of $H^{*}\left( X\right) $ defined in \ref
{smallq}. $V$ is equipped with the non-degenerate bilinear form of
intersection pairing; we choose an orthonormal basis $\left\{
T_0,...,T_m\right\} $ and coordinates $\left\{ t_0,...,t_m\right\} $, and
recall that we have defined the potential $\phi $ and the quantum product on 
$QH^{*}\left( X\right) =V\otimes {\Bbb Q}\left[ \left[ t_0,...,t_m\right]
\right] $, namely $T_i*T_j=\sum_{k=0}^m\phi _{ijk}T_k$.

\noindent The tangent bundle $TV$ is obviously trivial, and we can take $\left\{
T_0,...,T_m\right\} $ as a global frame; we give the formal definition of a
family of connections on $TV$:

\begin{definition}
\[
\nabla _\hbar :=d-\frac 1\hbar \sum_{i=0}^mdt_iT_i*\label{nab1}
\]
\end{definition}

\begin{remark}
Since $\nabla _\hbar \left( T_j\right) =-\frac 1\hbar \sum_{i,k=0}^m\phi
_{ikj}dt_iT_k$, the structure constants of the connections are exactly the
third derivatives of the potential, divided by $\hbar $.

\end{remark}

\begin{proposition}
$\nabla _\hbar $ is flat.
\end{proposition}

{\bf Proof. }Flatness of the connection is equivalent to associativity of
the quantum product. In fact, the matrix of $1$-forms associated with the
connection is $\Omega _{jk}^\hbar =\frac 1\hbar \sum_{i=0}^m\phi
_{ikj}dt_i=\frac 1\hbar d\left( \phi _{kj}\right) $; the matrix of the
curvature $\nabla _\hbar ^2$ is $K^\hbar =d\Omega ^\hbar -\Omega ^\hbar \cup
\Omega ^\hbar =-\Omega ^\hbar \cup \Omega ^\hbar $; since 
\begin{eqnarray*}
K_{jk}^\hbar &=&-\frac 1{\hbar ^2}\sum_{i=0}^md\left( \phi _{ji}\right) \cup
d\left( \phi _{ik}\right) =-\frac 1{\hbar ^2}\sum_{i=0}^m\left( \left(
\sum_{p=0}^m\phi _{pij}dt_p\right) \cup \left( \sum_{q=0}^m\phi
_{qik}dt_q\right) \right) =   \\
\ &=&-\frac 1{\hbar ^2}\sum_{i=0}^m\left( \left( \sum_{p=0}^m\phi
_{pij}dt_p\right) \cup \left( \sum_{q=0}^m\phi _{qik}dt_q\right) \right) =
 \\
\ &=&-\frac 1{\hbar ^2}\sum_{p,q=0}^m\left( \sum_{i=0}^m\phi _{pij}\phi
_{qik}\right) dt_p\cup dt_q= \\
\ &=&\sum\begin{Sb} p,q=0  \\ p<q  \end{Sb} ^m\left( -\frac 1{\hbar
^2}\sum_{i=0}^m\left( \phi _{pij}\phi _{qik}-\phi _{qij}\phi _{pik}\right)
\right) dt_p\cup dt_q\text{,}
\end{eqnarray*}
then 
\begin{equation*}
K^\hbar =0\Longleftrightarrow \sum_{i=0}^m\left( \phi _{pij}\phi _{qik}-\phi
_{qij}\phi _{pik}\right) \text{ }\forall p,q,j,k\text{.}\label{curv0} 
\end{equation*}

\noindent On the other side, 
\begin{eqnarray*}
\left( T_j*T_p\right) *T_k-T_j*\left( T_p*T_k\right) &=&\sum_i\phi
_{jpi}\left( T_i*T_k\right) -\sum_i\phi _{kpi}\left( T_j*T_i\right) =
 \\
\ &=&\sum_{q=0}^m\left( \sum_{i=0}^m\left( \phi _{jpi}\phi _{ikq}-\phi
_{kpi}\phi _{jiq}\right) \right) T_q\text{,}  
\end{eqnarray*}
and the quantum product is associative if and only if $\sum_{i=0}^m\left(
\phi _{pij}\phi _{qik}-\phi _{qij}\phi _{pik}\right) $ $\forall p,q,j,k$=0.

$\square $

\begin{corollary}
For every $\hbar $, the differential equation $\nabla _\hbar s(t,\hbar )$ $=0
$ formally admits $m+1$ linearly independent solutions.
\end{corollary}

\ 

A remarkable construction of the solutions will be given soon using
equivariant quantum cohomology. However, let us show first how we can
recover the coefficients of the differential equation, i.e the derivatives of the
potential, once we have a basis of independent solutions.

\noindent A vector field $s(t,\hbar )=\sum_{i=0}^ms_i(t,\hbar )T_i$ is a solution if
and only if it satisfies: 
\begin{eqnarray*}
\hbar \frac{\partial s_j}{\partial t_k}(t,\hbar ) &=&\sum_{i=0}^ms_i\left(
t,\hbar \right) \phi _{ijk}(t)  \label{nabcord} \\
\forall j,k &=&0,...,m\text{.}
\end{eqnarray*}

\noindent Conversely, given a basis $\left\{ s_\beta (t,\hbar )=\sum_{i=0}^ms_{\beta
i}(t,\hbar )T_i\right\} _{\beta =0,...,m}$ of independent solutions, then
the $\phi _{ijk\text{ }}(t)$ can be easily recovered with linear algebra
computations.

\smallskip\ 

Every time we deal with differential equations, we skip problems of
existence, that is, our solutions are formal solutions; typically the
existence in a formal neighborhood of the initial conditions is ensured. For
details about the analytical point of view, see \cite{Du}.

\bigskip\ 

\bigskip\ 

\section{ Equivariant potential and solutions of $\nabla _\hbar s=0\label
{solnab}$}

We consider $Y:=X\times {\Bbb P}^1$ with the $S^1$ action trivial on $X$ as
in \ref{xp1s1}: 
\begin{equation*}
\begin{array}{ccccc}
S^1 & \times & X\times {\Bbb P}^1 & \rightarrow & X\times {\Bbb P}^1 \\ 
(t & , & \left( x,\left[ x_0,x_1\right] \right) ) &  & \left( x,\left[
t^2x_0,x_1\right] \right)
\end{array}
\text{.} \label{acxp1} 
\end{equation*}

\noindent From now on, in order to simplify calculations, we shall often consider
cohomology with complex coefficients instead of rational ones; nevertheless,
we keep the same notation, since it will be clear which field are we
working on.

\noindent Recall that $H_T^{*}(Y){\Bbb \cong }H^{*}(X){\Bbb \otimes }\frac{{\Bbb C}%
\left[ \hbar ,p\right] }{\left( p^2-p\hbar \right) }$, and that, localizing
at fixed point components, we get: 
\begin{equation*}
H^{*}(X)\otimes \frac{{\Bbb C}\left( \hbar \right) \left[ p\right] }{%
\left( p^2-p\hbar \right) }\cong \left( H^{*}(X)\otimes {\Bbb C}\left(
\hbar \right) \right) {\Bbb \oplus }\left( H^{*}(X){\Bbb \otimes C}%
\left( \hbar \right) \right) \text{;} \label{comxp1} 
\end{equation*}
moreover, if $x\in H_T^{*}(Y{\Bbb )}$ and $\delta (x)=(t,\tau )$, the
equivariant integral is 
\begin{equation*}
\int_Yx=\frac{\int_Xt-\int_X\tau }\hbar \text{; } \label{inteq} 
\end{equation*}
finally, we denote $e_0:=\frac p\hbar $ , $e_\infty :=\frac{\hbar -p}\hbar $.

Now, let us fix a basis $\left\{ \beta _1,...,\beta _k\right\} $ of $H_2(X)$
and let $\beta _0$ be the generator of $H_2({\Bbb P}^1)$, i.e. the dual of
the hyperplane class; let $\left\{ d=\left( d_1,...,d_k\right) ,d_0\right\} $
be coordinates in $H_2(X\times {\Bbb P}^1)$. Insert ``counters '' $\left\{
q=\left( q_1,...,q_k\right) ,q_0\right\} $ in the equivariant potential,
writing it as 
\begin{equation*}
{\cal F}\left( x\right) =\sum_{n\geq 3}\frac 1{n!}\sum_{\left( d,d_0\right)
}q_0^{d_0}q^dI_{\left( d,d_0\right) }^T\left( \underset{n}{\underbrace{%
x,...,x}}\right) \text{,} \label{poteq} 
\end{equation*}
and develop ${\cal F}=\sum_{i=0}^\infty q_0^i{\cal F}^i$.

\noindent Our goal is to prove that {\em solutions for the equation }$\nabla _\hbar
s=0 ${\em \ are given by second order derivatives of }${\cal F}^1$.

\bigskip\ 

\subsection{Finding solutions of $\nabla _\hbar s=0$}

Let us fix some notation: $Y_{k,\left( d,d_0\right) }$ is the ``virtual ''
moduli space $\overline{{\cal M}}_{0,k}(X\times {\Bbb P}^1,(d,d_0))$: by
this we mean that 
\begin{equation*}
\int_{Y_{k,\left( d,d_0\right) }}^T\omega :=\int_{\overline{{\cal M}}_{0,k}(%
{\Bbb P}^n\times {\Bbb P}^1,(d,d_0))}^T\omega \cup {\cal E}_{k,(d,d_0)}^T%
\text{,} \label{gwx} 
\end{equation*}
where ${\cal E}_{k,(d,d_0)}^T=Euler^T(W_{k,(d,d_0)})$, and 
\begin{equation*}
W_{k,(d,d_0)}\left( \left[ C,x_1,...,x_k,(\psi _1,\psi _2)\right] \right)
:=H^0\left( C,\psi _1^{*}\left( \oplus_{i=1}^r{\cal O}(l_i)\right)
\right) \text{.} \label{wkdd0} 
\end{equation*}

\noindent Let $F$ be the non-equivariant potential for $X$; we need to prove a
technical lemma:

\begin{lemma}
\begin{eqnarray}
{\cal F}^0\left( x\right) =\frac{F\left( t\right) -F\left( \tau \right) }%
\hbar \text{,}  \label{f0} \\
\partial _{e_0}{\cal F}^1\left( x\right) =-\partial _{e_\infty }{\cal F}%
^1\left( x\right) =\frac 1\hbar \partial _p{\cal F}^1\left( x\right) =\frac
1\hbar {\cal F}^1\left( x\right) \text{.}  \label{f1}
\end{eqnarray}
\end{lemma}

{\bf Proof.} Since ${\cal F}^0\left( x\right) =\sum_{k\geq 3}\frac
1{k!}\sum_dq^dI_{\left( d,0\right) }^T\left( \underset{k}{\underbrace{x,...,x%
}}\right) $, we have to compute 
\begin{eqnarray*}
I_{\left( d,0\right) }^T\left( \underset{k}{\underbrace{x,...,x}}\right)
&=&\int_{Y_{k,\left( d,d_0\right) }}^T\rho _1^{*}\left( x\right) \cup
...\cup \rho _k^{*}(x)= \\
\ &=&\int_{\overline{{\cal M}}_{0,k}({\Bbb P}^n\times {\Bbb P}%
^1,(d,0))}^T\rho _1^{*}\left( x\right) \cup ...\cup \rho _k^{*}(x)\cup {\cal %
E}_{k,(d,0)}^T\text{.}
\end{eqnarray*}

\noindent By definition, for every $\left[ C,x_1,...,x_k,(\psi _1,\psi _2)\right] \in 
\overline{{\cal M}}_{0,k}({\Bbb P}^n\times {\Bbb P}^1,(d,0))$, the map $\psi
_2$ is constant, and therefore $\overline{{\cal M}}_{0,k}({\Bbb P}^n\times 
{\Bbb P}^1,(d,0))$ can be identified in a $T$ equivariant way to $\overline{%
{\cal M}}_{0,k}({\Bbb P}^n,d)\times {\Bbb P}^1$ ; the connected components
of fixed points are $\overline{{\cal M}}_{0,k}({\Bbb P}^n,d)\times \left\{
0\right\} $ and $\overline{{\cal M}}_{0,k}({\Bbb P}^n,d)\times \left\{
\infty \right\} $.

\noindent Thus we can compute the GW invariants $I_{\left( d,0\right) }^T\left(
x^k\right) $ by the localization theorem.

\noindent By definition, $W_{k,(d,0)}$, restricted to the two copies of $\overline{%
{\cal M}}_{0,k}({\Bbb P}^n,d)$, is $W_{k,d}$ itself. If we decompose $%
x=t+\tau $, and observe that $\rho _i^{*}\left( t\right) $ restricts to $0$
on $\overline{{\cal M}}_{0,k}({\Bbb P}^n,d)\times \left\{ \infty \right\} $
and $\rho _i^{*}\left( \tau \right) $ restricts to $0$ on $\overline{{\cal M}%
}_{0,k}({\Bbb P}^n,d)\times \left\{ 0\right\} $, we get 
\begin{equation*}
I_{\left( d,0\right) }^T\left( \underset{k}{\underbrace{x,...,x}}\right) = 
\end{equation*}
\begin{eqnarray*}
&=&\frac{\int_{\overline{{\cal M}}_{0,k}({\Bbb P}^n,d)}\rho _1^{*}\left(
t\right) \cup ...\cup \rho _k^{*}(t)\cup {\cal E}_{k,d}-\int_{\overline{%
{\cal M}}_{0,k}({\Bbb P}^n,d)}\rho _1^{*}\left( \tau \right) \cup ...\cup
\rho _k^{*}(\tau )\cup {\cal E}_{k,d}}\hbar = \\
&=&\frac{I_d\left( t^k\right) -I_d\left( \tau ^k\right) }\hbar \text{.}
\end{eqnarray*}
From this the first formula follows .

Now write down ${\cal F}^1\left( x\right) =\sum_{k\geq 3}\frac
1{k!}\sum_dq^dI_{\left( d,1\right) }^T\left( \underset{k}{\underbrace{x,...,x%
}}\right) $, and observe that $\partial _1{\cal F}^1\left( x\right)
=\sum_{k\geq 2}\frac 1{k!}\sum_dq^dI_{\left( d,1\right) }^T\left( \underset{k%
}{\underbrace{x,...,x}},1\right) $ $=0$ by properties of GW invariants,
since $(d,1)\neq 0$.

\noindent We can see that 
\begin{eqnarray*}
\frac 1\hbar \partial _p{\cal F}^1\left( x\right) &=&\frac 1\hbar
\sum_k\frac 1{k!}\sum_dq^dI_{\left( d,1\right) }^T\left( \underset{k}{%
\underbrace{x,...,x}},p\right) = \\
&=&\frac 1\hbar \sum_k\frac 1{k!}\sum_dq^d\int p\cdot I_{\left( d,1\right)
}^T\left( \underset{k}{\underbrace{x,...,x}}\right) = \\
&=&\frac 1\hbar \sum_k\frac 1{k!}\sum_dq^dI_{\left( d,1\right) }^T\left( 
\underset{k}{\underbrace{x,...,x}}\right) =\frac 1\hbar {\cal F}^1\left(
x\right) \text{.}
\end{eqnarray*}

\noindent Furthermore, since $e_0=\frac p\hbar $, we have 
\begin{eqnarray*}
\frac 1\hbar \partial _p{\cal F}^1\left( x\right) &=&\frac 1\hbar
\sum_k\frac 1{k!}\sum_dq^dI_{\left( d,1\right) }^T\left( \underset{k}{%
\underbrace{x,...,x}},p\right) = \\
&=&\sum_k\frac 1{k!}\sum_dq^d\frac 1\hbar I_{\left( d,1\right) }^T\left( 
\underset{k}{\underbrace{x,...,x}},p\right) = \\
&=&\sum_k\frac 1{k!}\sum_dq^dI_{\left( d,1\right) }^T\left( \underset{k}{%
\underbrace{x,...,x}},\frac p\hbar \right) =\partial _{e_0}{\cal F}^1\left(
x\right) \text{;}
\end{eqnarray*}

\noindent By definition $e_0+e_\infty =1$, and therefore by a similar computation we
get: $0=\partial _1{\cal F}^1\left( x\right) =\partial _{e_0}{\cal F}%
^1\left( x\right) +\partial _{e_\infty }{\cal F}^1\left( x\right)
\Rightarrow \partial _{e_0}{\cal F}^1\left( x\right) =-\partial _{e_\infty }%
{\cal F}^1\left( x\right) $ .
\begin{flushright}
$\Box$
\end{flushright}

\smallskip\ 

Let us take an orthonormal basis $\left\{ T_0,...,T_m\right\} $ of $%
V\subseteq H^{*}(X)$; by the localization isomorphism, an orthogonal basis
for the equivariant pairing in $H_T^{*}(X)$ is 
\begin{equation*}
\left\{ \frac p\hbar T_0,...,\frac p\hbar T_m,\frac{\hbar -p}\hbar T_0,...,%
\frac{\hbar -p}\hbar T_m\right\} \text{;} \label{basis} 
\end{equation*}
the coordinates with respect to the first $m+1$ classes, the so called
classes of type $0,$ are $\left\{ t_i\right\} _{i=0,...,m}$, the coordinates
of classes of type $\infty $ are $\left\{ \tau _i\right\} _{i=0,...,m}$;
each vector in the basis has squared norm $\frac 1\hbar $.

\noindent Therefore, if $x=t+\tau $, this means that $t=\sum_{i=0}^mt_i\left( \frac
p\hbar T_i\right) $ and $\tau =\sum_{i=0}^m\tau _i\left( \frac{\hbar -p}%
\hbar T_i\right) $.

\begin{proposition}
The vector field $\sum_{k=0}^m\frac{\partial ^2{\cal F}^1\left( x\right) }{%
\partial t_k\partial \tau _u}T_k$ is a solution of the differential equation 
$\nabla _\hbar s(t,\hbar )=0$, for every $u$.
\end{proposition}

{\bf Proof.} Let us write down the asociativity identity: 
\begin{equation*}
\sum_k{\cal F}_{xyk}{\cal F}_{kzu}=\sum_k{\cal F}_{xzk}{\cal F}_{kyu}\text{,}
\label{assequiv} 
\end{equation*}
and take the coefficient of $q_0$: 
\begin{equation*}
\sum_k{\cal F}_{xyk}^0{\cal F}_{kzu}^1+{\cal F}_{xyk}^1{\cal F}%
_{kzu}^0=\sum_k{\cal F}_{xzk}^0{\cal F}_{kyu}^1+{\cal F}_{xzk}^1{\cal F}%
_{kyu}^0\text{.} \label{assdeg1} 
\end{equation*}
Next specialize $y=e_0=\frac p\hbar T_0,x$ and $z$ of type $0$, $u$ of type $%
\infty $, and notice from formula \ref{f0} that every time we differentiate $%
{\cal F}^0$ once by a variable of type $0$, and once by a variable of type $%
\infty $, we get $0$, therefore ${\cal F}_{kzu}^0={\cal F}_{kyu}^0=0$. Also $%
{\cal F}_{kyu}^1={\cal F}_{kue_0}^1=\frac 1\hbar {\cal F}_{ku}^1$ by \ref{f1}%
, and 
\begin{equation*}
{\cal F}_{xyk}^0=\left\{ 
\begin{array}{c}
0\text{ if }k\text{ is of type }\infty \\ 
\frac{F_{xk1}(t)}\hbar =\frac{\left\langle x,t_k\right\rangle _X}\hbar \text{
if }k\text{ is of type }0
\end{array}
\right. \text{,} 
\end{equation*}
hence, on the left hand side we get $\frac 1\hbar \sum_k{\cal F}%
_{kzu}^1\left\langle x,t_k\right\rangle _X=\frac 1\hbar {\cal F}_{xzu}^1$

\noindent Hence substituting we obtain: 
\begin{equation*}
\frac 1\hbar {\cal F}_{xzu}^1=\frac 1\hbar \text{ }\sum_k\frac{F_{xzk}}\hbar 
{\cal F}_{ku}^1=\text{ }\frac 1{\hbar ^2}\sum_kF_{xzk}{\cal F}_{ku}^1\text{,}
\end{equation*}
for every $x,z$ of type $0$, $u$ of type $\infty $.

\noindent On the other side, $\nabla _\hbar \left( \sum_k{\cal F}_{ku}^1T_k\right) =0$
if and only if 
\begin{equation*}
\hbar \frac{\partial {\cal F}_{xu}^1}{\partial tz}(t,\hbar )=\sum_{k=0}^m%
{\cal F}_{ku}^1\left( t,\hbar \right) F_{kzx}(t). 
\end{equation*}
\begin{flushright}
$\Box$
\end{flushright}

\bigskip\ 

\subsection{Separating variables} \label{separate}

We want to transform ${\cal F}_{ku}^1$ by mean of the localization formula
in order to ``separate'' the dependence on the classes $t$ and $\tau $.

We have to understand, in the ``virtual '' space $Y_{k,\left( d,1\right) }$,
the nature of the $S^1$-fixed points, that is, we have to understand the $S^1$
fixed points in $\overline{{\cal M}}_{0,k}({\Bbb P}^n\times {\Bbb P}%
^1,(d,1)) $, and the behaviour of the Euler class.

\noindent From now on, assume $X={\Bbb P}^n$, and then modify GW invariants with
suitable Euler classes.

\begin{figure}[htb]
\hspace{-2cm} 
\begin{center} 
\mbox{\epsfig{file=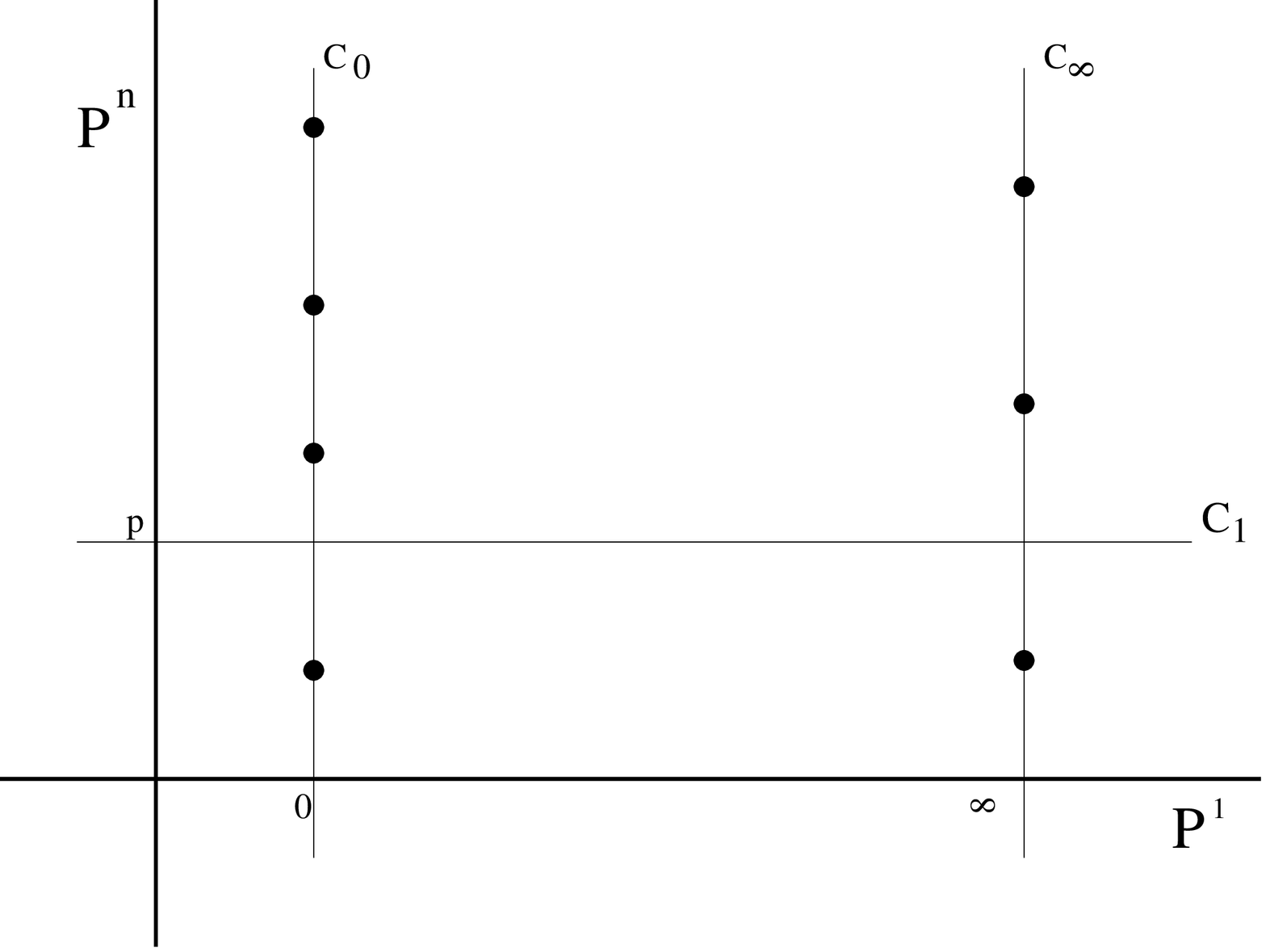,width=8cm,height=10cm,angle=270}}
\caption[]{
\label{dis1}}
\end{center}
\end{figure}

\noindent Remind that the fixed points have the following description (see figure \ref
{dis1}): 
\begin{equation*}
Y_{k,\left( d,1\right) }=\left\{ 
\begin{array}{c}
\left[ C,x_1,...,x_k,(\psi _1,\psi _2)\right] : \\ 
\left[ C,x_1,...,x_k,(\psi _1,\psi _2)\right] =\left[ C,x_1,...,x_k,(\psi
_1,t\psi _2)\right] \forall t\in S^1
\end{array}
\right\} \text{.}\label{ipst} 
\end{equation*}

\noindent This implies immediately that all the special points of $C$ are mapped by $%
\psi _2$ on the fixed points of ${\Bbb P}^1$. In particular we have $k_0$
marked points mapping to $0$ and $k_\infty $ marked points mapping to $%
\infty $, $k=k_0+k_\infty $.

\noindent Furthermore, by definition $\psi _2$ has degree $1$, hence in general there
is one special irreducible component $C_1$ of $C$ which maps isomorphically
to ${\Bbb P}^1$, while the image under $\psi _2$ of each one of the other
components is a fixed point in ${\Bbb P}^1$. Thus the curve $C$ is the union
of three irreducible components $C_1,C_0$ and $C_\infty $, such that: 
\begin{eqnarray*}
\psi _{2\mid C_1} &:&C_1\overset{\sim }{\rightarrow }{\Bbb P}^1 \\
\psi _2\left( C_0\right) &=&0 \\
\psi _2\left( C_\infty \right) &=&\infty \text{.}
\end{eqnarray*}

\noindent We claim that $\psi _{1\mid C_1}$ is constant: in fact, if $p\in C_1$ is
generic and $\psi \left( p\right) =(x,a)$, then we obtain $\left(
x,ta\right) $ by acting with the group. For every $t\in S^1$, there must be a
lifting of this action to $C_1$ which maps $p$ to some other generic point
of $C_1$ in $\psi _1^{-1}\left( x\right) $, hence $\psi _1^{-1}\left(
x\right) $ contains an entire $S^1$-orbit, and therefore the whole $C_1$. It
follows that $d=\deg \psi _1=\deg \psi _{1\mid C_1}+\deg \psi _{1\mid
C_0} + \deg \psi _{1\mid C_\infty }=0+d_0+d_\infty $.

\noindent Finally remark that the marked points are necessarily on the two components
(possibly degenerating to a point) $C_0$ and $C_\infty $.

\begin{definition}
$$
M_{d_0,d_\infty ,k_0,k_\infty }:=\left\{ 
\begin{array}{c}
\left[ C,x_1,...,x_k,(\psi _1,\psi _2)\right] \in \left( Y_{k,\left(
d,1\right) }\right) ^T: \\ 
\begin{array}{c}
\psi _{2\mid C_1}:C_1\stackrel{\symbol{126}}{\longrightarrow }{\Bbb P}%
^1,\psi _2\left( C_0\right) =0,\psi _2\left( C_\infty \right) =\infty , \\ 
\deg \psi _{1\mid C_0}=d_0,\deg \psi _{1\mid C_\infty }=d_\infty , \\ 
\#\left\{ \left\{ x_1,...,x_k\right\} \cap C_0\right\} =k_0,\#\left\{
\left\{ x_1,...,x_k\right\} \cap C_\infty \right\} =k_\infty \text{.}
\end{array}
\end{array}
\right\}
$$
\end{definition}

We just proved:

\begin{proposition}
\label{fixs1} 
\[
\left( Y_{k,\left( d,1\right) }\right) ^T=\cup \Sb d_0+d_\infty =d \\ %
k_0+k_\infty =k \endSb M_{d_0,d_\infty ,k_0,k_\infty }\text{.}
\]
\end{proposition}

\noindent In order to apply the localization formula we have now to compute:

\begin{itemize}
\item  the Euler classes of the normal bundle of the components in the $%
M_{d_0,d_\infty ,k_0,k_\infty }$'s,

\item  the restriction of the vector bundle $W_{k,\left( d,1\right) }$.
\end{itemize}

We denote by $x_0$ (resp. $x_\infty $) the point $C_0\cap C_1$ (resp. $%
C_\infty \cap C_1$); we shall say that $C_0$ (resp. $C_\infty $) is
degenerate if $\psi \left( x_0\right) $ (resp.$\psi \left( x_\infty \right) $%
) is smooth in the image of the curve; moreover, in this case $d_0=0$ (resp. 
$d_\infty =0$).

\noindent We have to distinguish four cases:

\begin{enumerate}
\item  \label{caso1}Both $C_0$ and $C_\infty $ are degenerate : this implies
that $d=0$, $k\leq 2$.

\item  Only $C_0$ is degenerate; we have two cases according if $k_0=0$, $%
k_0=1$.

\begin{itemize}
\item  $k_0=0$ ; the variety $M_{0,d,0,k}$ is isomorphic to $X_{k+1,d}$ : 
\begin{equation*}
\begin{array}{ccc}
M_{0,d,0,k} & \rightarrow  & X_{k+1,d} \\ 
\left[ C,x_1,...,x_k,(\psi _1,\psi _2)\right]  &  & \left[ C_\infty
,x_\infty ,x_1,...,x_k,\psi _1\right] 
\end{array}
\text{.}\label{iso1}
\end{equation*}

\item  \label{caso2}$k_0=1$ ; hence $x_0=x_i$ for some $i\in \left\{
1,...,k\right\} $, and the variety $M_{0,d,1,k-1}$ is a union of components
isomorphic to $X_{k,d}$: 
\begin{equation*}
\begin{array}{ccc}
M_{0,d,1,k-1}^i & \rightarrow  & X_{k,d} \\ 
\left[ C,x_1,...,x_k,(\psi _1,\psi _2)\right]  &  & \left[ C_\infty
,x_\infty ,x_1,...,x_{i-1},x_{i+1},...,x_k,\psi _1\right] 
\end{array}
\text{.}\label{iso2}
\end{equation*}
\end{itemize}

\item  Similarly for $C_\infty $ degenerate:

\begin{itemize}
\item  $M_{d,0,k,0}\cong X_{k+1,d}$

\item  \label{caso3}$M_{d,0,k-1,1}^i\cong X_{k,d}$.
\end{itemize}

\item  \label{caso4}$C_0,C_\infty $ are both non degenerate. For every
partition of $\left\{ 1,...,k\right\} $ in two subsets $A,B$ of cardinality $%
k_0,k_\infty $, we get a component $M_{d_0,d_\infty ,k_0,k_\infty }^{A,B}$,
and the following map is an isomorphism: 
\begin{equation*}
\begin{array}{ccc}
M_{d_0,d_\infty ,k_0,k_\infty }^{A,B} & \rightarrow  & X_{k_0+1,d_0}\times
_XX_{k_\infty +1,d_\infty } \\ 
\left[ C,x_1,...,x_k,(\psi _1,\psi _2)\right]  &  & \left[ \left[
C_0,x_0,x_{a_1},...,x_{a_{k_0}},\psi _1\right] ,\left[ C_\infty ,x_\infty
,x_{b_1},...,x_{b_{k_\infty }},\psi _1\right] \right] 
\end{array}
\text{,}\label{iso3}
\end{equation*}
where $X_{k_0+1,d_0}\overset{\rho _0}{\rightarrow }X$ and $X_{k_\infty
+1,d_\infty }\overset{\rho _\infty }{\rightarrow }X$ are the evaluation maps
on $x_0$ and $x_\infty $.
\end{enumerate}

\bigskip\ 

Recall that, for every $Z$ convex, on $Z_{k,d}$ we have line bundles $F_i$ , 
$i=1,...,k$; by definition $F_i\left( \left[ C,x_1,...,x_k,(\psi _1,\psi
_2)\right] \right) =T_{x_i}C$, the tangent space to the curve $C$ at the
point $x_i$. In case of a group action on $Z$, these line bundle are
canonically linearized by the pull back of the action on the tangent space
at the point $\psi \left( x_i\right) $ to the image $\psi \left( C\right) $.

\noindent We shall denote by $c_i$ the equivariant Chern class of $F_i$. All the same,
in the case where $C_0$ and $C_\infty $ are non degenerate, we get $%
F_0\left( \left[ C,x_1,...,x_k,(\psi _1,\psi _2)\right] \right) :=T_{x_0}C_0$%
, and $F_\infty \left( \left[ C,x_1,...,x_k,(\psi _1,\psi _2)\right] \right)
:=T_{x_\infty }C_\infty $, with equivariant Chern classes $c_0$ and $%
c_\infty $.

\noindent We want to compute the functions $\frac{\partial ^2{\cal F}^1\left( x\right) 
}{\partial t_k\partial \tau _u}$ by localization, thus we are only
interested in mixed derivatives with a variable of type $0$ and a variable
of type $\infty $, which vanish in the cases where we have $k_0=0$ or $%
k_\infty =0$; in fact, in these cases, the classes we should integrate to
obtain GW invariants depend only on variables of type $\infty $ (or of type $%
0$).

\bigskip\ 

\begin{remark}
Deformations of stable maps\label{norbu}.
\end{remark}

We need to describe normal bundles to closed subvarieties in the moduli
space of stable maps, hence we need a description of the tangent space to
the moduli space $\overline{{\cal M}}_{0,k}\left( X,\beta \right) $ in terms
of infinitesimal deformation of stable maps. For details and proofs see \cite{FP} . Let $%
\left[ C,x_1,...,x_k,(\mu )\right] $ be a stable map, and let $Def\left( \mu
\right) $ be the space of equivalence classes of first order infinitesimal deformations,
i.e. the tangent space to $\overline{{\cal M}}_{0,k}\left( X,\beta \right) $
at $\left[ C,x_1,...,x_k,(\mu )\right] $. If $\left\{ q_1,...q_d\right\} $
are the nodes on the curve $C$, and $C_i^{\prime },C_i^{\prime \prime }$ are
the two branches of the curve meeting in $q_i$, we can decompose $Def\left(
\mu \right) =Def_G\left( \mu \right) \oplus \left(
\oplus_{i=1}^dT_{q_i}C_i^{\prime }\otimes T_{q_i}C_i^{\prime \prime
}\right) $; by $Def_G\left( \mu \right) $ we mean {\em infinitesimal deformations of the
map which preserve the combinatorial type} (i.e. the dual graph) of the
curve; there is a surjective map $Def_G\left( \mu \right) \rightarrow
Def_G\left( C\right) \rightarrow 0$ on the space of infinitesimal deformations of the
curve which preserve the combinatorial type, and we obtain the exact
sequence 
\begin{equation*}
0\rightarrow Def_C\left( \mu \right) \rightarrow Def_G\left( \mu \right)
\rightarrow Def_G\left( C\right) \rightarrow 0\text{;} 
\end{equation*}
by $Def_C\left( \mu \right) $ we denote the space of{\em infinitesimal deformations of the map
whose restriction to the curve is trivial}. This space plays a role in another
sequence: 
\begin{equation*}
0\rightarrow H^0\left( C,TC\left( -\sum_{i=1}^dq_i\right) \right)
\rightarrow Def_R\left( \mu \right) \rightarrow Def_C\left( \mu \right)
\rightarrow 0\text{;} 
\end{equation*}
here, $Def_R\left( \mu \right) $ is the space of {\em infinitesimal deformations of the map
which keep the curve rigidly fixed}, and the kernel consists of
reparametrization of the curve fixing the nodes. Looking at $Hom\left(
C,X\right) $ as an open subset of the Hilbert scheme of graphs in $X\times C$%
, since $T_{\left[ \mu \right] }$ $Hom\left( C,X\right) \cong H^0\left(
C,\mu ^{*}\left( TX\right) \right) $, we build the last exact sequence in
order to obtain information on $Def_R\left( \mu \right) $: 
\begin{equation*}
0\rightarrow \oplus_{j=1}^kDef_i\left( \mu \right) \rightarrow
Def_R\left( \mu \right) \rightarrow H^0\left( C,\mu ^{*}\left( TX\right)
\right) \rightarrow 0\text{;} 
\end{equation*}
$Def_j\left( \mu \right) $ $\cong T_{x_j}C$ is the infinitesimal deformation which just
moves the $j-$th marking.

\bigskip\ 

\begin{proposition}
Let $D$ be a connected component of fixed points in $Y_{k,\left( d,1\right) }
$. Then 
\[
Euler\left( {\cal N}_{D/Y_{k,\left( d,1\right) }}\right) =
\]

\[
=\left\{ 
\begin{array}{cc}
-\hbar ^2 & \text{if }D\subset M_{0,0,1,1}\text{(\ref{caso1})} \\ 
-\hbar ^2\left( -\hbar +c_\infty \right)  & \text{if }D\subset M_{0,d,1,k-1}%
\text{(\ref{caso2})} \\ 
-\hbar ^2\left( \hbar +c_0\right)  & \text{if }D\subset M_{d,0,k-1,1}\text{(%
\ref{caso3})} \\ 
-\hbar ^2\left( \hbar +c_0\right) \left( -\hbar +c_\infty \right)  & \text{%
if }D\subset M_{d_0,d_\infty ,k_0,k_\infty }\text{, }d_0,d_\infty \neq 0%
\text{ (\ref{caso4})}
\end{array}
\right. \text{.}\label{euler}
\]
\end{proposition}

{\bf Proof.}

{\bf Case 1) }$Y_{2,\left( 0,0\right) }\cong X\times {\Bbb P}^1\times {\Bbb P%
}^1$, $T$ -equivariantly. The isomorphism is 
\begin{equation*}
\left[ C,x_0,x_\infty ,(\psi _1,\psi _2)\right] \rightarrow \left( \psi
_1\left( C\right) ,\psi _2\left( x_0\right) ,\psi _2\left( x_\infty \right)
\right) \text{,} 
\end{equation*}
and the two connected components of fixed points we need to consider, say $D$
and $D^{\prime }$ corresponds to $X\times \left\{ 0\right\} \times \left\{
\infty \right\} $ and $X\times \left\{ \infty \right\} \times \left\{
0\right\} $. The linearized normal bundle at some point $\left\{ x\right\}
\times \left\{ 0\right\} \times \left\{ \infty \right\} $ is identified with
the two dimensional vector space $\left\{ x\right\} \times {\Bbb P}%
^1\backslash \left\{ \infty \right\} \times {\Bbb P}^1\backslash \left\{
0\right\} $, and $S^1$ acts with characters $(-\hbar ,\hbar )$, hence the
equivariant Euler class is $-\hbar ^2$; similarly, $S^1$ acts with
characters $(\hbar ,-\hbar )$ on the normal bundle at some point $\left\{
x\right\} \times \left\{ \infty \right\} \times \left\{ 0\right\} $, and the
Euler class is still $-\hbar ^2$.

{\bf Case 2) }Consider the following subvarieties in $Y_{k,\left( d,1\right)
}$ : 
\begin{equation*}
E_i=\left\{ 
\begin{array}{c}
\left[ C,x_1,...,x_k,(\psi _1,\psi _2)\right] :C=C_1\cup C_2, \\ 
\deg \psi _{2\mid C_1}=1,\deg \psi _{2\mid C_2}=0, \\ 
\left\{ \left\{ x_1,...,x_k\right\} \cap C_1\right\} =\left\{ x_i\right\}
\end{array}
\right\} \text{.} \label{ei} 
\end{equation*}

\noindent The preimage of the divisor $D_{\left\{ i\right\} \left\{
1,...,i-1,i+1,...,k\right\} }$ under the map $Y_{k,\left( d,1\right)
}\rightarrow {\cal M}_{0,k}$, $\left[ C,x_1,...,x_k,(\psi _1,\psi _2)\right]
\rightarrow \left[ C,x_1,...,x_k\right] $, consists of two components,
according to $\deg \psi _{2\mid C_1}=1$ or $\deg \psi _{2\mid C_1}=0$ ; the
first one is exactly $E_i$, hence it is a divisor in $Y_{k,\left( d,1\right)
}$.

\noindent This divisor maps equivariantly to ${\Bbb P}^1\times {\Bbb P}^1$: 
\begin{equation*}
\begin{array}{ccc}
E_i & \rightarrow & {\Bbb P}^1\times {\Bbb P}^1 \\ 
\left[ C,x_1,...,x_k,(\psi _1,\psi _2)\right] &  & \left( \psi _2\left(
x_i\right) ,\psi _2\left( x_\infty \right) \right)
\end{array}
\text{,} 
\end{equation*}
and the connected components of fixed points $M_{0,d,1,k-1}^i$ are the
preimages of the point $\left( 0,\infty \right) $. By the same reasoning of
the previous case, ${\cal N}$ $_{M_{0,d,1,k-1}^i/E_i}=-\hbar ^2$. In order
to obtain ${\cal N}$ $_{M_{0,d,1,k-1}^i/Y_{k,\left( d,1\right) }}$, we have
to multiply by $Euler\left( {\cal N}_{E_i/Y_{k,\left( d,1\right) }}\right) $.

\noindent By the remark \ref{norbu}, the fiber of the bundle ${\cal N}%
_{E_i/Y_{k,\left( d,1\right) }}$is the tensor product of the two tangent
spaces to the two irreducible components of the curve in the join point,
since the combinatorial type of the curve is fixed in $E_i$.

\noindent In our case, the tangent space to $C_1$ at the point $x_i$ is acted on by $S^1$%
with character $-\hbar $, the line bundle on $E_i$ is trivial with
equivariant Chern class ; furthermore the tangent space to $C_\infty $ at $%
x_\infty $ has trivial action, hence the equivariant Chern class of the line
bundle is $c_\infty $. Thus $Euler\left( {\cal N}_{E_i/Y_{k,\left(
d,1\right) }}\right) =-\hbar +c_\infty $, and $Euler\left( {\cal N}%
_{D/Y_{k,\left( d,1\right) }}\right) =-\hbar ^2\left( -\hbar +c_\infty
\right) $.

{\bf Case 3) }The computation is exactly the same as in case 2, except that
the action on the tangent space to $C_1$ at the point $x_i$ has character $%
\hbar $, and the equivariant Chern class of the line bundle of the tangent
space to $C_0$ at $x_0$ is $c_0$, hence $Euler\left( {\cal N}_{D/Y_{k,\left(
d,1\right) }}\right) =-\hbar ^2\left( \hbar +c_0\right) $.

{\bf Case 4)} Again, consider the map $\pi :Y_{k,\left( d,1\right)
}\rightarrow \overline{{\cal M}}_{0,k}$; for every partition of $\left\{
1,...,k\right\} $ in two subsets $A,B$ of cardinality $k_0,k_\infty $, $\pi
^{-1}\left( D_{A,B}\right) $ is a divisor, and contains as codimension $1$
subvariety 
\begin{equation*}
E_{A,B}=\left\{ 
\begin{array}{c}
\left[ C,x_1,...,x_k,(\psi _1,\psi _2)\right] \in \pi ^{-1}\left(
D_{A,B}\right) :C=C_1\cup C_2\cup C_3, \\ 
\deg \psi _{2\mid C_1}=1,\deg \psi _{2\mid C_2}=\deg \psi _{2\mid C_3}=0
\end{array}
\right\} \text{,} \label{eab} 
\end{equation*}
which can be seen as the complete intersection of divisors 
\begin{equation*}
\left\{ 
\begin{array}{c}
\left[ C,x_1,...,x_k,(\psi _1,\psi _2)\right] \in \pi ^{-1}\left(
D_{A,B}\right) \\ 
\deg \psi _{2\mid C_1}=1,\deg \psi _{2\mid C_2}=0
\end{array}
\right\} \cap \left\{ 
\begin{array}{c}
\left[ C,x_1,...,x_k,(\psi _1,\psi _2)\right] \in \pi ^{-1}\left(
D_{A,B}\right) \\ 
\deg \psi _{2\mid C_1}=0,\deg \psi _{2\mid C_2}=1
\end{array}
\right\} \text{.} 
\end{equation*}

\noindent Hence, by the same computations, $Euler\left( {\cal N}_{E_{A,B}/Y_{k,\left(
d,1\right) }}\right) =\left( -\hbar +c_\infty \right) \left( \hbar
+c_0\right) $.

\noindent Analogously, we get an equivariant map 
\begin{equation*}
\begin{array}{ccc}
E_{A,B} & \rightarrow & {\Bbb P}^1\times {\Bbb P}^1 \\ 
\left[ C,x_1,...,x_k,(\psi _1,\psi _2)\right] &  & \left( \psi _2\left(
x_0\right) ,\psi _2\left( x_\infty \right) \right)
\end{array}
\text{,} 
\end{equation*}
and our connected component of fixed points is the preimage of $(0,\infty )$%
, hence has normal bundle in $E_{A,B}$ with equivariant Euler class $-\hbar
^2$. Thus $Euler\left( {\cal N}_{D/Y_{k,\left( d,1\right) }}\right) =-\hbar
^2\left( \hbar +c_0\right) \left( -\hbar +c_\infty \right) $.
\begin{flushright}
$\Box$
\end{flushright}

\bigskip\ 

Recall that by definition we have 
\begin{equation*}
\frac{\partial ^2{\cal F}^1\left( x\right) }{\partial t_i\partial \tau _j}%
=\sum_{k\geq 1}\frac 1{k!}\sum_dq^dI_{\left( d,1\right) }^T\left(
x^k,T_i\frac p\hbar ,T_j\frac{\hbar -p}\hbar \right) \text{,} \label{dersec} 
\end{equation*}
hence, in order to separate variables, we apply the integration formula to 
\begin{equation*}
I_{\left( d,1\right) }^T\left( x^k,T_i\frac p\hbar ,T_j\frac{\hbar -p}\hbar
\right) \text{.} 
\end{equation*}

\noindent First of all, let us fix a suitable basis for $H^{*}\left( X\right) $; let $%
v:=Euler\left( \oplus_{i=1}^r{\cal O}\left( l_i\right) \right) $, and let 
$\left\{ T_0,...,T_s\right\} $ be a set of independent elements such that $%
\left\{ \left[ T_0\right] ,...,\left[ T_s\right] \right\} $ form an
orthonormal basis in $\frac{H^{*}\left( X\right) }{Ann(v)}$, hence $%
\int_XT_i\cup T_jv=\delta _{ij}$. Let us denote $T_iv$ with $T_i^{\prime }$
and complete to orthonormal dual basis of $H^{*}(X)$, $\left\{
T_0,...,T_s,T_{s+1},...,T_m\right\} $ and $\left\{ T_0^{\prime
},...,T_s^{\prime },T_{s+1}^{\prime },...,T_m^{\prime }\right\} $, such
that the class of the diagonal in $X$ is 
\begin{equation*}
\Delta =\sum_{i=0}^mT_i\otimes T_i^{\prime }\text{.} 
\end{equation*}

\noindent We also define $W_{k+1,d}^{\prime }$ as the kernel of the map 
\begin{eqnarray*}
W_{k+1,d} &\rightarrow &\rho _1^{*}\left( \oplus_{i=1}^r{\cal O}%
(l_i)\right)   \\
\sigma &\rightarrow &\sigma \left( \psi \left( x_{k+1}\right) \right) \text{,%
}
\end{eqnarray*}
therefore ${\cal E}_{k+1,d}^{\prime }:=Euler(W_{k+1,d}^{\prime })$ satisfies 
${\cal E}_{k+1,d}={\cal E}_{k+1,d}^{\prime }\rho _{k+1}^{*}\left( v\right) $.

\begin{lemma}
\label{contributo}Let $D=M_{d_0,d_\infty ,n_0,n_\infty }^{A,B}$, then 
\[
\int_D\frac{\rho _1^{*}\left( x\right) ...\rho _k^{*}\left( x\right) {\cal E}%
_{k,\left( d,1\right) }}{Euler\left( {\cal N}_{D/Y_{k,\left( d,1\right)
}}\right) }=\label{contrab}
\]
\[
=\sum_{i=0}^m\int_{X_{k_0+1,d_0}}\rho _1^{*}\left( t\right) ...\rho
_{k_0}^{*}\left( t\right) \frac{\rho _0^{*}\left( T_i\right) {\cal E}%
_{k_0+1,d}}{\hbar (\hbar +c_0)}\cdot
\]
\[
\cdot
\int_{X_{k_\infty +1,d_\infty }}\rho
_1^{*}\left( \tau \right) ...\rho _{k_\infty }^{*}\left( \tau \right) \frac{%
\rho _\infty ^{*}\left( T_i\right) {\cal E}_{k_\infty +1,d}}{\hbar (-\hbar
+c_\infty )}\text{.}
\]
\end{lemma}

{\bf Proof. }We have $\binom k{k_0}$ such components in $M_{k_0,k_\infty
,d_0,d_\infty }$, on varying the sets $A$ and $B$, and each one gives the
same contribution, thus we assume without loss of generality that $A=\left\{
1,...,k_0\right\} $; recall that we have the isomorphism:

$
\begin{array}{ccc}
M_{d_0,d_\infty ,k_0,k_\infty }^{A,B} & \overset{\phi }{\rightarrow } & 
X_{k_0,d_0+1}\times _XX_{k_\infty ,d_\infty +1} \\ 
\left[ C,x_1,...,x_k,(\psi _1,\psi _2)\right] &  & \left[ \left[
C_0,x_0,x_1,...,x_{k_0},\psi _1\right] ,\left[ C_\infty ,x_\infty
,x_{k_0+1},...,x_k,\psi _1\right] \right]
\end{array}
$.

\noindent As a first step, let us compute the restriction of ${\cal E}_{k,d}$ to this
component. Let $\nu $ be the evaluation map on the meeting point, then we
get the exact sequence: 
\begin{equation*}
\begin{array}{ccccccccc}
0 & \rightarrow & W_{k_0+1,d_0}^{\prime }\oplus W_{k_\infty +1,d_\infty
}^{\prime } & \rightarrow & W_{k,\left( d,1\right) } & \rightarrow & \nu
^{*}\left( \oplus_i{\cal O}\left( l_i\right) \right) & \rightarrow & 0 \\ 
&  &  &  & \left[ \psi _1^{*}\left( \sigma \right) ,\psi _1^{*}\left( \tau
\right) \right] & \rightarrow & \sigma \left( x_0\right) =\tau \left(
x_\infty \right) &  & 
\end{array}
\text{.}  
\end{equation*}

\noindent Therefore ${\cal E}_{k,\left( d,1\right) }={\cal E}_{k_0+1,d_0}^{\prime }%
{\cal E}_{k_\infty +1,d_\infty }^{\prime }\nu ^{*}\left( v\right) $.

\noindent Putting everything toghether, and denoting with $\mu =\left( \rho _0,\rho
_\infty \right) $ the product of evaluation maps on the last marked points, $%
X_{k_0,d_0+1}\times X_{k_\infty ,d_\infty +1}\rightarrow X\times X$, we get 
\begin{equation*}
\int_D\frac{\rho _1^{*}\left( x\right) ...\rho _k^{*}\left( x\right) {\cal E}%
_{k,\left( d,1\right) }}{Euler\left( {\cal N}_{D/Y_{k,\left( d,1\right)
}}\right) }= \label{contrib2} 
\end{equation*}
\begin{eqnarray*}
\ &=&\int_{D}\frac{\rho
_1^{*}\left( t\right) ...\rho _{k_0}^{*}\left( t\right) \rho
_{k_0+1}^{*}\left( \tau \right) ...\rho _k^{*}\left( \tau \right) {\cal E}%
_{k_0+1,d_0}^{\prime }{\cal E}_{k_\infty +1,d_\infty }^{\prime }\mu
^{*}\left( v\left( 1\otimes 1\right) \right) \mu ^{*}\left( \Delta
\right) }{-\hbar ^2\left( \hbar +c_0\right) \left( -\hbar +c_\infty \right) }%
= \\
\ &=&\int_{D}\frac{\rho
_1^{*}\left( t\right) ...\rho _{k_0}^{*}\left( t\right) \rho
_{k_0+1}^{*}\left( \tau \right) ...\rho _k^{*}\left( \tau \right) {\cal E}%
_{k_0+1,d_0}^{\prime }{\cal E}_{k_\infty +1,d_\infty }^{\prime }\mu
^{*}\left( v\sum_{i=0}^mT_i\otimes T_i^{\prime }\right) }{-\hbar ^2\left(
\hbar +c_0\right) \left( -\hbar +c_\infty \right) }= \\
\ &=&\int_{D}\frac{\rho
_1^{*}\left( t\right) ...\rho _{k_0}^{*}\left( t\right) \rho
_{k_0+1}^{*}\left( \tau \right) ...\rho _k^{*}\left( \tau \right) {\cal E}%
_{k_0+1,d_0}^{\prime }{\cal E}_{k_\infty +1,d_\infty }^{\prime }\mu
^{*}\left( \left( v\otimes v\right) \sum_{i=0}^mT_i\otimes T_i\right) 
}{-\hbar ^2\left( \hbar +c_0\right) \left( -\hbar +c_\infty \right) }= \\
\ &=&-\frac 1{\hbar ^2}\sum_i\int_{X_{k_0,d_0+1}}\frac{\rho _1^{*}\left(
t\right) ...\rho _{k_0}^{*}\left( t\right) \rho _0^{*}\left( T_i\right) 
{\cal E}_{k_0+1,d_0}}{\left( \hbar +c_0\right) }\cdot \\
& &\cdot
\int_{X_{k_\infty ,d_\infty
+1}}^T\frac{\rho _1^{*}\left( \tau \right) ...\rho _{k_\infty }^{*}\left(
\tau \right) \rho _\infty ^{*}\left( T_i\right) {\cal E}_{k_\infty
+1,d_\infty }}{\left( -\hbar +c_\infty \right) }
\end{eqnarray*}
\begin{flushright}
$\Box$
\end{flushright}

\noindent With similar arguments we compute the contribution of components of
type 2 and 3; for example, if $D=M_{d,0,k-1,1}^j$, then 
\begin{equation*}
\int_D\frac{\rho _1^{*}\left( x\right) ...\rho _k^{*}\left( x\right) {\cal E}%
_{k,\left( d,1\right) }}{Euler\left( {\cal N}_{D/Y_{k,\left( d,1\right)
}}\right) }=-\frac 1{\hbar ^2}\sum_i\int_{X_{k,d}}\frac{\rho _1^{*}\left(
t\right) ...\rho _{i-1}^{*}\left( t\right) \rho _{i+1}^{*}\left( t\right)
...\rho _k^{*}\left( t\right) \rho _0^{*}\left( T_i\right) {\cal E}_{k,d}}{%
(\hbar +c_0)}\text{.} \label{contrib3} 
\end{equation*}

\noindent By definition, we set $\int_{X_2,_0}\frac{\rho _1^{*}\left( T_i\right) \rho
_2^{*}\left( T_j\right) }{\hbar +c}$ $:=\left\langle T_i,T_j\right\rangle $.

We now set

\begin{definition}
\[
\psi _{ij}\left( t,\hbar \right) =\sum_{k=0}^\infty \frac
1{k!}\sum_{d=0}^\infty q^d\int_{X_{k+2,d}}\frac{\rho _1^{*}\left( t\right)
...\rho _k^{*}\left( t\right) \rho _{k+1}^{*}\left( T_i\right) \rho
_{k+2}^{*}\left( T_j\right) {\cal E}_{k+2,d}}{(\hbar +c)}\text{.}
\]
\end{definition}

\noindent Now, lemma \ref{contributo} easily implies:

\begin{theorem}
\label{decomp}
\[
-\hbar ^2\frac{\partial ^2{\cal F}^1\left( x\right) }{\partial t_i\partial
\tau _j}=\sum_k\psi _{ik}\left( t,\hbar \right) \psi _{kj}\left( \tau
,-\hbar \right) 
\]
\end{theorem}

{\bf Proof. }Take coordinates $\left\{ t_0,...,t_m,\tau _0,...,\tau
_m\right\} $ in $H${\bf \ }$_T^{*}(X)$ relative to the basis $\left\{
T_0\frac p\hbar ,...,T_m\frac p\hbar ,T_0\frac{\hbar -p}\hbar ,...,T_m\frac{%
\hbar -p}\hbar \right\} $.

\noindent Recall that 
\begin{equation*}
\frac{\partial ^2{\cal F}^1\left( x\right) }{\partial t_i\partial \tau _j}%
=\sum_{k\geq 1}\frac 1{k!}\sum_dq^dI_{\left( d,1\right) }^T\left(
x^k,T_i\frac p\hbar ,T_j\frac{\hbar -p}\hbar \right) \text{.} 
\end{equation*}

\noindent Summing up contributions of components of type 1, 2, 3, 4, 
\begin{equation*}
I_{\left( d,1\right) }^T\left( x^k,T_i\frac p\hbar ,T_j\frac{\hbar -p}\hbar
\right) = \label{sumcont} 
\end{equation*}
\begin{eqnarray*}
&=&-\frac 1{\hbar ^2}\left\langle T_i,T_j\right\rangle -\frac 1{\hbar
^2}k\sum_l\int_{X_{k+2,d}}\frac{\rho _1^{*}\left( t\right) ...\rho
_k^{*}\left( t\right) \rho _{k+1}^{*}\left( T_i\right) \rho _{k+2}^{*}\left(
T_l\right) {\cal E}_{k,d}}{(\hbar +c_0)}\left\langle T_l,T_j\right\rangle +
\\
&&\ \ -\frac 1{\hbar ^2}k\sum_l\int_{X_{k+2,d}}\frac{\rho _1^{*}\left( \tau
\right) ...\rho _k^{*}\left( \tau \right) \rho _{k+1}^{*}\left( T_j\right)
\rho _{k+2}^{*}\left( T_l\right) {\cal E}_{k,d}}{(-\hbar +c_\infty )}%
\left\langle T_l,T_i\right\rangle + \\
&&\ \ -\frac 1{\hbar ^2}\sum\begin{Sb} d_0,d_\infty >0  \\ k_0,k_\infty >0 
\end{Sb} \binom k{k_0}\sum_l\int_{X_{k_0+2,d}}\frac{\rho _1^{*}\left(
t\right) ...\rho _{k_0}^{*}\left( t\right) \rho _{k_0+1}^{*}\left(
T_i\right) \rho _{k_0+2}^{*}\left( T_l\right) {\cal E}_{k_0+2,d_0}}{(\hbar
+c_0)}\cdot \\
&&\ \cdot \int_{X_{k_\infty +2,d}}\frac{\rho _1^{*}\left( \tau \right)
...\rho _{k_\infty }^{*}\left( \tau \right) \rho _{k_\infty +1}^{*}\left(
T_l\right) \rho _{k_\infty +2}^{*}\left( T_j\right) {\cal E}_{k_\infty
+2,d_\infty }}{(-\hbar +c_\infty )}\text{.}
\end{eqnarray*}

\noindent Substituting in the formula for $\frac{\partial ^2{\cal F}^1\left( x\right) 
}{\partial t_i\partial \tau _j}$, and recalling that $\left\{
T_0,...,T_m\right\} $ is an orthonormal basis, the result follows.
\begin{flushright}
$\Box$
\end{flushright}

\begin{theorem}
\label{soluz}The matrix $\Psi \left( t,\hbar \right) $ is a fundamental
matrix of solutions for the differential equation 
\[
\nabla _\hbar s(x,t)=0\text{.}
\]
\end{theorem}

{\bf Proof.} $\left( \frac{\partial ^2{\cal F}^1\left( x\right) }{\partial
t_i\partial \tau _j}\right) =-\frac 1{\hbar ^2}\Psi \left( t,\hbar \right)
\Psi \left( \tau ,-\hbar \right) $ is a fundamental matrix of solutions.
Since the part of degree $0$ in $q$ of $\Psi \left( \tau ,-\hbar \right) $
is $-\frac 1{\hbar ^2}\left\langle T_i,T_j\right\rangle $ , this matrix is
invertible, and we can write 
\begin{equation*}
\Psi \left( t,\hbar \right) =-\hbar ^2\left( \frac{\partial ^2{\cal F}%
^1\left( x\right) }{\partial t_i\partial \tau _j}\right) \Psi ^{-1}\left(
\tau ,-\hbar \right) \text{.} 
\end{equation*}

\noindent Since $\Psi ^{-1}\left( \tau ,-\hbar \right) $ does not depend on $t$, and we
are differentiating with respect to the variables $\left\{ t_i\right\} $,
we have that $\Psi \left( t,\hbar \right) $ is another fundamental matrix of
solutions.
\begin{flushright}
$\Box$
\end{flushright}

In general, given any basis $\left\{ v_0,...,v_n\right\} $ for $H^{*}\left(
X\right) $, with $g_{ij}$ the corresponding intersection matrix, we can
define 
\begin{equation*}
\psi _{ij}\left( t,\hbar \right) =\sum_{k=0}^\infty \frac
1{k!}\sum_{d=0}^\infty q^d\int_{X_{k+2,d}}\frac{\rho _1^{*}\left( t\right)
...\rho _k^{*}\left( t\right) \rho _{k+1}^{*}\left( v_i\right) \rho
_{k+2}^{*}\left( v_j\right) {\cal E}_{k+2,d}}{(\hbar +c)}\text{,}  
\end{equation*}
and then prove in a similar way that 
\begin{equation*}
-\hbar ^2\frac{\partial ^2{\cal F}^1\left( x\right) }{\partial t_i\partial
\tau _j}=\sum_{k,l}\psi _{ik}\left( t,\hbar \right) g^{kl}\psi _{lj}\left(
\tau ,-\hbar \right) \text{.}  
\end{equation*}
This implies that theorem \ref{soluz} holds for any choice of the basis.

\noindent A quicker way to prove this is to observe that for any $\omega \in H^{*}\left(
X\right) $, the element 
\begin{eqnarray*}
\psi \left( \omega ,t,\hbar \right) &:&=  \label{psiomega} \\
\ &=&\sum_{i=0}^m\left( \sum_{k=0}^\infty \frac 1{k!}\sum_{d=0}^\infty
q^d\int_{X_{k+2,d}}\frac{\rho _1^{*}\left( t\right) ...\rho _k^{*}\left(
t\right) \rho _{k+1}^{*}\left( T_i\right) \rho _{k+2}^{*}\left( \omega
\right) {\cal E}_{k+2,d}}{(\hbar +c)}\right) T_i
\end{eqnarray*}
is a solution of the differential equation $\nabla _\hbar s(x,t)=0$, because
of the linearity of $\psi \left( \omega ,t,\hbar \right) $ in the variable $%
\omega $. From now on we will use a homogeneous basis with respect to the
grading, and denote the unity by $T_0$.

\bigskip\ 

\subsection{Restriction to $SQH^{*}\left( X\right) $}

We now pass to study how our solution $\Psi \left( t,\hbar \right) $
restricts for $t$ $=t_1T_1$, i.e. what happens if we consider our
connections restricted to the tangent bundle to $H^2\left( X\right) $ . We
recall that by our assumptions $H^2\left( X\right) $ (or $H^2\left( X\right)
\cap V$) is one-dimensional. The solutions of $\nabla _\hbar =0$ will be
considered as elements in the small quantum cohomology ring of $X$.

\noindent Our aim is to prove:

\begin{theorem}
\label{redsmall}Let $s_{ij}\left( t,\hbar \right) :=\psi _{ij}\left( t,\hbar
\right) _{\mid t\in H^2\left( X\right) }$. Then, for every $i,j$, we have 
\[
s_{ij}\left( t,\hbar \right) =\sum_{d\geq 0}q^d\exp \left( dt_1\right)
\int_{X_{2,d}}\exp \left( \frac{\rho _1^{*}\left( t\right) }\hbar \right) 
\frac{\rho _1^{*}\left( T_i\right) \rho _2^{*}\left( T_j\right) }{\hbar +c}%
{\cal E}_{2,d}\text{.}\label{esseij}
\]
\end{theorem}

\noindent The proof of this theorem requires some simple facts. Consider the map 
\begin{eqnarray*}
\pi &:&X_{k,d}\rightarrow X_{k-1,d} \\
\left[ C,x_1,...,x_k,\psi \right] &\rightarrow &\left[
C,x_1,...,x_{k-1},\psi \right]
\end{eqnarray*}
which forgets the last point, and the map $\sigma _1:X_{k-1,d}\rightarrow
X_{k,d}$ , the section given by the first marked point. Let $D:=\sigma
_1\left( X_{k-1,d}\right) $, and $L:={\cal O}\left( D\right) $. Set $F_1^k$
equal to the line bundle on $X_{k,d}$ whose fibers are the tangent lines at
the first marked points.

\begin{lemma}
\begin{enumerate}
\item  $\sigma _1^{*}\left(F_1^k\right) $ is trivial.

\item  $\pi ^{*}\left( F_1^{k-1}\right) =F_1^k\otimes L$.
\end{enumerate}
\end{lemma}

{\bf Proof.} 1) By definition of the section $\sigma _1$, the divisor $D$
consists of those stable maps $\left[ C,x_1,...,x_k,\psi \right] $ such that 
$C$ splits in two irreducible components $C_1$ and $C_2$, such that $C_1$
contains exactly the two marked points $x_1$ and $x_k$, and $\psi $
restricted to $C_1$ is constant. This immediately implies our claim since
the line bundle $\sigma _1^{*}\left( F_1^k\right) $ is the pull back of a
bundle on $\overline{{\cal M}}_{0,3}=\left\{ pt\right\} $.

\noindent 2) It is clear that $\pi ^{*}\left( F_1^{k-1}\right) $ and $F_1^k$ are
canonically isomorphic outside $D$; thus $\pi ^{*}\left( F_1^{k-1}\right)
\otimes \left( F_1^k\right) ^{-1}=L^a$, for some $a$, since it has
support only on $D$. To compute the integer $a$, let us restrict to $D$. We
have already used that the fiber of $L_{\mid D}={\cal N}_{D/X_{k,d}}\,$ is $%
T_pC_1\otimes T_pC_2$, where $p$ is the join point.

\noindent In this case, the first factor is trivial, and the second equals $\pi
^{*}\left( F_1^{k-1}\right) $, hence $L_{\mid D}=\pi ^{*}\left(
F_1^{k-1}\right) $.

\noindent On the other hand, by 1) we deduce that $\left( F_1^k\right) _{\mid
D}=\sigma _1^{*}\left( F_1^k\right) $ is trivial. Substituting, we get $%
\left( \pi ^{*}\left( F_1^{k-1}\right) \otimes \left( F_1^k\right)
^{-1}\right) _{\mid D}=\pi ^{*}\left( F_1^{k-1}\right) $, and $a=1$. 
\begin{flushright}
$\Box$
\end{flushright}

We set $c\left( k\right) :=c_1\left( F_1^k\right) \in H^2\left(
X_{k,d}\right) $ ; with abuse of notation, we will write only $c$ if it is
clear to which cohomology ring it belongs.

\begin{lemma}
Let $\pi :X_{k,d}\rightarrow X_{k-1,d}$ be the map which forgets the last
marked point, and let $a\in H^{*}\left( X_{k,d}\right) $. Then for any $%
r\geq 1$, 
\[
\pi _{*}\left( ac^r\right) =\pi _{*}a\left( c^r\right) -\sigma _1^{*}\left(
a\right) c^{r-1}\text{.}
\]
\end{lemma}

{\bf Proof. }From the previous lemma , if we set $x:=\sigma _{*}\left(
1\right) =c_1(L)$, we get: 
\begin{equation*}
\left\{ 
\begin{array}{c}
\sigma _1^{*}\left( c^r\right) =0 \\ 
\pi ^{*}\left( c\left( k-1\right) \right) =c\left( k\right) +x
\end{array}
\right. \text{.}  
\end{equation*}

\noindent For $r=1$ : 
\begin{eqnarray*}
\pi _{*}\left( ac\left( k\right) \right) &=&\pi _{*}\left\{ a\left[ \pi
^{*}c\left( k-1\right) -x\right] \right\} =\pi _{*}\left[ a\pi ^{*}c\left(
k-1\right) \right] -\pi _{*}\left( ax\right) =  \\
&=&\pi _{*}\left( a\right) c\left( k-1\right) -\pi _{*}\left( a\sigma
_{1*}\left( 1\right) \right) =\pi _{*}\left( a\right) c\left( k-1\right)
-\pi _{*}\sigma _{1*}\left( 1\cdot \sigma _1^{*}\left( a\right) \right) = \\
&=&\pi _{*}\left( a\right) c\left( k-1\right) -\sigma _1^{*}\left( a\right) 
\text{.}  
\end{eqnarray*}

\noindent By induction, if it is true for $r-1$, 
\begin{eqnarray*}
\pi _{*}\left( ac^r\right) &=&\pi _{*}\left( ac^{r-1}c\right) =\pi
_{*}\left[ ac^{r-1}\left( \pi ^{*}c-x\right) \right] =\pi _{*}\left(
ac^{r-1}\right) c-\pi _{*}\left[ ac^{r-1}\sigma _{1*}\left( 1\right) \right]
=   \\
&=&\pi _{*}\left( ac^{r-1}\right) c-\pi _{*}\sigma _{1*}\left[ 1\cdot \sigma
_1^{*}\left( ac^{r-1}\right) \right] =\pi _{*}\left( ac^{r-1}\right)
c-\sigma _1^{*}\left( ac^{r-1}\right) = \\
&=&\pi _{*}a\left( c^r\right) -\sigma _1^{*}\left( a\right) c^{r-1}-\sigma
_1^{*}\left( ac^{r-1}\right) =\pi _{*}a\left( c^r\right) -\sigma
_1^{*}\left( a\right) c^{r-1}\text{.}  
\end{eqnarray*}
\begin{flushright}
$\Box$
\end{flushright}

{\bf Proof }of \ref{redsmall}. Recall that 
\begin{eqnarray*}
s_{ij}(t,\hbar ) &=&\psi _{ij}\left( t,\hbar \right) _{\mid t\in H^2\left(
X\right) }\text{, and} \\
\psi _{ij}\left( t,\hbar \right) &=&\sum_{k=0}^\infty \frac
1{k!}\sum_{d=0}^\infty q^d\int_{X_{k+2,d}}\frac{\rho _1^{*}\left( T_i\right)
\rho _2^{*}\left( T_j\right) \rho _3^{*}\left( t\right) ...\rho
_{k+2}^{*}\left( t\right) {\cal E}_{k+2,d}}{(\hbar +c)}\text{.}
\end{eqnarray*}

\noindent Let $a\in H^{*}\left( X\right) $. Developing in power series, we have 
\begin{equation*}
\frac a{\hbar +c}=\sum_{r=0}^\infty \frac{\left( -1\right) ^r}{\hbar ^{r+1}}%
ac^r\text{,} 
\end{equation*}
and 
\begin{eqnarray*}
\pi _{*}\left( \frac a{\hbar +c}\right) &=&\sum_{r=0}^\infty \frac{\left(
-1\right) ^r}{\hbar ^{r+1}}\pi _{*}\left( ac^r\right) =\sum_{r=0}^\infty 
\frac{\left( -1\right) ^r}{\hbar ^{r+1}}\left[ \pi _{*}\left( a\right)
\left( c^r\right) -\sigma _1^{*}\left( a\right) c^{r-1}\right] =
 \\
\ &=&\frac{\pi _{*}\left( a\right) }{\hbar +c}-\frac{\sigma _{*}\left(
a\right) }{\hbar \left( \hbar +c\right) }\text{.}  
\end{eqnarray*}

\noindent Moreover, if $y\in H^{*}\left( X_{k,d}\right) $, and $y^{\prime }=\pi
_{*}\left( y\right) $, then 
\begin{eqnarray*}
\pi _{*}\left( \frac{y\rho _k^{*}\left( t\right) }{\hbar +c}\right) &=&\frac{%
\pi _{*}\left( y\rho _k^{*}\left( t\right) \right) }{\hbar +c}-\frac{\sigma
_{*}\left( y\rho _k^{*}\left( t\right) \right) }{\hbar \left( \hbar
+c\right) }=   \\
\ &=&\frac{y^{\prime }\pi _{*}\rho _k^{*}\left( t\right) }{\hbar +c}-\frac{%
\sigma _{*}\rho _k^{*}\left( t\right) }\hbar \frac{y^{\prime }}{\hbar +c}= \\
\ &=&\frac{y^{\prime }\left\langle d,t\right\rangle }{\hbar +c}-\frac{\rho
_1^{*}\left( t\right) }\hbar \frac{y^{\prime }}{\hbar +c}\text{.}
\end{eqnarray*}

\noindent Using that ${\cal E}_{k,d}=\pi ^{*}\left( {\cal E}_{k-1,d}\right) $, apply
this last formula to obtain: 
\begin{equation*}
\int_{X_{k+2,d}}\frac{\rho _1^{*}\left( T_i\right) \rho _2^{*}\left(
T_j\right) \rho _3^{*}\left( t\right) ...\rho _{k+2}^{*}\left( t\right) 
{\cal E}_{k+2,d}}{(\hbar +c)}=  
\end{equation*}
\begin{equation*}
=\int_{X_{k+1,d}}\frac{\pi _{*}\left( \rho _1^{*}\left( T_i\right) \rho
_2^{*}\left( T_j\right) \rho _3^{*}\left( t\right) ...\rho _{k+2}^{*}\left(
t\right) \right) {\cal E}_{k+1,d}}{(\hbar +c)}= 
\end{equation*}
\begin{eqnarray*}
\ &=&\int_{X_{k+1,d}}\frac{\rho _1^{*}\left( T_i\right) \rho _2^{*}\left(
T_j\right) \rho _3^{*}\left( t\right) ..\rho _{k+1}\left( t\right) }{\hbar +c%
}\left( dt_1-\frac{\rho _1^{*}\left( t\right) }\hbar \right) {\cal E}%
_{k+1,d}= \\
\ &=&...=\int_{X_{2,d}}\frac{\rho _1^{*}\left( T_i\right) \rho _2^{*}\left(
T_j\right) }{\hbar +c}\left( dt_1-\frac{\rho _1^{*}\left( t\right) }\hbar
\right) ^k{\cal E}_{2,d}= \\
\ &=&\int_{X_{2,d}}\frac{\rho _1^{*}\left( T_i\right) \rho _2^{*}\left(
T_j\right) }{\hbar +c}\sum_{u=0}^k\binom ku\left( dt_1\right) ^u\left( \frac{%
\rho _1^{*}\left( t\right) }\hbar \right) ^{k-u}{\cal E}_{2,d}\text{.}
\end{eqnarray*}

\noindent Hence 
\begin{eqnarray*}
s_{ij}\left( t,\hbar \right) &=&\sum_{d=0}^\infty q^d\sum_{k=0}^\infty \frac
1{k!}\sum_{u=0}^k\binom ku\left( dt_1\right) ^u\int_{X_{2,d}}\frac{\rho
_1^{*}\left( T_i\right) \rho _2^{*}\left( T_j\right) }{\hbar +c}^u\left( 
\frac{\rho _1^{*}\left( t\right) }\hbar \right) ^{k-u}{\cal E}_{2,d}=
 \\
\ &=&\sum_{d\geq 0}q^d\exp \left( dt_1\right) \int_{X_{2,d}}\exp \left( 
\frac{\rho _1^{*}\left( t\right) }\hbar \right) \frac{\rho _1^{*}\left(
T_i\right) \rho _2^{*}\left( T_j\right) }{\hbar +c}{\cal E}_{2,d}\text{.}
\end{eqnarray*}
\begin{flushright}
$\Box$
\end{flushright}

\bigskip\ 

A special role will be played by the collection of functions 
\begin{equation*}
s_\beta \left( t,\hbar \right) :=s_{0 \beta}\left( t,\hbar \right) \text{,} 
\end{equation*}
i.e. by the first components of the vector solutions of $\nabla _\hbar s=0$.

\noindent These will be the functions to manipulate in order to get solutions for
Picard Fuchs equation, and these give enough information about enumerative
geometry in $X.$

\noindent Let $P(\hbar \frac d{dt_1},\exp t_1,\hbar )$ denote a non-commutative
polynomial differential operator in the variables $\hbar ,\frac d{dt_1},\exp
t_1$, that is 
\begin{equation*}
P(\hbar \frac d{dt_1},\exp t_1,\hbar )=\sum_rQ_r(\exp t_1,\hbar )\left(
\hbar \frac d{dt_1}\right) ^r\text{,} 
\end{equation*}

\begin{proposition}
\label{enumerative}
Take $P(\hbar \frac d{dt_1},\exp t_1,\hbar )$ as above. If 
\[
Ps_\beta \left( t,\hbar \right) =0\text{ }
\]
for every $\beta =0,...,m$, then the following relation holds in $%
SQH^{*}\left( X\right) $: 
\[
P\left( T_1,q,0\right) =0\text{.}
\]
\end{proposition}

{\bf Proof. }The functions $s_\beta $ form the first row in the fundamental
solution matrix $S=\left\{ s_{ij}\right\} $. Let $M_1$ be the matrix of
quantum multiplication by the class $T_1$. Obviously, 
\begin{equation*}
\left( M_1\cdot S\right) _{ij}=\sum_ks_{kj}\phi _{1ki}=\hbar \frac
d{dt_1}s_{ij}\text{;} 
\end{equation*}
on the other hand, as functions $\phi _{1ki}$ depend on $t_1$, we have 
\begin{equation*}
\left( \hbar \frac d{dt_1}\right) ^rs_{ij}=T_1^r*s_{ij}\text{ }+\hbar \cdot 
\widetilde{Q}\text{,} 
\end{equation*}
where $\widetilde{Q}$ is the multiplication by a matrix whose entries are
functions of the variable $T_1$ and depend polynomially on $\dot \hbar $.

\noindent Develop the polynomial $P$ with respect to the last variable, i.e. $%
P=P_0+\hbar P_1+...+\hbar ^rP_r$, and 
\begin{equation*}
P=P_0+\hbar P_1+...+\hbar ^rP_r\text{,} 
\end{equation*}
reducing modulo $\hbar $, we get: 
\begin{equation*}
\left( P_0(M_1,\exp t_1)\cdot S\right) _{0 \beta}=\left( P(M_1,\exp
t_1,0)\cdot S\right) _{0 \beta }=0\text{.} 
\end{equation*}

\noindent We know that $S$ is a non-degenerate matrix, and therefore 
\begin{equation*}
\left( P(M_1,\exp t_1,0)\right) _{0 \beta}=0\text{,} 
\end{equation*}
for every $\beta =0,...,m$.

Next goal is to prove that the first column of the matrix $P(M_1,\exp t_1,0)$ vanishes. 
This would imply that this quantum operator, applied to the
unit $T_{0}$, vanishes, and therefore $P(T_1,q,0)$ in $SQH^{*}\left( X\right) $.
Since we can choose an orthogonal basis, for every index $i$, writing down explicitely
the expression of $\left( P(M_1,\exp t_1,0)\right) T_{0} * T_{1}$ and using the fact that
$T_{0}$ is a unit for the quantum product, we obtain that the $i$-th element of the first 
column of our matrix is $0$.$\Box$

\newpage

\part{Solutions of Picard-Fuchs equation}

From now on let $X$ be a smooth projective complete intersection in ${\Bbb P}%
^n$ given by $r$ equations of degrees $l_i$, $1\leq i\leq r$ and denote by $W
$ the vector bundle $\oplus_{i=1}^r{\cal O}(l_i)$ over ${\Bbb P}^n$%
. We will be interested in the following cases:\\ 

i) $\sum\limits_{i=1}^rl_i<n$;

ii) $\sum\limits_{i=1}^rl_i=n;$

iii) $\sum\limits_{i=1}^rl_i=n+1$.

\noindent Observe that case iii) implies that $X$ is Calabi-Yau. It is useful for our
purpose to recall that the only Calabi-Yau threefolds are:

\begin{itemize}
\item  a quintic in ${\Bbb P}^4$,

\item  the intersection of a quadric and a quartic in ${\Bbb P}^5$,

\item  the intersection of two cubics in ${\Bbb P}^5$,

\item  the intersection of two quadrics and a cubic in ${\Bbb P}^6$.
\end{itemize}

\bigskip

\noindent In each case we can consider a differential equation as explained in the
Introduction. We are going to prove how solutions of this equation can be
recovered from functions $s_\beta $ introduced in the previous section, i.e.
from ``manipulated''solutions of the equation $\nabla _\hbar s=0$. To this
extent, we need to find an explicit form of a suitable class $S(t,\hbar )$
in $H^{*}{\bf (}{\Bbb P}^n)$: we will then make use of equivariant
cohomology theory and graph theory.

\section{Fixed points on the moduli spaces of stable maps.}

The action of the $(n+1)-$ dimensional torus $T^{n+1}$ on ${\Bbb P}^n$
induces an action of $T^{n+1}$ on $\overline{{\cal M}}_{0,2}({\Bbb P}^n,d)$.
We will describe the set of fixed points. Denote by $p_i$, $0\leq i\leq 0$,
the fixed points of $T^{n+1}$ acting on ${\Bbb P}^n$: $p_i$ is the
projectivization of the $i-th$ coordinate line in ${\Bbb C}^{n+1}$. Also
denote $l_{ij}=l_{ji}$, $i\neq j$, the line in ${\Bbb P}^n$ through $p_i$
and $p_j$. Finally, we call $\Phi $ the configuration of fixed points $p_i$
and lines $l_{ij}$.

\noindent Let $\left[ C,x_0,x_1;f\right] $, with $C=\bigcup_\alpha C_\alpha $, be a
fixed point in $Y_{2,d}$. By the definition of the action on $Y_{2,d}$, the
geometric image $f(C)$ should be contained in $\Phi $, so it is a union of
lines $l_{ij}$. Obviously, the images of all special points, as being in
finite number, are points $p_i$. Suppose now $C_\alpha $ has less than two
special points. Because of the stability of $f$, this means that $f_\alpha =f_{\mid
C_\alpha }:C_\alpha \rightarrow l_{ij}={\Bbb P}^1$ is a degree $d_\alpha $
covering. When $d_\alpha >1$, notice that the only branch points of $%
f_\alpha $ are $p_i$ and $p_j$, since an element $t\in T^{n+1}$sends the
fiber over a point $q\notin \{p_i,p_j\}$ in the fiber over the point $t\cdot
q$. On the other hand, since the fiber of $f_\alpha $ over a fixed point
is an orbit with respect to the action of the torus, $f_\alpha $ is a
totally ramified covering with two ramification points.

\noindent If we introduce homogeneous coordinates $z_1,z_2$ on $C_\alpha $, then we
have the following diagram 
\begin{equation*}
\begin{array}{ccc}
\left[ z_1,z_2\right] & \longrightarrow & [t_i^{\frac 1{d_\alpha
}}z_1,t_j^{\frac 1{d_\alpha }}z_2] \\ 
\downarrow &  & \downarrow \\ 
\left[ 0,\ldots ,X_i,\ldots ,X_j,\ldots ,0\right] & \longrightarrow & \left[
0,\ldots ,t_iX_i,\ldots ,t_jX_j,\ldots ,0\right]
\end{array}
\text{,} 
\end{equation*}

\noindent with $X_i=z_1^{d_\alpha }$ and $X_j=z_2^{d_\alpha }$.

\noindent By what we have just said, it is also clear that irreducible components with
more than two special points are mapped to fixed points in ${\Bbb P}^n$.

\smallskip\ 

We now describe components of fixed points in terms of a special class of graphs.
In the sequel, we shall consider the family ${\cal G}$ of finite graphs $%
\Gamma $ with $V_\Gamma $ and $E_\Gamma $ respectively its set of vertices
and edges. Each graph $\Gamma $ is also equipped with three additional data:

\begin{itemize}
\item  a map $\delta :V_\Gamma \rightarrow \{0,\ldots ,n\}$,

\item  a map $\lambda :E_\Gamma \rightarrow {\Bbb Z}^{+}$,

\item  a map $P:V\rightarrow {\cal P}(\left\{ 0,1\right\} )$ such that $%
\left\{ 0,1\right\} $ is the disjoint union of the $P(v)$, for $v\in
V_\Gamma $.
\end{itemize}

\noindent We set $h_v=\#P\left( v\right) $, and denote by $l_v$ the valency of $v$,
that is, the cardinality of the set $L_v$ of edges issuing from $v$.
Finally, we say that two graphs $\Gamma _1$ and $\Gamma _2$ are isomorphic
if there exist bijections between $V_{\Gamma _1}$ and $V_{\Gamma _2}$, $%
E_{\Gamma _1}$ and $E_{\Gamma _2}$ which preserve adjacency and additional
data. In the sequel, we shall denote by $[\Gamma ]$ an equivalence class of
graphs in ${\cal G}$.

\smallskip\ 

We now construct a graph $\Gamma $ associated to a fixed point $\left[
C,x_0,x_1;f\right] $ (see \cite{K}). The vertices $v\in V_\Gamma $
correspond to the connected components $C_v$ of $f^{-1}(p_0,\ldots ,p_n)$,
i.e. points or non-empty union of contracted components of $C$; the edges $%
\alpha \in E_\Gamma $ correspond to irreducible components $C_\alpha $ of
genus zero mapping to lines $l_{ij}$. Moreover we define the map $P$
so that it assigns to each vertex the indices of marked points lying in $C_v$%
. It is quite obvious that $\Gamma $ is a connected tree (there are no
simple loops for the choice of the set of vertices). Notice that the map $%
\delta $ is well defined!

\begin{figure}[htb]
\hspace{-2cm} 
\begin{center} 
\mbox{\epsfig{file=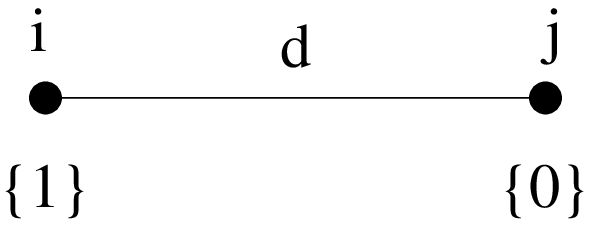,width=2cm,height=4cm,angle=270}}
\caption[]{
\label{dis2}}
\end{center}
\end{figure}

\noindent Alternatively, given an element $\left[ C,x_0,x_1;f\right] \in \overline{%
{\cal M}}_{0,2}({\Bbb P}^n,d)^{T^{n+1}}$, we associate to it a graph in $%
{\cal G}$ . Each edge $\alpha $ connecting two vertices, with labels $i$ and 
$j$, corresponds to an irreducible component mapped to $l_{ij}$ with degree $%
\lambda (\alpha )$, whereas each vertex $v$ with $\delta (v)=k$ corresponds
to connected components $C_v$ mapped to $p_k$. Finally, if $P(v)\neq
\emptyset $, then $C_v$ contains marked points with indices prescribed by $%
P(v)$. Notice that we don't fix the number of connected components
contracted to points. For example, if $\Gamma $ is the graph of figure \ref
{dis2},\ then the corresponding element $\left[ C,x_0,x_1;f\right] $ is a
point, since degrees $d$ maps from $C$ to ${\Bbb P}^1$ with two ramification
points define the same element in $\overline{{\cal M}}_{0,2}({\Bbb P}^1,d)$.

\smallskip\ 

Let $\Sigma _\Gamma $ be the locus of fixed points in $\overline{{\cal M}}%
_{0,2}({\Bbb P}^n,d)$ whose associated graph is $\Gamma $. Denote by $%
\widetilde{V}_\Gamma $ the set of vertices for which $h_{v_i}+l_{v_i}-3\geq
0 $ and choose, for each $v\in \widetilde{V}_\Gamma $, an ordering of $L_v$.

\begin{proposition}
\label{graph} 
\[
\Sigma _\Gamma \cong \frac{\prod\limits_{i=1}^s\overline{{\cal M}}%
_{0,h_{v_i}+l_{v_i}}}{Aut(\Gamma )}\text{.}
\]
\end{proposition}

{\bf Proof. }Consider the morphism 
\begin{equation*}
\xi _\Gamma :\prod\limits_{i=1}^s\overline{{\cal M}}_{0,h_{v_i}+l_{v_i}}%
\longrightarrow \overline{{\cal M}}_{0,2}({\Bbb P}^n,d) 
\end{equation*}
which is defined as follows. The image of a point 
\begin{equation*}
\left( \lbrack C_1;x_{v_1,1},\ldots ,x_{v_1,h_{v_1+l_{v_1}}}],\ldots
,[C_s;x_{v_s,1},\ldots ,x_{v_s,h_{v_s}+l_{v_s}}]\right) 
\end{equation*}
is an element associated with the graph $\Gamma $ whose contracted
components are given by $C_1,\ldots ,C_s$. One can easily verify that $%
\Sigma _\Gamma $ is the image of $\xi _\Gamma $ and does not depend on the
choice of orderings of the $L_v$, $v\in \widetilde{V}_\Gamma $. Since two
elements in the same fiber of $\xi _\Gamma $ differ by an automorphism of $%
\Gamma $, the Proposition is proved.
\begin{flushright}
$\Box$
\end{flushright}

\bigskip\ 

\bigskip\ 

\section{Projective complete intersections with $l_1+\ldots +l_r<n$.}

Let $X$ be a smooth projective complete intersection in ${\Bbb P}^n$ given
by $r$ equations of degrees $l_i$, $1\leq i\leq r$, with $%
\sum\limits_{i=1}^rl_i<n$, and let $m=n-r$. This is the easier case: in
fact, in this section we are going to prove that the functions $s_\beta
(t,\hbar ):=s_{0 \beta}(t,\hbar )$, $\beta =0,\ldots ,m$, (components of
vector solutions $\left( s_{0 \beta}(t,\hbar ),\ldots s_{m \beta}(t,\hbar
)\right) $ of the equation $\nabla _\hbar s=0$) themselves satisfy the
Picard-Fuchs equation for the mirror symmetric family of $X$. To this
purpose, we shall give an explicit form of $s_\beta (t,\hbar )$. Fix a basis 
$\{T_\beta =P^\beta \}$ of $H^{*}({\Bbb P}^n)$ and recall that 
\begin{equation*}
s_\beta (t,\hbar )=\sum\limits_{d\geq 0}e^{dt}\int\limits_{\overline{{\cal M}%
}_{0,2}({\Bbb P}^n,d)}(e_0)^{*}\left( P^\beta e^{\frac{Pt}\hbar }\right) 
\frac{{\cal E}_{2,d}}{\hbar +c}\text{, }\beta =0,\ldots ,m\text{,} 
\end{equation*}
with $e_0$ the evaluation on the first marked point. We can also view $%
s_\beta (t,\hbar )$ as components of the function

\begin{equation}
S(t,\hbar )=e^{\frac{Pt}\hbar }\sum\limits_{d\geq 0}e^{dt}(e_0)_{*}\left( 
\frac{{\cal E}_{2,d}}{\hbar +c}\right) \in H^{*}({\Bbb P}^n)\text{,}
\label{esse}
\end{equation}
since by the integration formula, we obtain 
\begin{equation}
s_\beta (t,\hbar )=\int\limits_{{\Bbb P}^n}P^\beta S(t,\hbar )\text{.}
\label{essebeta}
\end{equation}

For $d=0$, when $\overline{{\cal M}}_{0,2}({\Bbb P}^n,d)$ is not defined, we
let $Euler(\oplus_{i=1}^r{\cal O}(l_i))$ play the role of $%
(e_0)_{*}({\cal E}_{2,0})$.

\medskip\ 

\begin{theorem}
\label{one}{\sl \ }Suppose that $l_1+\ldots +l_r<n$. Then 
\[
S(t,\hbar )=e^{\frac{Pt}\hbar }\sum_{d\geq 0}e^{dt}\frac{%
\prod_{m=0}^{dl_1}(l_1P+m\hbar )\ldots \prod_{m=0}^{dl_r}(l_rP+m\hbar )}{%
\prod_{m=1}^d(P+m\hbar )^{n+1}}\text{.}
\]
\end{theorem}

\noindent We will postpone the proof of this Theorem at the end of this section.

\medskip\ 

\begin{corollary}
\label{two}{\sl \ }The components $s_\beta (t,\hbar )$, $\beta =0,\ldots ,m$%
, form a basis of solutions of the Picard-Fuchs equation
\end{corollary}

\begin{equation*}
(\hbar \frac d{dt})^{n+1-r}F(t,\hbar
)=e^t\prod_{j=1}^rl_j\prod_{m=1}^{l_j-1}\hbar (l_j\frac d{dt}+m)F(t,\hbar ) 
\end{equation*}
for the mirror symmetric family of $X$.

{\bf Proof.\ }Set

\begin{equation*}
\Gamma (d)=\frac{\prod_{m=0}^{dl_1}(l_1P+m\hbar )\ldots
\prod_{m=0}^{dl_r}(l_rP+m\hbar )}{\prod_{m=1}^d(P+m\hbar )^{n+1}}\text{.} 
\end{equation*}
Applying differential operators $(\hbar \frac d{dt})^{n+1-r}$ and $%
e^t\prod_{j=1}^rl_j\prod_{m=1}^{l_j-1}\hbar (l_j\frac d{dt}+m)$ to $%
S(t,\hbar )$, we obtain 
\begin{equation*}
\sum\limits_{d\geq 0}(P+d\hbar )^{n+1-r}\Gamma (d)e^{(\frac P\hbar +d)t} 
\end{equation*}

\noindent and 
\begin{equation*}
\prod_{j=1}^rl_j\sum\limits_{d\geq 0}\prod\limits_{m=1}^{l_j-1}\left[
l_j(P+d\hbar )+m\hbar \right] e^{(\frac P\hbar +d+1)t}\Gamma (d)\text{.} 
\end{equation*}

On the other hand we have, 
\begin{equation*}
\left[ P+(d+1)\hbar \right] ^{n+1-r}\Gamma
(d+1)=\prod_{j=1}^rl_j\sum\limits_{d\geq 0}\prod\limits_{m=1}^{l_j-1}\left[
l_j(P+d\hbar )+m\hbar \right] \Gamma (d)\text{,} 
\end{equation*}
since
\begin{equation*}
\begin{array}{c}
\left[ P+(d+1)\hbar \right] ^{n+1-r}\dfrac{\prod_{m=0}^{(d+1)l_1}(l_1P+m%
\hbar )\ldots \prod_{m=0}^{(d+1)l_r}(l_rP+m\hbar )}{\prod_{m=1}^{d+1}(P+m%
\hbar )^{n+1}}= \\ 
=\Gamma (d)\dfrac{\prod_{m=1}^{l_1-1}(l_1P+(dl_1+m)\hbar )\ldots
\prod_{m=1}^{l_r-1}(l_rP+(dl_r+m)\hbar )}{\left[ P+(d+1)\hbar \right] ^r}%
\prod_{i=1}^r\left[ l_iP+(dl_i+l_i)\hbar \right] = \\ 
=\prod_{j=1}^rl_j\sum\limits_{d\geq 0}\prod\limits_{m=1}^{l_j-1}\left[
l_j(P+d\hbar )+m\hbar \right] \Gamma (d)\text{.}
\end{array}
\end{equation*}
\begin{flushright}
$\Box$
\end{flushright}

\medskip\ 

\begin{corollary}
\label{three}{\sl \ }If $\dim X\neq 2$, the cohomology class $p$ of
the hyperplane section satisfies  the relation 
$$p^{n+1-r}=l_1^{l_1}\ldots l_r^{l_r}qp^{l_1+\ldots +l_r-r}$$
in the quantum cohomology of $X$.
\end{corollary}

{\bf Proof. }It follows easily from observations made in section 2.
\begin{flushright}
$\Box$
\end{flushright}

\noindent This corollary is consistent with results of Beauville (\cite[Be]{Be}) and
Jinzenji (\cite[Ji]{Ji}) on quantum cohomology of projective hypersurfaces.

\medskip\ 

\begin{example}
If $l_1=\ldots =l_r=1$, then the relation (cfr. section 1) for the small
quantum cohomology of ${\Bbb P}^{n-r}$is
\end{example}

\begin{equation*}
p^{n+1-r}=qp\text{.} 
\end{equation*}

\bigskip\ 

Givental deduces Theorem \ref{one} from its equivariant generalization. To
this purpose, let $T=(S^1)^k$ act diagonally on ${\Bbb P}^n$, and $T^{\prime
}=(S^1)^r$ acts trivially on ${\Bbb P}^n$ and diagonally on the fiber of the
vector bundle $\oplus_{i=1,\ldots ,r}{\cal O}(l_i)$. Recall that 
\begin{equation*}
H_{T\times T^{\prime }}^{*}({\Bbb P}^n)=\frac{{\Bbb Q}\left[ p,\lambda
_0,\ldots ,\lambda _k\right] }{\prod\limits_{i=0}^k(p-\lambda _i))}%
\otimes {\Bbb Q}\left[ \mu _1,\ldots ,\mu _r\right] 
\end{equation*}
and the equivariant Euler class of the vector bundle $\oplus_{i=1,\ldots
,r}{\cal O}(l_i)$ is equal to $\prod\limits_{i=1}^r(l_ip-\mu _i)$. Introduce
the equivariant counterpart $S^{\prime }$ of the class $S$ in the $T\times
T^{^{\prime }}$ equivariant cohomology of ${\Bbb P}^n$. This means that we
use the equivariant class $p$ instead of $P$ and replace the Euler classes $%
{\cal E}_{2,d}$ and $c$ by their equivariant partners.

\medskip\ 

\begin{theorem}
\label{four}Let $S^{\prime }(t,\hbar )=e^{\frac{pt}\hbar }\sum\limits_{d\geq
0}e^{dt}(e_0)_{*}\left( \frac{{\cal E}_{2,d}}{\hbar +c}\right) $.{\sl \ }%
Suppose that $l_1+\ldots +l_r<n$. Then
\end{theorem}

\begin{equation*}
S^{\prime }(t,h)=e^{\frac{pt}\hbar }\sum\limits_{d\geq 0}e^{dt}\frac{%
\prod_{m=0}^{dl_1}(l_1p-\mu _1+m\hbar )\ldots \prod_{m=0}^{dl_r}(l_rp-\mu
_r+m\hbar )}{\prod_{m=1}^d(p-\lambda _0+m\hbar )\ldots
\prod_{m=1}^d(p-\lambda _n+m\hbar )}\text{.} 
\end{equation*}

\medskip\ 

\begin{remark}
\label{eqnoneq}
\end{remark}

\noindent The inclusion of a fiber 
\begin{equation*}
i:{\Bbb P}^n\hookrightarrow {\Bbb P}_{T\times T^{\prime }}^n 
\end{equation*}
induces a morphism 
\begin{equation*}
\begin{array}{cccc}
i^{*}: & H_{T\times T^{\prime }}^{*}\left( {\Bbb P}^n\right) & \rightarrow & 
H^{*}\left( {\Bbb P}^n\right) \\ 
& q(p,\lambda ,\mu ) & \rightarrow & q(p,0,0)
\end{array}
\text{,} 
\end{equation*}
which is the cup product with the fundamental class of the fiber. Moreover
the commutative diagram 
\begin{equation*}
\begin{array}{ccc}
{\Bbb P}^n & \hookrightarrow & {\Bbb P}_{T\times T^{\prime }}^n \\ 
\downarrow \pi _{\mid {\Bbb P}^n} &  & \downarrow \pi \\ 
\left\{ pt\right\} & \hookrightarrow & \left( {\Bbb P}^\infty \right)
^{n+1+r}
\end{array}
\end{equation*}
induces a diagram in cohomology: 
\begin{equation*}
\begin{array}{ccc}
H_{T\times T^{\prime }}^{*}\left( {\Bbb P}^n\right) & \rightarrow & 
H^{*}\left( {\Bbb P}^n\right) \\ 
\downarrow \pi _{*} &  & \downarrow \pi _{*} \\ 
H^{*}\left( {\Bbb P}^\infty \right) ^{n+1+r} & \rightarrow & H^{*}\left(
\left\{ pt\right\} \right)
\end{array} 
\end{equation*}
and this implies that if we simply put $\lambda _i=0$, $1\leq i\leq n$, $\mu
_j=0$, $1\leq j\leq r$ we pass from equivariant to non equivariant
cohomology, from equivariant to non equivariant integral.

In the proof of Theorem \ref{four} we write down all formulas for $r=1$ (it
serves the case when $X$ is a hypersurface in ${\Bbb P}^n$ of degree $l<n)$.
The proof for $r>1$ differs only for longer product formulas. So our aim is
to prove Theorem \ref{four} in the following form

\begin{equation}
S^{\prime }(t,h)=e^{\frac{pt}\hbar }\sum\limits_{d\geq 0}e^{dt}\frac{%
\prod_{m=0}^{dl}(lp-\mu +m\hbar )}{\prod_{m=1}^d(p-\lambda _0+m\hbar )\ldots
\prod_{m=1}^d(p-\lambda _n+m\hbar )}\text{.}  \label{five}
\end{equation}

From now on, with abuse of notation, we shall consider the ring $%
H_{T^{n+1}}^{*}({\Bbb P}^n{\bf )}$ tensorized with ${\Bbb Q}(\lambda
_0,\ldots ,\lambda _n)$; this is due to the fact that we need to localize on
fixed points components. It becomes a vector space on this field with basis $%
\phi _i=\frac{\prod_{j\neq i}(p-\lambda _j)}{\prod_{j\neq i}(\lambda
_i-\lambda _j)}$, $0\leq i\leq n$. Moreover 
\begin{equation*}
\left\langle \phi _i,\phi _j\right\rangle =\frac{\delta _{ij}}{\prod_{j\neq
i}(\lambda _i-\lambda _j)}\text{.} 
\end{equation*}

\noindent We recall that $W_{2,d}$ is the vector bundle on $\overline{{\cal M}}_{0,2}(%
{\Bbb P}^n,d)$ whose fiber over the point $\left[ C,x_0,x_1;f\right] $ is $%
H^0(C,f^{*}({\cal O}(l)))$ and ${\cal E}_{2,d}$ its equivariant Euler class,
with respect to the given linearization of the action. Besides, we define $%
W_{2,d}^{^{\prime }}$ via the following exact sequence of vector bundles

\begin{equation*}
0\longrightarrow W_{2,d}^{^{\prime }}\longrightarrow W_{2,d}\longrightarrow
(e_0)^{*}({\cal O}(l))\longrightarrow 0 
\end{equation*}
\noindent and ${\cal E}_{2,d}^{\prime }$ its Euler class.

\medskip\ 

\begin{lemma}
\label{six}{\bf \ }For{\sl \ }$0\leq i\leq n$
\end{lemma}

\begin{equation*}
\prod\limits_{j\neq i}\left( \lambda _i-\lambda _j\right) \left\langle \phi
_i,S^{^{\prime }}\right\rangle =e^{\frac{\lambda _it}\hbar }(l\lambda _i-\mu
)\sum\limits_{d\geq 0}e^{dt}\int_{\overline{{\cal M}}_{0,2}({\Bbb P}%
^n,d)}(e_0)^{*}(\phi _i)\frac{{\cal E}_{2,d}^{\prime }}{\hbar +c}\text{.} 
\end{equation*}

{\bf Proof. }Obviously, ${\cal E}_{2,d}={\cal E}_{2,d}^{\prime
}(e_0)^{*}(lp-\mu )$. Since $(e_0)_{*}(\frac{{\cal E}_d^{\prime }}{\hbar +c}%
)\in H_{T^{n+1}}^{*}({\Bbb P}^n{\bf )}$ may be written as $%
\sum_{k=0}^nf_k\phi _k$, for $f\in {\Bbb Q}\left[ p,\lambda _0,\ldots
,\lambda _n\right] $, we have 
\begin{equation*}
\begin{array}{ll}
\prod\limits_{j\neq i}\left( \lambda _i-\lambda _j\right) \left\langle \phi
_i,S^{\prime }\right\rangle & =\dfrac 1{2\pi \sqrt{-1}}\displaystyle%
\sum\limits_{d\geq 0}e^{dt}\displaystyle\int \dfrac{e^{\frac{pt}\hbar }\phi
_i(e_0)_{*}\left( \dfrac{{\cal E}_{2,d}^{\prime }(e_0)^{*}(lp-\mu )}{\hbar +c%
}\right) }{(p-\lambda _0)\ldots (p-\lambda _n)}dp= \\ 
& =\dfrac 1{2\pi \sqrt{-1}}\displaystyle\sum\limits_{d\geq 0}e^{dt}%
\displaystyle\int \dfrac{e^{\frac{pt}\hbar }\phi _i(lp-\mu )(e_0)_{*}\left( 
\dfrac{{\cal E}_{2,d}^{\prime }}{\hbar +c}\right) }{(p-\lambda _0)\ldots
(p-\lambda _n)}dp= \\ 
& =\dfrac 1{2\pi \sqrt{-1}}\displaystyle\sum\limits_{k=0}^nf_k\phi _{k\mid
p=\lambda _i}\sum\limits_{d\geq 0}e^{dt}e^{\frac{\lambda _it}\hbar
}(l\lambda _i-\mu )= \\ 
& =e^{\frac{\lambda _it}\hbar }(l\lambda _i-\mu )\displaystyle%
\sum\limits_{d\geq 0}e^{dt}\dfrac 1{2\pi \sqrt{-1}}\displaystyle\int \dfrac{%
\phi _i(e_0)_{*}\left( \dfrac{{\cal E}_{2,d}^{\prime }}{\hbar +c}\right) }{%
(p-\lambda _0)\ldots (p-\lambda _n)}dp= \\ 
& =e^{\frac{\lambda _it}\hbar }(l\lambda _i-\mu )\displaystyle%
\sum\limits_{d\geq 0}e^{dt}\displaystyle\int_{\overline{{\cal M}}_{0,2}(%
{\Bbb P}^n,d)}(e_0)^{*}(\phi _i)\dfrac{{\cal E}_{2,d}^{\prime }}{\hbar +c}%
\text{.}
\end{array}
\end{equation*}
\begin{flushright}
$\Box$
\end{flushright}

\medskip\ 

For the computation made in Lemma \ref{six}, we introduce functions 
\begin{equation*}
Z_i(q,\hbar ):=1+\sum\limits_{d>0}q^d\prod\limits_{j\neq i}(\lambda
_i-\lambda _j)\displaystyle\int_{\overline{{\cal M}}_{0,2}({\Bbb P}%
^n,d)}\left( e_0\right) ^{*}(\phi _i)\frac{{\cal E}_{2,d}^{\prime }}{\hbar +c%
},0\leq i\leq n. \label{for} 
\end{equation*}
\label{zed}

\noindent We will determine explicitly these functions by finding '{\sl linear
recursive relations}' satisfied by them. These relations stem from the
computation of the equivariant integral in \ref{zed} by means of the formula
of integration over connected components of fixed points in $\overline{{\cal %
M}}_{0,2}({\Bbb P}^n,d)$.

\noindent We restrict only to those components of elements $[C,x_0,x_1;f]$ whose
marked point $x_0$ is mapped to the $i$-th point in ${\Bbb P}^n$: indeed,
when computing the integral

\begin{equation}
\int\limits_{\overline{{\cal M}}_{0,2}({\Bbb P}^n,d)}\left( e_0\right)
^{*}(\phi _i)\frac{{\cal E}_{2,d}^{\prime }}{\hbar +c}\text{,}
\label{integral}
\end{equation}
$\phi _i$ has zero localizations at all the fixed points $p_j$, $j\neq i$.

\medskip\ 

\begin{proposition}
\label{zetai} 
\[
Z_i(q,\hbar ):=1+\sum\limits_{d>0}\left( \frac q{\hbar ^{n+1-l}}\right)
^d\prod\limits_{j\neq i}(\lambda _i-\lambda _j)\int\limits_{\overline{{\cal M%
}}_{0,2}({\Bbb P}^n,d)}\left( e_0\right) ^{*}(\phi _i)\frac{{\cal E}%
_{2,d}^{\prime }}{1+\frac c\hbar }\left( -c\right) ^{(n+1-l)d-1}\text{.}
\]
\end{proposition}

{\bf Proof. } Since $\frac 1{c+\hbar }=\sum\limits_j\frac{\left( -c\right) ^j%
}{\hbar ^{j+1}}$ and dim $\overline{{\cal M}}_{0,2}({\Bbb P}^n,d)=n+(n+1)d-1$%
, the equivariant integral 
\begin{equation*}
\int\limits_{\overline{{\cal M}}_{0,2}({\Bbb P}^n,d)}\left( e_0\right)
^{*}(\phi _i){\cal E}_{2,d}^{\prime }\left( -c\right) ^j 
\end{equation*}

\noindent gives non-zero contribution only if

\begin{equation*}
n+(n+1)d-1\leq j+n+\text{deg}({\cal E}_{2,d}^{\prime })=j+n+ld-1 
\end{equation*}
or 
\begin{equation*}
j\geq (n+1-l)d-1\text{.} 
\end{equation*}
\begin{flushright}
$\Box$
\end{flushright}

\medskip\ 

Consider now a fixed point curve $C$ whose marked point $x_0$ is indeed
mapped to the $i$-th fixed point in ${\Bbb P}^n$. There are two cases:

(i) the marked point $x_0$ is situated on an irreducible component of $C$
mapped with some degree $d^{\prime }$, $d^{\prime }\geq 1$ onto the line
joining the $i$-th fixed point with the $j$-th fixed point in ${\Bbb P}^n$
with $i\neq j$;

(ii) the marked point $x_0$ is situated on a component of $C$ mapped to the $%
i$-th fixed point and carrying two or more other special points.

\begin{proposition}
\label{type2}A type (ii) fixed point component ${\cal C}$ in $\overline{%
{\cal M}}_{0,2}({\Bbb P}^n,d)$ gives zero contribution to the integration
formula for

\[
\int\limits_{\overline{{\cal M}}_{0,2}({\Bbb P}^n,d)}\left( e_0\right)
^{*}(\phi _i){\cal E}_{2,d}^{\prime }\left( -c\right) ^k\text{, with }k\geq
(n+1-l)d-1\text{.\label{nine}}
\]
\end{proposition}

{\bf Proof. }Let $\Gamma $ be a graph associated to ${\cal C}$ and denote
by $v_0$ the vertex having label $\{0\}$. Since ${\cal C}$ is a type (ii)
component, $h_{v_0}+l_{v_0}-3\geq 0$. By Proposition \ref{graph}, there
exists an isomorphism $\xi _\Gamma $ between ${\cal C}$ and $\prod\limits_{%
\underset{h_v+l_v-3\geq 0}{v\in V_\Gamma }}\frac{\overline{{\cal M}}%
_{0,h_v+l_v}}{Aut(\Gamma )}$, with $Aut(\Gamma )$ the group of automorphism
of $\Gamma $. Consider the projection

\begin{equation*}
\pi _{v_0}:\prod\limits_{\underset{h_v+l_v-3\geq 0}{v\in V_\Gamma }}\frac{%
\overline{{\cal M}}_{0,h_v+l_v}}{Aut(\Gamma )}\longrightarrow \frac{%
\overline{{\cal M}}_{0,h_{v_0}+l_{v_0}}}{Aut(\Gamma )} 
\end{equation*}
and denote by $c^{\prime }$ the Chern class of the $Aut(\Gamma )$-invariant
line bundle ${\cal L}$ on $\overline{{\cal M}}_{0,h_{v_0}+l_{v_0}}$ whose
fiber at the point $[C,x_{v_0,1}=x_0,\ldots ,x_{h_{v_0}+l_{v_0}}]$ is the
tangent space to $C$ at $x_0$. By the definiton of the morphism $\xi _\Gamma 
$, we see that $\pi _{v_0}^{*}(c^{\prime })=\xi _\Gamma ^{*}(c)$. Thus if we
compute \ref{integral} over the connected component ${\cal C}$, we get

\begin{eqnarray*}
\int\limits_{{\cal C}}\frac{\left( e_0\right) ^{*}(\phi _i){\cal E}%
_{2,d}^{\prime }}{{\cal E}({\cal N}_{{\cal C}})}(-c)^k &=&\int_{\prod_v%
\overline{{\cal M}}_{0,h_v+l_v}}\xi _\Gamma ^{*}\left( (-c)^k\left( \frac{%
\left( e_0\right) ^{*}(\phi _i){\cal E}_{2,d}^{\prime }}{{\cal E}({\cal N}_{%
{\cal C}})}\right) \right) = \\
&=&\frac 1{\left| Aut(\Gamma )\right| }\int_{\overline{{\cal M}}%
_{0,h_{v_0}+l_{v_0}}}(-c^{\prime })^k(\pi _{v_0})_{*}\left( \xi _\Gamma
^{*}\left( \frac{\left( e_0\right) ^{*}(\phi _i){\cal E}_{2,d}^{\prime }}{%
{\cal E}({\cal N}_{{\cal C}})}\right) \right)
\end{eqnarray*}
with ${\cal E}({\cal N}_{{\cal C}})$ the Euler class of the normal bundle $%
{\cal N}_{\overline{{\cal M}}_{0,2}({\Bbb P}^n,d)/{\cal C}}$. On the other
hand, dim$\overline{{\cal M}}_{0,h_{v_0}+l_{v_0}}=$ $%
h_{v_0}+l_{v_0}-3<(n+1-l)d-1\leq k$, since $l<n$ and $l_{v_0}\leq d$. Since
the nilpotency degree of $c^{\prime }$ should not exceed $h_{v_0}+l_{v_0}-3$%
, the Lemma follows.
\begin{flushright}
$\Box$
\end{flushright}

\medskip\ 

\begin{figure}[htb]
\hspace{-2cm} 
\begin{center} 
\mbox{\epsfig{file=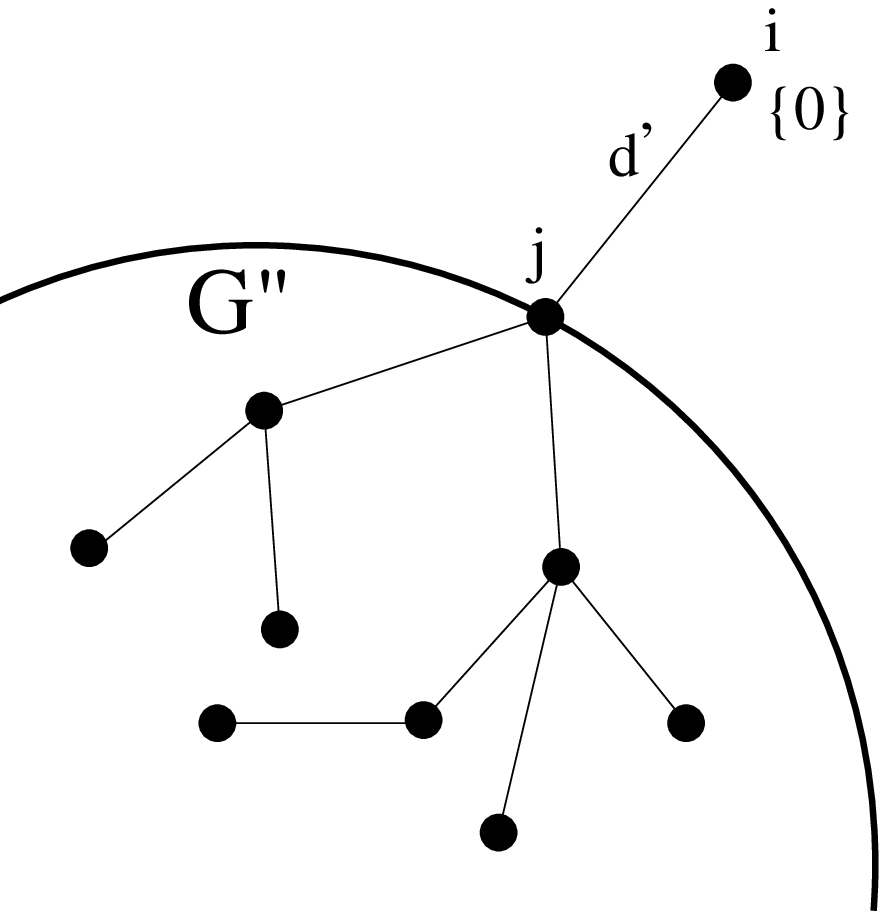,width=6cm,height=6cm,angle=270}}
\caption[]{
\label{dis4}}
\end{center}
\end{figure}

Consider now a fixed point $[C,x_0,x_1;f]$ in a type (i) component ${\cal K}$%
. By assumption, $C$ can be regarded as the union of a connected curve $%
C^{\prime \prime }$ and an irreducible component $C^{\prime }$, which
carries the marked point $x_0$ and maps onto $l_{ij}$, $i\neq j$, with
degree $d^{\prime }$, $1\leq d^{\prime }\leq d$. Moreover, if $f^{\prime
\prime }$ denotes the restriction of $f$ to $C^{\prime \prime }$, then $%
[C^{\prime \prime },y,x_1;f^{\prime \prime }]$ ($y=C^{\prime \prime }\cap
C^{\prime }$) is obviously a fixed point in $\overline{{\cal M}}_{0,2}({\Bbb %
P}^n,d-d^{\prime })$. (Observe that if $d=d^{\prime }$, then $C^{\prime
\prime }=\emptyset $.) In addition, the graph $G$ associated with ${\cal K}$
contains a subgraph $G_{d^{\prime }}^{\prime }$ with one edge $\alpha $
labelled $d^{\prime }$ and two vertices $v_1$, $v_2$ having respectively $j$
and $\left\{ \{0\},i\right\} $ as labels. If we remove from $G$ the edge $%
\alpha $ and the vertex $v_2$, and add to $v_1$ the label $\{0\}$, we obtain
a graph $G^{\prime \prime }$ which corresponds to the connected component in 
$\overline{{\cal M}}_{0,2}({\Bbb P}^n,d-d^{\prime })$ containing $[C^{\prime
\prime },y,x_1;f^{\prime \prime }]$. Notice that $G^{\prime \prime }$ is
empty when $d^{\prime }=d$. In this case ${\cal K}$ is a point in $\overline{%
{\cal M}}_{0,2}({\Bbb P}^n,d)$. See figure \ref{dis4}.

The fundamental idea is to reduce ourselves to integrate on connected
component corresponding to ``smaller'' graphs, i.e., in this case, to write
such a formula 
\begin{equation*}
\int_{{\cal K}(G)}...=\omega \left( d^{\prime }\right) \int_{{\cal K}\left(
G^{\prime \prime }\right) }... 
\end{equation*}
and from this to build recursive relations among the $Z_i$. It is easy to check that 
the automorphism groups of the graphs behave well with respect to this decomposition.
Let us first consider the easier case, namely $d=d^{\prime }$.

\begin{lemma}
\label{ten}Suppose $G^{\prime \prime }=\emptyset $ and denote by ${\cal N}_{%
{\cal K}}$ the normal bundle ${\cal N}_{\overline{{\cal M}}_{0,2}({\Bbb P}%
^n,d)/{\cal K}}$. Then
\end{lemma}

\begin{equation*}
\prod\limits_{j\neq i}(\lambda _i-\lambda _j)\int_{{\cal K}}\frac{%
e_0^{*}(\phi _i){\cal E}_{2,d}^{\prime }(-c)^{(n+1-l)d-1}}{(1+\frac c\hbar )%
{\cal E}({\cal N}_{{\cal K}})}= 
\end{equation*}

\begin{equation*}
=\frac{\prod\limits_{m=1}^{ld}[l\lambda _i-\mu +\frac md(\lambda _j-\lambda
_i)]\left( \frac{\lambda _j-\lambda _i}d\right) ^{(n+1-l)d-1}}{d\left( 1+%
\frac{\lambda _i-\lambda j}{d\hbar }\right) \prod\limits_{\underset{(\alpha
,m)\neq (j,d)}{\alpha =0}}^n\prod\limits_{m=1}^d[\lambda _i-\lambda _\alpha
+\frac md(\lambda _j-\lambda _i)]}:=coeff_i^j(d)\text{.} 
\end{equation*}

{\bf Proof.} Set ${\cal K}=\{[C^{\prime },x_0^{\prime },x_1^{\prime
};f^{\prime }]\}$. We will prove the Lemma in several steps.

{\bf STEP 1)} The localization of $e_0^{*}(\phi _i)$ at ${\cal K}$ is $1$.
This follows easily from the following commutative diagram: 
\begin{equation*}
\begin{array}{ccc}
e_0^{*}:H_{T^{n+1}}^{*}({\Bbb P}^n) & \longrightarrow & H_{T^{n+1}}^{*}(%
\overline{{\cal M}}_{0,2}({\Bbb P}^n,d)) \\ 
\delta _{p_i}^{*}\downarrow &  & \downarrow \delta _{{\cal K}}^{*} \\ 
H_{T^{n+1}}^{*}(\{p_i\}) & \longrightarrow & H_{T^{n+1}}^{*}(\{{\cal K}\})
\end{array}
\end{equation*}
with $\delta _{p_i}:\{p_i\}\rightarrow {\Bbb P}^n$ and $\delta _{{\cal K}}:\{%
{\cal K}\}\rightarrow \overline{{\cal M}}_{0,2}({\Bbb P}^n,d)$.

\smallskip\ 

{\bf STEP 2)} The localization of $c$ at ${\cal K}$ is $\frac{\lambda
_i-\lambda _j}d$.

\noindent Let $T^{(0)}$ be the equivariant line bundle whose fiber at $[C,x_0,x_1;f]$
is the tangent line to $C$ at $x_0$. With the same notation adopted in {\bf %
STEP 1)} we need to compute $c_1(\delta _{{\cal K}}^{*}T^{(0)})$, namely the
Chern class of the line bundle

\begin{equation*}
\begin{array}{c}
T_{x_0^{^{\prime }}}C^{\prime }\times _{T^{n+1}}(S^\infty )^{n+1} \\ 
\downarrow \\ 
({\Bbb P}^\infty )^{n+1}
\end{array}
\text{,} 
\end{equation*}
with $T_{x_0^{^{\prime }}}C^{\prime }$ the tangent line to $C^{\prime }$ at $%
x_0^{^{\prime }}$. On the other hand, if $z_1,z_2$ denote homogeneous
coordinates on $C^{\prime }$, the following diagram holds:

\begin{equation*}
\begin{array}{ccc}
\lbrack z_1,z_2] & \longrightarrow & [t_i^{\frac 1d}z_1,t_j^{\frac 1d}z_2]
\\ 
\downarrow &  & \downarrow \\ 
\lbrack 0,\ldots ,z_1^d,z_2^d,\ldots ,0] & \longrightarrow & [0,\ldots
,t_iz_1^d,t_jz_2^d,\ldots ,0]\text{.}
\end{array}
\end{equation*}

\noindent At the point $x_0=[1,0]$, the action on $T_{x_0^{^{\prime }}}C^{\prime }$ is
given by

\begin{equation*}
\frac \partial {\partial z}\mapsto (t_jt_i^{-1})^{\frac 1d}\frac \partial
{\partial z}\text{, }z=\frac{z_2}{z_1}\text{.} 
\end{equation*}

\smallskip\ 

{\bf STEP 3)} The localization of ${\cal E}_{2,d}^{\prime }$ at ${\cal K}$
is $\prod\limits_{m=1}^{ld}[(l\lambda _i-\mu )+\frac md(\lambda _j-\lambda
_i)]$.

\noindent By the same reasoning adopted in {\bf STEP\ 2)} we need to compute the Euler
class of the vector bundle 
\begin{equation*}
\begin{array}{c}
W_{2,d\mid _{[C^{\prime },x_0^{\prime },x_1^{\prime };f^{\prime }]}}^{\prime
}\times _{T^{n+1}}(S^\infty )^{n+1}\times _{S^1}(S^\infty ) \\ 
\downarrow \\ 
({\Bbb P}^\infty )^{n+2}\text{.}
\end{array}
\end{equation*}

\noindent By definition of $W_{2,d\mid _{[C^{\prime },x_0^{\prime },x_1^{\prime
};f^{\prime }]}}^{\prime }$, a basis for this fiber is clearly $%
\{z_2^mz_1^{dl-m}\}_{k=1}^{dl}$. Since the group $T^{\prime }=S^1$ acts on
global sections of $f^{\prime *}({\cal O}(l))$ by multiplication by $%
t^{^{\prime }-1}$, the action on the basis is

\begin{equation*}
z_2^mz_1^{dl-m}\longrightarrow t^{\prime -1}t_i^{\left( l-\frac md\right)
}t_j^{\frac md}z_2^mz_1^{dl-m}\text{, }m\in \{1,\ldots ,ld\}\text{.} 
\end{equation*}
Hence the Euler class is

\begin{equation*}
\prod\limits_{m=1}^{ld}[(l\lambda _i-\mu )+\frac md(\lambda _j-\lambda _i)]%
\text{.} 
\end{equation*}

\smallskip\ 

{\bf STEP 4)} The localization of ${\cal E}({\cal N}_{{\cal K}})$ at ${\cal K%
}$ is 
$$
\prod\limits_{j\neq i}(\lambda _i-\lambda _j)\prod\limits_{\underset{%
(\alpha ,m)\neq (j,d)}{\alpha =0}}^n\prod\limits_{m=1}^d\left( \lambda
_i-\lambda _\alpha +\frac md(\lambda _j-\lambda _i)\right) \text{.} 
$$

\noindent Remind that ${\cal K}$ is a point and thus the normal bundle is the tangent
space to $\overline{{\cal M}}_{0,2}({\Bbb P}^n,d)$ at ${\cal K}$. We next
define the vector bundle ${\cal W}_{f^{\prime }}$ via the exact sequence

\begin{equation*}
0\rightarrow T_{C^{\prime }}\otimes {\cal O}(-x_0^{\prime }-x_1^{\prime
})\rightarrow f^{\prime *}T_{{\Bbb P}^n}\rightarrow {\cal W}_{f^{\prime
}}\rightarrow 0\text{.} 
\end{equation*}
\label{norm1}

\noindent By Remark \ref{norbu}, the tangent space to $\overline{{\cal M}}_{0,2}({\Bbb %
P}^n,d)$ at ${\cal K}$ is exactly $H^0(C^{\prime },{\cal W}_{f^{\prime }})$.
Hence,

\begin{equation*}
{\cal E}({\cal N}_{{\cal K}})=\frac{{\cal E}(H^0(C^{\prime },f^{\prime *}T_{%
{\Bbb P}^n}))}{{\cal E}(H^0(C^{\prime },T_{C^{\prime }}\otimes {\cal O}%
_{C^{\prime }}(-x_0^{\prime }-x_1^{\prime }))}\text{.} 
\end{equation*}

\noindent We need to compute a basis for $H^0(C^{\prime },f^{\prime *}T_{{\Bbb P}^n})$%
. Since $f^{\prime }:C^{\prime }\rightarrow l_{ij}\subseteq {\Bbb P}^n$, the
following exact sequence holds

\begin{equation*}
0\rightarrow T_{l_{ij}}\rightarrow T_{{\Bbb P}^n\mid _{l_{ij}}}\rightarrow 
{\cal N}_{l_{ij}}\rightarrow 0\text{.} 
\end{equation*}
As ${\Bbb P}^n$ is convex, this implies

\begin{equation*}
0\rightarrow H^0(C^{\prime },f^{^{\prime }*}T_{l_{ij}})\rightarrow
H^0(C^{\prime },f^{\prime *}T_{{\Bbb P}^n\mid _{l_{ij}}})\rightarrow
H^0(C^{\prime },f^{\prime *}{\cal N}_{l_{ij}})\rightarrow 0\text{.} 
\end{equation*}
A basis for $H^0(C^{\prime },f^{\prime *}{\cal N}_{l_{ij}})$ is given by $%
z_1^cz_2^mf^{\prime *}\left( \frac \partial {\partial X_\alpha }\right) $,
with $\alpha \in \{0,\ldots \widehat{i},\widehat{j},\ldots ,n\}$ and $c+m=d$%
, $m,c\geq 0$, since $f^{\prime *}{\cal N}_{l_{ij}}=\overset{(n-1)-times}{%
\overbrace{{\cal O}_{C^{^{\prime }}}(d)\oplus \ldots \oplus {\cal O}%
_{C^{^{\prime }}}(d)}}$.

\noindent On the other hand, a basis for $H^0(C^{\prime },f^{^{\prime }*}T_{l_{ij}})$
is $\left\{ \left( \frac{z_2}{z_1}\right) ^az_2\frac \partial {\partial
z_2}\right\} _{a=-d}^d$, since $f^{*}T_{l_{ij}}=\frac{{\cal O}_{C^{^{\prime
}}}(d)\oplus {\cal O}_{C^{^{\prime }}}(d)}{{\Bbb C}}$.

\noindent Finally, we observe that a basis for $H^0(C^{\prime },T_{C^{\prime
}}\otimes {\cal O}_{C^{\prime }}(-x_0^{\prime }-x_1^{\prime })$ is $%
z_2\frac \partial {\partial z_2}$. At this point, the action of the torus $%
T^{n+1}$ on $H^0(C^{\prime },f^{\prime *}T_{{\Bbb P}^n\mid _{l_{ij}}})$ is

\begin{equation*}
z_1^cz_2^mf^{\prime *}\left( \frac \partial {\partial X_\alpha }\right)
\rightarrow t_i^{\frac ad}t_j^{\frac md}t_k^{-1}z_1^cz_2^mf^{\prime *}\left(
\frac \partial {\partial X_\alpha }\right) \text{, }\alpha \in \{0,\ldots 
\widehat{i},\widehat{j},\ldots ,n\} 
\end{equation*}
\noindent and $c+m=d,\text{ }c,m\geq 0$,
\begin{equation*}
\left( \frac{z_2}{z_1}\right) ^az_2\frac \partial {\partial z_2}\rightarrow
t_j^{\frac ad}t_i^{-\frac ad}t_j^{\frac ad}t_j^{-\frac ad}\left( \frac{z_2}{%
z_1}\right) ^az_2\frac \partial {\partial z_2}\text{.} 
\end{equation*}

\noindent In other words, we have

\begin{equation*}
{\cal E}({\cal N}_{{\cal K}})=\prod\limits_{\underset{a\neq 0}{a=-d}%
}^d\left( \frac{a(\lambda _j-\lambda _i)}d\right) \prod\limits_{\underset{%
\alpha \neq j}{\alpha \neq i}}^n\prod_{\underset{m,c\geq 0}{c+m=d}}[\frac
cd\lambda _i+\frac md\lambda _j-\lambda _\alpha ]= 
\end{equation*}

\begin{equation*}
=\prod\limits_{\alpha \neq \{i,j\}}(\lambda _i-\lambda _\alpha
)\prod\limits_{\underset{a\neq 0}{a=-d}}^d\left( \frac{a(\lambda _j-\lambda
_i)}d\right) \prod\limits_{\underset{\alpha \neq j}{\alpha \neq i}%
}^n\prod_{m=1}^d[\lambda _i-\lambda _\alpha +\frac md(\lambda _j-\lambda
_i)]= 
\end{equation*}

\begin{equation*}
=\prod\limits_{\alpha \neq \{i,j\}}(\lambda _i-\lambda _\alpha
)\prod\limits_{\underset{a\neq 0}{a=-d}}^{-1}\left( \frac{a(\lambda
_j-\lambda _i)}d\right) \prod\limits_{\alpha \neq j}^n\prod_{m=1}^d[\lambda
_i-\lambda _\alpha +\frac md(\lambda _j-\lambda _i)]= 
\end{equation*}

\begin{equation*}
=\prod\limits_{\alpha \neq \{i,j\}}(\lambda _i-\lambda _\alpha )(\lambda
_i-\lambda _j)\prod\limits_{\underset{(\alpha ,m)\neq (j,d)}{\alpha =0}%
}^n\prod_{m=1}^d[\lambda _i-\lambda _\alpha +\frac md(\lambda _j-\lambda
_i)]= 
\end{equation*}

\begin{equation*}
=\prod\limits_{j\neq i}(\lambda _i-\lambda _j)\prod\limits_{\underset{%
(\alpha ,m)\neq (j,d)}{\alpha =0}}^n\prod\limits_{m=1}^d\left( \lambda
_i-\lambda _\alpha +\frac md(\lambda _j-\lambda _i)\right) \text{.} 
\end{equation*}

\smallskip\ 

To summarize, we have

\begin{equation*}
\prod\limits_{j\neq i}(\lambda _i-\lambda _j)\int_{{\cal K}}\frac{%
e_0^{*}(\phi _i){\cal E}_{2,d}^{\prime }(-c)^{(n+1-l)d-1}}{(1+\frac c\hbar )%
{\cal E}({\cal N}_{{\cal K}})}= 
\end{equation*}

\begin{equation*}
=\frac{\prod\limits_{m=1}^{ld}[l\lambda _i-\mu +\frac md(\lambda _j-\lambda
_i)]\left( \frac{\lambda _j-\lambda _i}d\right) ^{(n+1-l)d-1}}{\left( 1+%
\frac{\lambda _i-\lambda j}{d\hbar }\right) \prod\limits_{\underset{(\alpha
,m)\neq (j,d)}{\alpha =0}}^n\prod\limits_{m=1}^d[\lambda _i-\lambda _\alpha
+\frac md(\lambda _j-\lambda _i)]}\int\limits_{{\cal K}}1, 
\end{equation*}
and by the orbifold structure of $\overline{{\cal M}}_{0,2}({\Bbb P}^n,d)$
this equals

\begin{equation*}
\frac{\prod\limits_{m=1}^{ld}[l\lambda _i-\mu +\frac md(\lambda _j-\lambda
_i)]\left( \frac{\lambda _j-\lambda _i}d\right) ^{(n+1-l)d-1}}{\left( 1+%
\frac{\lambda _i-\lambda j}{d\hbar }\right) \prod\limits_{\underset{(\alpha
,m)\neq (j,d)}{\alpha =0}}^n\prod\limits_{m=1}^d[\lambda _i-\lambda _\alpha
+\frac md(\lambda _j-\lambda _i)]}\frac 1{\#Aut(f^{^{\prime }}:C^{\prime
}\rightarrow l_{ij})}= 
\end{equation*}

\begin{equation*}
=\frac{\prod\limits_{m=1}^{ld}[l\lambda _i-\mu +\frac md(\lambda _j-\lambda
_i)]\left( \frac{\lambda _j-\lambda _i}d\right) ^{(n+1-l)d-1}}{\left( 1+%
\frac{\lambda _i-\lambda j}{d\hbar }\right) \prod\limits_{\underset{(\alpha
,m)\neq (j,d)}{\alpha =0}}^n\prod\limits_{m=1}^d[\lambda _i-\lambda _\alpha
+\frac md(\lambda _j-\lambda _i)]}\frac 1d\text{.} 
\end{equation*}
So the Lemma is proved.
\begin{flushright}
$\Box$
\end{flushright}

\medskip\ 

Consider now a connected component whose associated graph is the union of
two subgraphs as described before Lemma \ref{ten}. Fix $d^{\prime }$, $1\leq
d^{\prime }\leq d$ and $j$ (the label of the vertex $v_1$). Denote by ${\cal %
K}_{G^{\prime \prime },j}$ the connected component in $\overline{{\cal M}}%
_{0,2}({\Bbb P}^n,d-d^{\prime })$ associated with $G^{\prime \prime }$. We
can now prove

\begin{lemma}
\label{eleven}For a connected component ${\cal K}$ as above,
\end{lemma}

\begin{equation*}
\prod\limits_{j\neq i}(\lambda _i-\lambda _j)\int_{{\cal K}}\frac{%
e_0^{*}(\phi _i){\cal E}_{2,d}^{\prime }(-c)^{(n+1-l)d-1}}{(1+\frac c\hbar )%
{\cal E}({\cal N}_{{\cal K}})}= 
\end{equation*}

\begin{equation*}
=\frac{\prod\limits_{m=1}^{ld^{\prime }}[l\lambda _i-\mu +\frac m{d^{\prime
}}(\lambda _j-\lambda _i)]\left( \frac{\lambda _j-\lambda _i}{d^{\prime }}%
\right) ^{(n+1-l)d-1}}{d^{\prime }\left( 1+\frac{\lambda _i-\lambda j}{%
d^{\prime }\hbar }\right) \prod\limits_{\underset{(\alpha ,m)\neq
(j,d^{\prime })}{\alpha =0}}^n\prod\limits_{m=1}^{d^{\prime }}[\lambda
_i-\lambda _\alpha +\frac m{d^{\prime }}(\lambda _j-\lambda _i)]}%
\cdot 
\end{equation*}
\begin{equation*}
\cdot
\prod\limits_{k\neq j}(\lambda _j-\lambda _k)\int\limits_{{\cal K}%
_{G^{\prime \prime },j}}\frac{{\cal E}_{2,d-d^{\prime }}^{\prime
}e_0^{*}(\phi _j)}{{\cal E}({\cal N}_{{\cal K}_{G^{\prime \prime },j}})(%
\frac{\lambda _j-\lambda _i}{d^{\prime }}+c)}= 
\end{equation*}

\begin{equation*}
=\left( \frac{\lambda _j-\lambda _i}{d^{\prime }}\right)
^{(n+1-l)(d-d^{\prime })}coeff_i^j(d^{\prime })\prod\limits_{k\neq
j}(\lambda _j-\lambda _k)\int\limits_{{\cal K}_{G^{\prime \prime },j}}\frac{%
{\cal E}_{2,d-d^{\prime }}^{\prime }e_0^{*}(\phi _j)}{{\cal E}({\cal N}_{%
{\cal K}_{G^{\prime \prime },j}})(\frac{\lambda _j-\lambda _i}{d^{\prime }}%
+c)}\text{.} 
\end{equation*}

{\bf Proof.} By our assumptions, we can regard ${\cal K}$ in the following
way. Consider $\overline{{\cal M}}_{0,2}({\Bbb P}^n,d^{\prime })$ and $%
\overline{{\cal M}}_{0,2}({\Bbb P}^n,d-d^{\prime })$, and denote by $%
e_0^{^{\prime \prime }}$ and $e_1^{\prime }$ the evaluation maps

\begin{equation*}
e_1^{\prime }:\overline{{\cal M}}_{0,2}({\Bbb P}^n,d^{\prime
})\longrightarrow {\Bbb P}^n\text{,} 
\end{equation*}

\begin{equation*}
e_0^{^{\prime \prime }}:\overline{{\cal M}}_{0,2}({\Bbb P}^n,d-d^{\prime
})\longrightarrow {\Bbb P}^n\text{.} 
\end{equation*}

\begin{figure}[htb]
\hspace{-2cm} 
\begin{center} 
\mbox{\epsfig{file=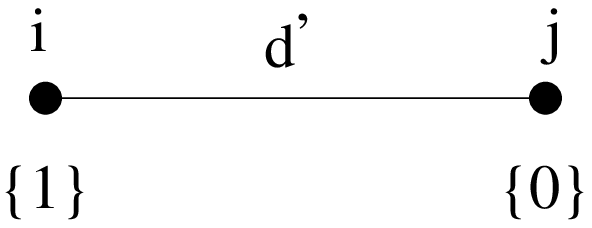,width=2cm,height=4cm,angle=270}}
\caption[]{
\label{dis3}}
\end{center}
\end{figure}

\noindent Next let ${\cal K}_{d^{\prime },j}$ the connected component of fixed points
in $\overline{{\cal M}}_{0,2}({\Bbb P}^n,d^{\prime })$ associated with the
graph of figure \ref{dis3}, and ${\cal K}_{G^{\prime \prime },j}$ the
connected component in $\overline{{\cal M}}_{0,2}({\Bbb P}^n,d-d^{\prime })$
associated with the graph $G^{\prime \prime }$. Then ${\cal K}\simeq {\cal K}%
_{d^{\prime },j}\times {\cal K}_{G^{\prime \prime },j}\subseteq \overline{%
{\cal M}}_{0,2}({\Bbb P}^n,d^{\prime })\times _{{\Bbb P}^n}\overline{{\cal M}%
}_{0,2}({\Bbb P}^n,d-d^{\prime })$, where the fiber product is defined via
the evaluation maps $e_0^{^{\prime \prime }}$ and $e_1^{\prime }$. Even in
this Lemma, we proceed in several steps.

{\bf STEP 1)} The localization of $e_0^{*}(\phi _i)$ at ${\cal K}$ is $1$.
This is similar to {\bf STEP 1)} of Lemma \ref{ten}.

\smallskip\ 

{\bf STEP 2)}${\Bbb \ }$Since the line bundle $T^{(0)}$ is trivial on ${\cal %
K}$, the localization of $c$ is computed in the same manner as in {\bf STEP
2)} of Lemma \ref{ten}.

\smallskip\ 

{\bf STEP 3)} Consider the projections

\begin{equation*}
\pi _1:{\cal K}\longrightarrow {\cal K}_{d^{\prime },j}\text{,} 
\end{equation*}

\begin{equation*}
\pi _2:{\cal K}\longrightarrow {\cal K}_{G^{\prime \prime },j}\text{,} 
\end{equation*}
and define the vector bundles

\begin{equation*}
{\cal F}_{2,d^{\prime }}:=\pi _1^{*}(W_{2,d^{^{\prime }}})\text{ and} 
\end{equation*}

\begin{equation*}
{\cal F}_{2,d-d^{\prime }}:=\pi _2^{*}(W_2,_{d-d^{\prime }})\text{.} 
\end{equation*}
If 
\begin{equation*}
\nu :{\cal K}\longrightarrow {\Bbb P}^n 
\end{equation*}
is the evaluation map on the point corresponding to the vertex $v_1$, then
we have the following exact sequence

\begin{equation*}
0\rightarrow W_{2,d}\rightarrow {\cal F}_{2,d^{\prime }}\oplus {\cal F}%
_{2,d-d^{\prime }}\rightarrow \nu ^{*}({\cal O}(l))\rightarrow 0\text{.} 
\end{equation*}
Moreover,

\begin{equation*}
0\rightarrow W_{2,d}^{\prime }\rightarrow W_{2,d}\rightarrow e_0^{*}({\cal O}%
(l))\rightarrow 0\text{,} 
\end{equation*}

\begin{equation*}
0\rightarrow W_{2,d^{\prime }}^{\prime }\rightarrow W_{2,d^{\prime
}}\rightarrow e_0^{\prime *}({\cal O}(l))\rightarrow 0\text{,} 
\end{equation*}

\begin{equation*}
0\rightarrow W_{2,d-d^{\prime }}^{\prime }\rightarrow W_{2,d-d^{\prime
}}\rightarrow e_0^{^{\prime \prime }*}({\cal O}(l))\rightarrow 0\text{.} 
\end{equation*}
Since $e_0^{\prime }\circ \pi _1=e_0$ and $e_0^{\prime \prime }\circ \pi
_2=\nu $, the Euler class of $W_{2,d}^{\prime }$ is 
\begin{equation*}
{\cal E}_{2,d}^{\prime }={\cal E}_{2,d^{\prime }}^{^{\prime }}{\cal E}%
_{2,d-d^{\prime }}^{\prime }\frac{{\cal E}(\pi _1^{*}e_0^{^{\prime }*}({\cal %
O}(l))){\cal E}(\pi _2^{*}e_0^{\prime \prime *}({\cal O}(l)))}{{\cal E}%
(e_0^{*}({\cal O}(l))){\cal E}(\nu ^{*}({\cal O}(l)))}={\cal E}_{2,d^{\prime
}}^{^{\prime }}{\cal E}_{2,d-d^{\prime }}^{\prime }\text{.} 
\end{equation*}
Observe that the localization of ${\cal E}_{2,d^{\prime }}^{\prime }$ is
computed in {\bf STEP 3)} of Lemma \ref{ten}.

\smallskip\ 

{\bf STEP 4) }By our hypotheses on elements $[C,x_0,x_1,f]$ of ${\cal K}$, $%
C=C^{\prime }\cup C^{\prime \prime }$, with $C^{\prime }$ mapped to $l_{ij}$
with degree $d^{\prime }$. Denote by $p$ (mapped to $p_j$ in ${\Bbb P}^n$)
the intersection point of $C^{\prime }$ and $C^{\prime \prime }$ and by $%
\left\{ q_1,\ldots ,q_s\right\} $ the remaining nodes of $C$. We define the
vector bundle ${\cal L}$ via the exact sequence

\begin{equation*}
0\rightarrow {\cal N}_{{\cal K}}\rightarrow \pi _1^{*}({\cal N}_{{\cal K}%
_{j,d^{\prime }}})\oplus \pi _2^{*}({\cal N}_{{\cal K}_{G^{\prime \prime
},j}})\rightarrow {\cal L}\rightarrow 0\text{.} 
\end{equation*}

\noindent We also have

\begin{equation*}
0\rightarrow H^0\left( C,f^{*}T_{{\Bbb P}^n}\right) \rightarrow
H^0(C^{\prime },f^{^{\prime }*}T_{{\Bbb P}^n})\oplus H^0(C^{\prime \prime
},f^{\prime \prime *}T_{{\Bbb P}^n})\rightarrow f^{*}T_{p_j}\rightarrow 0%
\text{,} 
\end{equation*}

\begin{equation*}
0\rightarrow H^0\left( C^{\prime },T_{C^{\prime }}(-p)\right) \rightarrow
H^0\left( C^{\prime },T_{C^{\prime }}\right) \rightarrow T_p\rightarrow 0%
\text{,} 
\end{equation*}
\begin{equation*}
0\rightarrow H^0\left( C^{\prime },T_{C^{\prime
}}(-p-\sum_{t=1}^sq_t)\right) \rightarrow H^0\left( C^{\prime \prime
},T_{C^{\prime \prime }}(-\sum_{t=1}^sq_t)\right) \rightarrow T_p\rightarrow
0\text{,} 
\end{equation*}
where $f^{\prime }=f_{\mid _{C^{\prime }}}$ and $f^{\prime \prime }=f_{\mid
_{C^{\prime \prime }}}$ and, by abuse of notation, we still denote by $p$
the points on $C^{\prime }$ and $C^{\prime \prime }$ which correspond to $%
p=C^{\prime }\cap C^{\prime \prime }$.

\noindent Moreover, the deformations which preserve the combinatorial type of $C$ are
exactly the deformations which preserve the combinatorial type of $C^{\prime
}$ and those of $C^{\prime \prime }$. Keeping in mind Remark \ref{norbu},
the Euler class of ${\cal L}$ is 
\begin{equation*}
\frac{{\cal E}(T_pC^{\prime }\otimes T_pC^{\prime \prime })}{{\cal E}%
(f^{*}T_{p_j})}\text{.} 
\end{equation*}

\noindent Let $c$ denote again the Euler class of $T_pC^{\prime \prime }$ and observe
that ${\cal E}(T_pC^{\prime })$ is given by computing the character of $%
T^{n+1}$ on the tangent space to $C^{\prime }$ at $p$, since $T_pC^{\prime }$
is trivial on ${\cal K}_{j,d^{\prime }}$. If we introduce homogeneous
coordinates $z_1,z_2$ on $C^{\prime }$, we remind that the action of $%
T^{n+1} $ is 
\begin{equation*}
\lbrack z_1,z_2]\longrightarrow [t_i^{\frac 1{d^{\prime }}}z_1,t_j^{\frac
1{d^{\prime }}}z_2]\text{ } 
\end{equation*}
and if $p$ $=[0,1]$, then we have 
\begin{equation*}
\frac{z_1}{z_2}\rightarrow t_i^{\frac 1{d^{\prime }}}t_j^{-\frac 1{d^{\prime
}}}\frac{z_1}{z_2}\text{.} 
\end{equation*}
Hence ${\cal E}(T_pC^{\prime })=\frac{\lambda _j-\lambda _i}{d^{\prime }}$.
Finally, ${\cal E}(f^{*}T_{p_j})=\prod\limits_{k\neq j}(\lambda _j-\lambda
_k)$, since a basis for $f^{*}T_{p_j}$ is given $\frac \partial {\partial
X_k}$, $k\neq j$.
We have come to the conclusion 
\begin{equation*}
{\cal E}({\cal N}_{{\cal K}})=\frac{{\cal E}({\cal N}_{{\cal K}_{G^{\prime
\prime },j}})\frac{\lambda _j-\lambda _i}{d^{\prime }}\prod\limits_{j\neq
i}(\lambda _i-\lambda _j)\prod\limits_{\underset{(\alpha ,m)\neq (j,d)}{%
\alpha =0}}^n\prod\limits_{m=1}^d\left( \lambda _i-\lambda _\alpha +\frac
md(\lambda _j-\lambda _i)\right) }{\prod\limits_{k\neq j}(\lambda _j-\lambda
_k)}\text{,} 
\end{equation*}
because of {\bf STEP 4)} of Lemma \ref{ten}.

\smallskip\ 

To summarize, we have

\begin{equation*}
\prod\limits_{j\neq i}(\lambda _i-\lambda _j)\int_{{\cal K}}\frac{%
e_0^{*}(\phi _i){\cal E}_{2,d}^{\prime }(-c)^{(n+1-l)d-1}}{(1+\frac c\hbar )%
{\cal E}({\cal N}_{{\cal K}})}= 
\end{equation*}

\begin{equation*}
=\frac{\prod\limits_{m=1}^{ld^{\prime }}[l\lambda _i-\mu +\frac m{d^{\prime
}}(\lambda _j-\lambda _i)]\left( \frac{\lambda _j-\lambda _i}{d^{\prime }}%
\right) ^{(n+1-l)d-1}}{\left( 1+\frac{\lambda _i-\lambda j}{d^{\prime }\hbar 
}\right) \prod\limits_{\underset{(\alpha ,m)\neq (j,d^{\prime })}{\alpha =0}%
}^n\prod\limits_{m=1}^{d^{\prime }}[\lambda _i-\lambda _\alpha +\frac
m{d^{\prime }}(\lambda _j-\lambda _i)]}\cdot 
\end{equation*}

\begin{equation*}
\cdot \prod\limits_{k\neq j}(\lambda _j-\lambda _k)\int\limits_{{\cal K}%
_{G^{\prime \prime },j}}\frac{{\cal E}_{2,d-d^{\prime }}^{\prime
}e_0^{*}(\phi _j)}{{\cal E}({\cal N}_{{\cal K}_{G^{\prime \prime },j}})(%
\frac{\lambda _j-\lambda _i}{d^{\prime }}+c)}\int\limits_{{\cal K}%
_{d^{\prime },j}}1= 
\end{equation*}

\begin{equation*}
=\frac{\prod\limits_{m=1}^{ld^{\prime }}[l\lambda ^{\prime }-\mu +\frac
m{d^{\prime }}(\lambda _j-\lambda _i)]\left( \frac{\lambda _j-\lambda _i}{%
d^{\prime }}\right) ^{(n+1-l)d-1}}{d^{\prime }\left( 1+\frac{\lambda
_i-\lambda j}{d^{\prime }\hbar }\right) \prod\limits_{\underset{(\alpha
,m)\neq (j,d^{\prime })}{\alpha =0}}^n\prod\limits_{m=1}^{d^{\prime
}}[\lambda _i-\lambda _\alpha +\frac m{d^{\prime }}(\lambda _j-\lambda _i)]}%
\cdot
\end{equation*}
\begin{equation*}
\cdot
\prod\limits_{k\neq j}(\lambda _j-\lambda _k)\int\limits_{{\cal K}%
_{G^{\prime \prime },j}}\frac{{\cal E}_{2,d-d^{\prime }}^{\prime
}e_0^{*}(\phi _j)}{{\cal E}({\cal N}_{{\cal K}_{G^{\prime \prime },j}})(%
\frac{\lambda _j-\lambda _i}{d^{\prime }}+c)}\text{.} 
\end{equation*}
\begin{flushright}
$\Box$
\end{flushright}

\medskip\ 

\begin{theorem}
\label{twelve} Set $z_i(Q,\hbar )=Z_i(\hbar ^{n+1-l}Q,\hbar )$. Then
\end{theorem}

\begin{equation*}
z_i(Q,\hbar )=1+\sum\limits_{d^{\prime }>0}Q^{d^{\prime }}\sum\limits_{j\neq
i}coeff_i^j(d^{\prime })z_j(Q,\frac{\lambda _j-\lambda _i}{d^{\prime }})%
\text{,} 
\end{equation*}
where

\begin{equation*}
coeff_i^j(d^{\prime })=\frac{\prod\limits_{m=1}^{ld^{\prime }}[l\lambda
_i-\mu +\frac m{d^{\prime }}(\lambda _j-\lambda _i)]\left( \frac{\lambda
_j-\lambda _i}{d^{\prime }}\right) ^{(n+1-l)d^{\prime }-1}}{d^{\prime
}\left( 1+\frac{\lambda _i-\lambda j}{d^{\prime }\hbar }\right)
\prod\limits_{\underset{(\alpha ,m)\neq (j,d^{\prime })}{\alpha =0}%
}^n\prod\limits_{m=1}^{d^{\prime }}[\lambda _i-\lambda _\alpha +\frac
m{d^{\prime }}(\lambda _j-\lambda _i)]}\text{.} 
\end{equation*}

\begin{figure}[htb]
\hspace{-2cm} 
\begin{center} 
\mbox{\epsfig{file=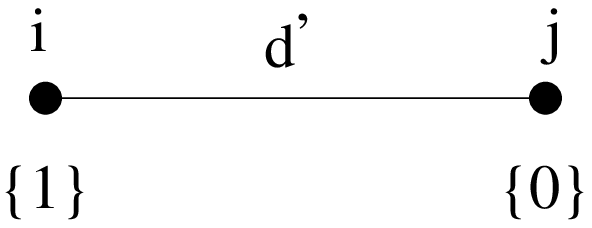,width=2cm,height=4cm,angle=270}}
\caption[]{
\label{dis31}}
\end{center}
\end{figure}

{\bf Proof. }By the integration formula over connected components ${\cal K}$
of fixed points, we obtain

\begin{equation*}
z_i(Q,\hbar )=1+\sum\limits_{d>0}Q^d\prod\limits_{j\neq i}(\lambda
_i-\lambda _j)\int\limits_{\overline{{\cal M}}_{0,2}({\Bbb P}^n,d)}\left(
e_0\right) ^{*}(\phi _i)\frac{{\cal E}_{2,d}^{\prime }}{1+\frac c\hbar }%
\left( -c\right) ^{(n+1-l)d-1}= 
\end{equation*}

\begin{equation*}
=1+\sum\limits_{d>0}Q^d\sum\limits_{{\cal K}}\prod\limits_{j\neq i}(\lambda
_i-\lambda _j)\int\limits_{{\cal K}}\left( e_0\right) ^{*}(\phi _i)\frac{%
{\cal E}_{2,d}^{\prime }}{\left( 1+\frac c\hbar \right) {\cal E}({\cal N}_{%
{\cal K}})}\left( -c\right) ^{(n+1-l)d-1}= 
\end{equation*}
\begin{equation*}
=1+\sum\limits_{d>0}Q^d\sum\limits_{j\neq i}\sum\limits_{d^{\prime
}=1}^d\sum\limits_{{\cal K}_{j,d^{\prime }}}\prod\limits_{j\neq i}(\lambda
_i-\lambda _j)\cdot
\end{equation*}
\begin{equation*}
\cdot
\int\limits_{{\cal K}_{j,d^{\prime }}}\left( e_0\right)
^{*}(\phi _i)\frac{{\cal E}_{2,d}^{\prime }}{\left( 1+\frac c\hbar \right) 
{\cal E}({\cal N}_{{\cal K}_{j,d^{\prime }}})}\left( -c\right) ^{(n+1-l)d-1}%
\text{,} 
\end{equation*}
where ${\cal K}_{j,d^{\prime }}$ is the connected component associated with
the graph of figure \ref{dis31}.

\ \ 

\noindent Therefore, by Lemma \ref{ten} and Lemma \ref{eleven} we have

\begin{equation*}
z_i(Q,\hbar )=1+\sum\limits_{d>0}Q^d\sum\limits_{j\neq i}\left\{
\sum\limits_{d^{\prime }=1}^{d-1}\left( \frac{\lambda _j-\lambda _i}{%
d^{\prime }}\right) ^{(n+1-l)(d-d^{\prime })}coeff_i^j(d^{\prime
})\prod\limits_{k\neq j}(\lambda _j-\lambda _k)\cdot \right. 
\end{equation*}

\begin{equation*}
\cdot \left. \sum\limits_{{\cal K}_{G^{\prime \prime },j}}\int_{{\cal K}%
_{G^{\prime \prime },j}}\frac{{\cal E}_{2,d-d^{\prime }}^{\prime
}e_0^{*}(\phi _j)}{{\cal E}({\cal N}_{{\cal K}_{G^{^{\prime \prime }},j}})(%
\frac{\lambda _j-\lambda _i}{d^{\prime }}+c)}+coeff_i^j(d)\right\} = 
\end{equation*}

\begin{equation*}
=1+\sum\limits_{d^{\prime }+(d-d^{\prime })>0}\sum\limits_{j\neq i}\left\{
\sum\limits_{d^{\prime }=1}^{d-1}Q^{d^{\prime }}\left( \frac{\lambda
_j-\lambda _i}{d^{\prime }}\right) ^{(n+1-l)(d-d^{\prime
})}coeff_i^j(d^{\prime })\prod\limits_{k\neq j}(\lambda _j-\lambda _k)\cdot
\right. 
\end{equation*}

\begin{equation*}
\cdot \left. \sum\limits_{{\cal K}_{G^{\prime \prime },j}}Q^{d-d^{\prime
}}\int\limits_{{\cal K}_{G^{\prime \prime },j}}\frac{{\cal E}_{2,d-d^{\prime
}}^{\prime }e_0^{*}(\phi _j)}{{\cal E}({\cal N}_{{\cal K}_{G^{^{\prime
\prime }},j}})(\frac{\lambda _j-\lambda _i}{d^{\prime }}+c)}%
+Q^dcoeff_i^j(d)\right\} = 
\end{equation*}

\begin{equation*}
=1+\sum\limits_{d^{\prime }+(d-d^{\prime })>0}\sum\limits_{j\neq i}\left\{
\sum\limits_{d^{\prime }=1}^{d-1}Q^{d^{\prime }}coeff_i^j(d^{\prime
})\prod\limits_{k\neq j}(\lambda _j-\lambda _k)\cdot \right. 
\end{equation*}

\begin{equation*}
\cdot \left. \sum\limits_{{\cal K}_{G^{\prime \prime },j}}Q^{d-d^{\prime
}}\int\limits_{{\cal K}_{G^{\prime \prime },j}}\frac{{\cal E}_{2,d-d^{\prime
}}^{\prime }e_0^{*}(\phi _j)(-c)^{(n+1-l)(d-d^{\prime })-1}}{{\cal E}({\cal N%
}_{{\cal K}_{G^{^{\prime \prime }},j}})(1+\frac{cd^{\prime }}{\lambda
_j-\lambda _i})}+Q^dcoeff_i^j(d)\right\} = 
\end{equation*}

\begin{equation*}
=1+\sum\limits_{d^{\prime }>0}Q^{d^{\prime }}\sum\limits_{j\neq
i}coeff_i^j(d^{\prime })\left\{ 1+\sum\limits_{d^{\prime \prime
}>0}Q^{d^{\prime \prime }}\prod\limits_{k\neq j}(\lambda _j-\lambda _k)\cdot
\right. 
\end{equation*}

\begin{equation*}
\cdot \left. \sum\limits_{{\cal K}_{G^{\prime \prime },j}}\int_{{\cal K}%
_{G^{\prime \prime },j}}\left( e_0\right) ^{*}(\phi _j)\frac{{\cal E}%
_{2,d^{\prime \prime }}^{\prime }}{{\cal E}({\cal N}_{{\cal K}_{G^{^{\prime
\prime }},j}})(1+\frac{cd^{\prime }}{\lambda _j-\lambda _i})}\left(
-c\right) ^{(n+1-l)d^{\prime \prime }-1}\right\} = 
\end{equation*}

\begin{equation*}
=1+\sum\limits_{d^{\prime }>0}Q^{d^{\prime }}\sum\limits_{j\neq
i}coeff_i^j(d^{\prime })z_j(Q,\frac{\lambda _j-\lambda _i}{d^{\prime }})%
\text{.} 
\end{equation*}
\begin{flushright}
$\Box$
\end{flushright}

We now determine explicitely the functions $z_i(Q,\hbar )$.

\begin{proposition}
\label{END} The functions $z_i(Q,1/\omega )$, $0\leq i\leq n$, $\omega
=1/\hbar $, are power series $\sum\limits_{d\geq 0}C_i(d,1/\omega )Q^d$ in $Q
$ with coefficients $C_i(d,1/\omega )$ which are rational functions of $%
\omega $ with poles of order at least one at $\omega =\frac{d^{\prime }}{%
\lambda _j-\lambda _i}$, with $d^{\prime }=1,\ldots ,d$. The functions $z_i$
are uniquely determined by these properties, the recursive relations of
Theorem \ref{twelve} and the initial condition $C_i(0)=1$.
\end{proposition}

{\bf Proof. }The proof is immediate from the recursion relations and from
the nature of the coefficients, which have poles only at the indicated
points.
\begin{flushright}
$\Box$
\end{flushright}

\begin{proposition}
\label{frac}The series 
\begin{equation}
z_i(Q,1/\omega )=\sum\limits_{d\geq 0}Q^d\frac{\prod\limits_{m=1}^{ld}[(l%
\lambda _i-\mu )\omega +m]}{d!\prod\limits_{\alpha \neq
i}\prod\limits_{m=1}^d[(\lambda _i-\lambda _\alpha )\omega +m]}
\label{formula}
\end{equation}
\end{proposition}

\noindent satisfy all the conditions of Proposition \ref{END}.

{\bf Proof.}We need to prove that the coefficients of functions in (\ref
{formula}) satisfy the following recursive relations:

\begin{equation*}
C_i(d,1/\omega )=\sum\limits_{j\neq i}\sum\limits_{a=1}^dcoeff_i^j(a)C_j(d-a,%
\frac{\lambda _j-\lambda _i}a)\text{.} 
\end{equation*}

\noindent To achieve this goal, we decompose the coefficients of the power series in (%
\ref{formula}) into sum of simple fractions, i.e.

\begin{equation*}
\frac{\prod\limits_{m=1}^{ld}[(l\lambda _i-\mu )\omega +m]}{%
d!\prod\limits_{\alpha \neq i}\prod\limits_{m=1}^d[(\lambda _i-\lambda
_\alpha )\omega +m]}= 
\end{equation*}

\begin{equation*}
=\sum\limits_{j\neq i}\sum\limits_{a=1}^d\frac 1{[(\lambda _i-\lambda
_j)\omega +a]}\frac{\prod\limits_{m=1}^{ld}\left( \frac{l\lambda _i-\mu }{%
\lambda _j-\lambda _i}a+m\right) }{d!\prod\limits_{\underset{(\alpha ,m)\neq
(j,a)}{\alpha \neq i}}\prod\limits_{m=1}^d\left( \frac{\lambda _i-\lambda
_\alpha }{\lambda _j-\lambda _i}a+m\right) }\text{.} 
\end{equation*}

\noindent On the other hand, we have 
\begin{equation*}
\sum\limits_{j\neq i}\sum\limits_{a=1}^dcoeff_i^j(a)C_j(d-a,\frac{\lambda
_j-\lambda _i}a)= 
\end{equation*}

\begin{eqnarray*}
&=&\sum\limits_{j\neq i}\sum\limits_{a=1}^d\frac 1{[(\lambda _i-\lambda
_j)\omega +a]}\frac{\prod\limits_{m=1}^{la}\left( \frac{l\lambda _i-\mu }{%
\lambda _j-\lambda _i}a+m\right) }{\prod\limits_{\underset{(\alpha ,m)\neq
(j,a)}{\alpha }}\prod\limits_{m=1}^a\left( \frac{\lambda _i-\lambda _\alpha 
}{\lambda _j-\lambda _i}a+m\right) }\cdot \\
&\cdot& \frac{\prod\limits_{m=1}^{l(d-a)}\left( 
\frac{l\lambda _j-\mu }{\lambda _j-\lambda _i}a+m\right) }{%
(d-a)!\prod\limits_{\underset{}{\alpha \neq j}}\prod\limits_{m=1}^{d-a}%
\left( \frac{\lambda _j-\lambda _\alpha }{\lambda _j-\lambda _i}a+m\right) }=
\\
&=&\sum\limits_{j\neq i}\sum\limits_{a=1}^d\frac 1{[(\lambda _i-\lambda
_j)\omega +a]}\frac{\prod\limits_{m=1}^{la}\left( \frac{l\lambda _i-\mu }{%
\lambda _j-\lambda _i}a+m\right) }{\prod\limits_{\underset{(\alpha ,m)\neq
(j,a)}{\alpha }}\prod\limits_{m=1}^a\left( \frac{\lambda _i-\lambda _\alpha 
}{\lambda _j-\lambda _i}a+m\right) }\frac{\prod\limits_{m=la+1}^{ld}\left( 
\frac{l\lambda _i-\mu }{\lambda _j-\lambda _i}a+m\right) }{%
\prod\limits_\alpha \prod\limits_{m=a+1}^d\left( \frac{\lambda _i-\lambda
_\alpha }{\lambda _j-\lambda _i}a+m\right) }=
\end{eqnarray*}
\begin{equation*}
=\sum\limits_{j\neq i}\sum\limits_{a=1}^d\frac 1{[(\lambda _i-\lambda
_j)\omega +a]}\frac{\prod\limits_{m=1}^{ld}\left( \frac{l\lambda _i-\mu }{%
\lambda _j-\lambda _i}a+m\right) }{(d-a)!\prod\limits_{\underset{}{\alpha
\neq j}}\prod\limits_{m=1}^d\left( \frac{\lambda _i-\lambda _\alpha }{%
\lambda _j-\lambda _i}a+m\right) \prod\limits_{m=1}^{a-1}(-a+m)}= 
\end{equation*}

\begin{equation*}
=\sum\limits_{j\neq i}\sum\limits_{a=1}^d\frac 1{[(\lambda _i-\lambda
_j)\omega +a]}\frac{\prod\limits_{m=1}^{ld}\left( \frac{l\lambda _i-\mu }{%
\lambda _j-\lambda _i}a+m\right) }{d!\prod\limits_{\underset{(\alpha ,m)\neq
(j,a)}{\alpha \neq i}}\prod\limits_{m=1}^d\left( \frac{\lambda _i-\lambda
_\alpha }{\lambda _j-\lambda _i}a+m\right) }=C_i(d,1/\omega )\text{.} 
\end{equation*}
\begin{flushright}
$\Box$
\end{flushright}

We can finally prove Theorem \ref{four} in the form \ref{five}.
\smallskip

{\bf Proof of Theorem \ref{four}.} By Proposition \ref{frac} and Lemma \ref
{six},

\begin{equation*}
\prod\limits_{j\neq i}\left( \lambda _i-\lambda _j\right) \left\langle \phi
_i,S^{^{\prime }}\right\rangle =e^{\frac{\lambda _it}\hbar }(l\lambda _i-\mu
)\sum\limits_{d\geq 0}\frac{e^{dt}}{\hbar ^{(n+1-l)d}}\frac{%
\prod\limits_{m=1}^{ld}[\frac{(l\lambda _i-\mu )}\hbar +m]}{%
d!\prod\limits_{\alpha \neq i}\prod\limits_{m=1}^d[\frac{(\lambda _i-\lambda
_\alpha )}\hbar +m]}\text{,} 
\end{equation*}
with $\omega =1/\hbar $. Hence

\begin{equation*}
\prod\limits_{j\neq i}\left( \lambda _i-\lambda _j\right) \left\langle \phi
_i,S^{^{\prime }}\right\rangle =e^{\frac{\lambda _it}\hbar }(l\lambda _i-\mu
)\sum\limits_{d\geq 0}\frac{e^{dt}}{\hbar ^d}\frac{\prod%
\limits_{m=1}^{ld}[(l\lambda _i-\mu )+m\hbar ]}{d!\prod\limits_{\alpha \neq
i}\prod\limits_{m=1}^d[(\lambda _i-\lambda _\alpha )+m\hbar ]}= 
\end{equation*}
\begin{equation*}
=\frac 1{2\pi \sqrt{-1}}\int \phi _ie^{\frac{pt}\hbar }\sum\limits_{d\geq 0}%
\frac{e^{dt}}{\hbar ^d}\sum\limits_{d\geq 0}e^{dt}\frac{\prod%
\limits_{m=1}^{ld}[(lp-\mu )+m\hbar ]}{\prod\limits_{\alpha
=0}^n\prod\limits_{m=1}^d[(p-\lambda _\alpha )+m\hbar ]\prod\limits_{\alpha
=0}^n(p-\lambda _\alpha )}dp 
\end{equation*}

\noindent This implies that

\begin{equation*}
S^{\prime }(t,\hbar )=e^{\frac{pt}\hbar }\sum\limits_{d\geq 0}e^{dt}\frac{%
\prod\limits_{m=1}^{ld}[(lp-\mu )+m\hbar ]}{\prod\limits_{\alpha
=0}^n\prod\limits_{m=1}^d[(p-\lambda _\alpha )+m\hbar ]}\text{.} 
\end{equation*}
\begin{flushright}
$\Box$
\end{flushright}

\bigskip\ 

\bigskip\ 

\section{Projective complete intersections with $l_1+\ldots +l_r=n$}

\bigskip\ 

Let $X\subset {\Bbb P}^n$ be a non-singular complete intersection given by
equations of degrees $l_1,\ldots ,l_r$ with $l_1+\ldots +l_r=n$. As in the
previous section, we look for recursive relations so as to determine
explicitly functions $s_\beta (t,\hbar )$.

\noindent However, in this case the recursive relations involve more terms, since we
have to modify Proposition \ref{type2} as explained in Lemma \ref{s5}, that
is, it is not any more true that these components give zero contribution.

\noindent This will imply that the two generating function, the one which contains
numbers counting curves, and the one containing solutions of Picard Fuchs
equation, do not coincide exactly, as they did in the previous case, but we
can pass from one another with suitable transformations.

\noindent Let us now explain in details; with the same notations introduced in the
previous section, let $S(t,\hbar )\in H^{*}({\Bbb P}^n)$ be defined as in 
\ref{esse} and denote by $S^{\prime }(t,\hbar )$ its equivariant
counterpart. The main theorem of this section is

\begin{theorem}
\label{s'2} Suppose $l_1+\ldots +l_r=n$. Then
\end{theorem}

\begin{equation*}
S^{\prime }(t,\hbar )=e^{\frac{pt-l_1!\ldots l_r!e^t}\hbar }\sum_{d\geq
0}e^{dt}\frac{\prod\limits_{m=0}^{dl_1}(l_1p-\mu _1+m\hbar )\ldots
\prod\limits_{m=0}^{dl_r}(l_rp-\mu _r+m\hbar )}{\prod\limits_{m=1}^d(p-%
\lambda _0+m\hbar )\ldots \prod\limits_{m=1}^d(l_1p-\lambda _n+m\hbar )}%
\text{.} 
\end{equation*}

\noindent As explained in Remark \ref{eqnoneq}, we deduce

\begin{theorem}
Suppose $l_1+\ldots +l_r=n${\bf . }Then
\end{theorem}

\begin{equation*}
S(t,\hbar )=e^{\frac{Pt-l_1!\ldots l_r!e^t}\hbar }\sum_{d\geq 0}e^{dt}\frac{%
\prod\limits_{j=1}^r\prod\limits_{m=1}^d(l_jP+m\hbar )}{\prod%
\limits_{m=1}^d(P+m\hbar )^{n+1}}\text{,} 
\end{equation*}
with $P$ the generator of the cohomology algebra of ${\Bbb P}^n$.

\noindent In this case, the functions $s_\beta (t,\hbar )=\int\limits_{{\Bbb P}%
^n}P^\beta S(t,\hbar )$ do not satisfy the Picard-Fuchs equation for the
mirror symmetric family of $X$. Nevertheless, if we multiply them by the
common factor $e^{\frac{l_1!\ldots l_r!e^t}\hbar }$ , we can prove that

\begin{equation*}
s_{_{(1)}\beta }^{\prime }(t,\hbar ):=e^{\frac{l_1!\ldots l_r!e^t}\hbar
}\int\limits_{{\Bbb P}^n}P^\beta S(t,\hbar )=e^{\frac{l_1!\ldots l_r!e^t}%
\hbar }s_\beta (t,\hbar )\text{, }\beta =0,\ldots ,m, 
\end{equation*}
satisfy the Picard - Fuchs differential equation

\begin{equation*}
\left[ \left( \hbar \frac d{dt}\right)
^{n+1-r}-e^t\prod\limits_{j=1}^rl_j\prod\limits_{m=1}^{l_j-1}\left( l_j\hbar
\frac d{dt}+\hbar m\right) \right] F(t,\hbar )=0\text{.} 
\end{equation*}

{\bf Proof. }It is exactly the same computation as in \ref{two}.
\begin{flushright}
$\Box$
\end{flushright}

Anyway, instead of modifying the functions, we prefer to modify the
differential operator, so to obtain information on $SQH^{*}\left( X\right) $.

\begin{proposition}
\label{stwo} The functions $s_\beta (t,\hbar )$ satisfy the differential
equation
\end{proposition}

\begin{equation*}
\left[ D^{n+1-r}-e^t\prod\limits_{j=1}^rl_j\prod\limits_{m=1}^{l_j-1}\left(
l_jD+\hbar m\right) \right] F(t,\hbar )=0\text{,} 
\end{equation*}
with $D=\hbar \frac d{dt}+l_1!\ldots l_r!e^t$.

{\bf Proof. }Let $f$ be a function of $t$. Then it is easy to prove by
induction on $m$ that

\begin{equation*}
\left[ D^mf(t)\right] e^{\frac{l_1!\ldots l_r!e^t}\hbar }=\left( \hbar \frac
d{dt}\right) ^m\left[ f(t)e^{\frac{l_1!\ldots l_r!e^t}\hbar }\right] \text{.}
\end{equation*}

\noindent This implies that

\begin{equation*}
\left( \hbar \frac d{dt}\right) ^m\left[ s_\beta (t,\hbar )e^{\frac{%
l_1!\ldots l_r!e^t}\hbar }\right] =\left[ D^ms_\beta (t,\hbar )\right] e^{%
\frac{l_1!\ldots l_r!e^t}\hbar }\text{.} 
\end{equation*}

\noindent Moreover, if we develop the differential operator

\begin{equation*}
e^t\prod\limits_{j=1}^rl_j\prod\limits_{m=1}^{l_j-1}\left( l_j\hbar \frac
d{dt}+\hbar m\right) 
\end{equation*}
as a polynomial in the powers of $\hbar \frac d{dt}$, we also obtain

\begin{equation*}
\begin{array}{c}
\left\{ e^t\prod\limits_{j=1}^rl_j\prod\limits_{m=1}^{l_j-1}\left( l_j\hbar
\frac d{dt}+\hbar m\right) \right\} \left[ s_\beta (t,\hbar )e^{\frac{%
l_1!\ldots l_r!e^t}\hbar }\right] = \\ 
=\left\{ \left\{ e^t\prod\limits_{j=1}^rl_j\prod\limits_{m=1}^{l_j-1}\left(
l_jD+\hbar m\right) \right\} \left[ s_\beta (t,\hbar )\right] \right\} e^{%
\frac{l_1!\ldots l_r!e^t}\hbar }\text{.}
\end{array}
\end{equation*}

\noindent So the result follows.
\begin{flushright}
$\Box$
\end{flushright}

\medskip\ 

\begin{corollary}
\label{enu}{\bf \ }In the small quantum cohomology algebra of $X$, $\dim X \neg 2$, the class $%
p$ of hyperplane sections satisfies the following relation: 
\[
(p+l_1!\ldots l_r!q)^{n+1-r}=l_1^{l_1}\ldots l_r^{l_r}q(p+l_1!\ldots
l_r!q)^{n-r}\text{.}
\]
\end{corollary}

{\bf Proof. }It follows directly from Proposition \ref{stwo} and
Proposition \ref{enumerative}
\begin{flushright}
$\Box$
\end{flushright}

\smallskip\ 

\begin{example}
{\bf \ }Let $X$ be a non singular cubic surface in ${\Bbb P}^3.$ Then the
relation found in Corollary \ref{enu} is
\end{example}

\begin{equation*}
(p+6q)^3=27q(6q+p)^2\text{,} 
\end{equation*}
or

\begin{equation*}
p^3=9qp^2+216q^2p+756q^3\text{.} 
\end{equation*}

\noindent If we denote by $\left\{ T_0=1,T_1=p,p^2\right\} $ a basis for $H^{*}(X)$,
then the intersection matrix is given by

\begin{equation*}
\left( 
\begin{array}{ccc}
0 & 0 & 1 \\ 
0 & 3 & 0 \\ 
1 & 0 & 0
\end{array}
\right) \text{.} 
\end{equation*}

In $SQH^{*}(X)$, we have

\begin{equation*}
T_1*T_1=\Phi _{110}T_2+\frac{\Phi _{111}}3T_1+\Phi _{112}T_0 
\end{equation*}
and

\begin{equation*}
T_1*T_1*T_1=\Phi _{110}(T_2*T_1)+\frac{\Phi _{111}}3(T_1*T_1)+\Phi _{112}T_1%
\text{,} 
\end{equation*}
with $\Phi $ the potential involving Gromov-Witten invariants. Since

\begin{equation*}
\Phi _{111}=qI_1(T_1,T_1,T_1)+\sum\limits_{d\geq 2}q^dI_d(T_1,T_1,T_1)\text{,%
} 
\end{equation*}
and $T_1*T_1*T_1=p^3$, this implies

\begin{equation*}
9q(T_1*T_1)=\frac q3I_1(T_1,T_1,T_1)(T_1*T_1), 
\end{equation*}
namely $I_1(T_1,T_1,T_1)=27$. In other words, by the very definiton of
Gromov-Witten invariants, there are $27$ lines on a non-singular cubic
surface in ${\Bbb P}^3$!

\medskip\ 

In order to prove Theorem \ref{s'2}, we proceed as in section 4. More
precisely, we will determine explicitely the functions $Z_i(q,\hbar )$
defined in the previous section. In this case, we have

\begin{proposition}
\label{s4}

\[
Z_i(q,\hbar )=\prod\limits_{j\neq i}(\lambda _i-\lambda
_j)\sum\limits_{d\geq 0}\frac{q^d}{\hbar ^d}\int_{\overline{{\cal M}}_{0,2}(%
{\Bbb P}^n,d)}\frac{e_0^{*}(\phi _i){\cal E}_{2,d}^{^{\prime }}(-c)^{d-1}}{%
1+\frac c\hbar },
\]
with $q=e^{dt}$.
\end{proposition}

{\bf Proof. }See Proposition \ref{zetai} and keep in mind that $l=n$.
\begin{flushright}
$\Box$
\end{flushright}

\medskip\ 

In this case the recursive relations are more complicated because of the
following

\begin{figure}[htb]
\hspace{-2cm} 
\begin{center} 
\mbox{\epsfig{file=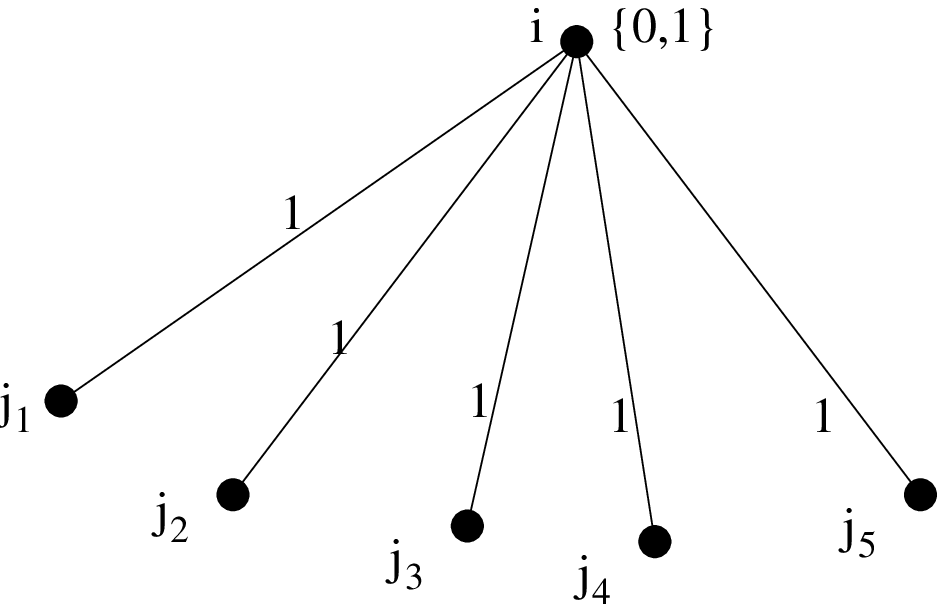,width=3cm,height=6cm,angle=270}}
\caption[]{
\label{dis5}}
\end{center}
\end{figure}

\begin{lemma}
\label{s5} Let $\Sigma $ $\in \overline{{\cal M}}_{0,2}({\Bbb P}%
^n,d)^{T^{n+1}}$ be a connected component of fixed points $\left[
C,x_0,x_1;f\right] $ whose marked point $x_0$ is situated on a component of $%
C$ with two or more special points.Then $\Sigma $ gives zero contribution to
the computation of the equivariant integrals 
\[
\int_{\overline{{\cal M}}_{0,2}({\Bbb P}^n,d)}(-c)^k{\cal E}_{2,d}^{\prime
}e_0^{*}(\phi _i)\text{, }k\geq d-1\text{,}
\]
unless $k=d$ and $C=C^{\prime }\cup C^{\prime \prime }$, where $C^{\prime }$
is mapped to a fixed point $p_i$ in ${\Bbb P}^n$ and carries both marked
points, and $C^{\prime \prime }$ is a disjoint union of $d$ irreducible
components (intersecting $C^{\prime }$ at $d$ special points) mapped (each
with multiplicity $1$) onto straight lines outgoing the point $p_i$ (see
figure \ref{dis5} for $d=5$).
\end{lemma}

{\bf Proof. }By Proposition \ref{type2}, the only non zero contribution
comes from those components whose associated graph $\Gamma $ (see section 4
for notations) satisfies the following requirements:

\begin{itemize}
\item  there exists a vertex $v$ with at least the label $\left\{ 0\right\} $%
;

\item  there are $l_v$ edges issuing from $v$ such that $h_v+l_v-3=d-1$.
\end{itemize}

\noindent This necessarily implies that $h_v=2$ and $l_v=d$. Then the Lemma follows.
\begin{flushright}
$\Box$
\end{flushright}

\medskip\ 

Set $z_i(Q,\hbar )=Z_i(\hbar Q,\hbar )$. We now compute the contribution of
these connected components.

\begin{proposition}
\label{s6}{\bf \ }The non-zero contribution to
$$
\prod\limits_{j\neq
i}(\lambda _i-\lambda _j)\sum\limits_{d\geq 0}Q^d \int_{\overline{{\cal M}}
_{0,2}({\Bbb P}^n,d)}(-c)^{d-1}{\cal E}_{2,d}^{\prime }e_0^{*}(\phi _i)
$$
given by connected components ${\cal K}$ described in Lemma \ref{s5} is
\end{proposition}

\begin{equation*}
\exp \left\{ Q\frac{\prod\limits_{a=1}^r(l_a\lambda _i-\mu _a)^{l_a}}{%
\prod\limits_{\alpha \neq i}(\lambda _i-\lambda _\alpha )}\right\} \exp
\left\{ -Ql_1!\ldots l_r!\right\} \text{.} 
\end{equation*}

{\bf Proof.} Let $\Gamma $ be the graph associated with ${\cal K}$. By our
assumptions, $\Gamma $ is the graph (see figure \ref{dis5} with $d=5$) with $%
j_s\neq i$, $s\in \left\{ 1,\ldots ,d\right\} $. Once we fix a $d-$ tuple $%
\left\{ j_1,\ldots ,j_d\right\} $, we denote by ${\cal K}_{j_1,\ldots ,j_d}$
the corresponding connected component and by $\Gamma _{j_1,\ldots ,j_d}$ its
associatd graph. Next recall that by ${\cal N}_{j_1,\ldots ,j_d}$ we mean
the normal bundle ${\cal N}_{\overline{{\cal M}}_{0,2}({\Bbb P}^n,d)/{\cal K}%
_{j_1,\ldots ,j_d}}$. Let us consider the contribution to

\begin{equation*}
\prod\limits_{j\neq i}(\lambda _i-\lambda _j)\sum\limits_{d\geq 0}Q^d\int_{%
\overline{{\cal M}}_{0,2}({\Bbb P}^n,d)}(-c)^{d-1}{\cal E}_{2,d}^{\prime
}e_0^{*}(\phi _i) 
\end{equation*}
given by all these connected components, i.e.

\begin{equation*}
\prod\limits_{j\neq i}(\lambda _i-\lambda _j)\sum\limits_{d\geq
0}Q^d\sum\limits_{(j_1,\ldots ,j_d)}\int_{{\cal K}_{j_1,\ldots ,j_d}}\frac{%
(-c)^{d-1}{\cal E}_{2,d}^{\prime }e_0^{*}(\phi _i)}{{\cal E}({\cal N}%
_{_{j_1,\ldots ,j_d}})} 
\end{equation*}
where the second sum ranges over all $d$-tuples of integers $j_s$, $0\leq
j_s\leq n$, none of which equals $i$. Consider now the morphism

\begin{equation*}
\delta :=\delta _{_{j_1,\ldots ,j_d}}:{\cal K}_{j_1,\ldots ,j_d}\rightarrow 
\overline{{\cal M}}_{0,2}({\Bbb P}^n,d)\text{.} 
\end{equation*}

\noindent We are going to make explicit computations for the case $r=1$ and then write
down formulas in the general case. As in {\bf STEP 1) }of Lemma \ref{ten},
the localization of $e_0^{*}(\phi _i)$ at ${\cal K}_{j_1,\ldots ,j_d}$ is 1.%
\smallskip\ 

Since the component $C^{\prime \prime }$ (where $x_0$ is situated) of
each point $\left[ C,x_0,x_1;f\right] \in {\cal K}_{j_1,\ldots ,j_d}$ is
mapped to the fixed point $p_i$ in ${\Bbb P}^n$, $T^{n+1}$ acts trivially on 
$\delta ^{*}T^{\left( 0\right) }$. By the localization formula, this means
that $c_1(\delta ^{*}T^{\left( 0\right) })$ ($c$ in the integral)$\in
H^{*}\left( {\cal K}_{j_1,\ldots ,j_d}\right) $. By Proposition \ref{graph} ${\cal K%
}_{_{j_1,\ldots ,j_d}}$ is isomorphic to $\overline{{\cal M}}_{0,d+2}$ via $%
\xi _{\Gamma _{_{j_1,\ldots ,j_d}}}$ and 

\noindent $\xi _{\Gamma _{_{j_1,\ldots
,j_d}}}^{*}\left( c_1(\delta ^{*}T^{\left( 0\right) })\right) =\widetilde{c}%
_0$, where $\widetilde{c}_0$ is the Chern class of {\sl the universal
tangent line bundle} of the marked point corresponding to $x_0$ on $%
\overline{{\cal M}}_{0,d+2}$.\smallskip\ 

We next compute the localization of ${\cal E}_{2,d}^{\prime }$ at $%
{\cal K}_{j_1,\ldots ,j_d}$, i.e. the Euler class of $\delta
^{*}W_{2,d}^{^{\prime }}$. By assumptions, we have $C=C_0\cup C_1\cup
\ldots C_d$, where each $C_s$ is mapped to $l_{ij_s}$ with degree one.
Besides, denote by $q_s=C_0\cap C_s$, $1\leq s\leq d$. Since $C_0$ is mapped
to the fixed point $p_i$, $\delta ^{*}W_{2,d}^{\prime }$ is trivial and thus
its equivariant Euler class belongs to $H_{T^{n+1}\times T^{^{\prime
}}}^{*}(B(T^{n+1}\times T^{\prime }))$. If $z_{1,s}$ and $z_{2,s}$, $%
s=1,\ldots ,d$, are homogeneous coordinates on $C_s$, then the action of $%
T^{n+1}$on $C^s$ is 
\begin{equation*}
[z_{1,s},z_{2,s}]\longrightarrow [t_iz_{1,s},t_jz_{j,s}]\text{.}
\end{equation*}

\noindent Therefore the character of the action of $T^{n+1}\times T^{\prime }$ on a
basis of $H^0(C_s,f_s^{*}{\cal O}(l))$, $f_s$ being the restriction of $f$
to $C_s$, is obtained as follows: 
\begin{equation*}
z_{1,s}^{l-m}z_{2,s}^m\rightarrow t_i^{l-m}t_j^mt^{\prime
-1}z_{1,s}^{l-m}z_{2,s}^m\text{,} 
\end{equation*}
with $t^{\prime }\in T^{\prime }$ and $0\leq m\leq l$. Hence the character
is 
\begin{equation*}
m\lambda _{j_s}-\mu +(l-m)\lambda _i=(l\lambda _i-\mu )+m(\lambda
_{j_s}-\lambda _i)\text{, }0\leq m\leq l\text{.} 
\end{equation*}

\noindent By the definition of $W_{2,d}^{\prime }$, the exact sequence holds 
\begin{equation*}
\begin{array}{ccccccccc}
0 & \rightarrow & H^0(C,f^{*}{\cal O}(l)) & \rightarrow & 
\oplus_{s=1}^dH^0(C_s,f_s^{*}({\cal O}(l))) & \rightarrow & 
\oplus_{s=1}^df_s^{*}({\cal O}(l)) & \rightarrow & 0 \\ 
&  &  &  & (\sigma _1,\ldots ,\sigma _d) & \rightarrow & (\sigma
_1(q_1),\ldots ,\sigma _d(q_d)) &  & 
\end{array}
\text{.} 
\end{equation*}

\noindent Then the localization of ${\cal E}_{2,d}^{\prime }$ at ${\cal K}_{j_1,\ldots
,j_d}$ is 
\begin{equation*}
\prod\limits_{s=1}^d\prod_{m=1}^l[(l\lambda _i-\mu )+m(\lambda
_{j_s}-\lambda _i)]\text{.} 
\end{equation*}

\noindent In general, for $r>1$, we have

\begin{equation*}
\prod\limits_{a=1}^r\prod\limits_{s=1}^d\prod_{m=1}^{l_a}[(l_a\lambda _i-\mu
_a)+m(\lambda _{j_s}-\lambda _i)]\text{.} 
\end{equation*}

\smallskip

We now compute the localization of ${\cal E}({\cal N}_{_{j_1,\ldots
,j_d}})$ at ${\cal K}_{j_1,\ldots ,j_d}$. To this end, we use the
observations made in Remark \ref{norbu}. Since the deformations which
preserve the combinatorial type of $C_0$ and $C_s$, $1\leq s\leq d$, are
trivial, we need to compute contributions of those deformations induced by
holomorphic vector fields on $C$. First of all, notice that the following
exact sequence holds: 
\begin{equation*}
0\rightarrow H^0(C,f^{*}T_{{\Bbb P}^n})\rightarrow H^0(C_0,f_0^{*}T_{{\Bbb P}%
^n})\oplus \oplus_{s=1}^dH^0(C_s,f_s^{*}T_{{\Bbb P}^n})\rightarrow
\oplus_{s=1}^df^{*}T_{p_i}\rightarrow 0\text{,}
\end{equation*}
with $f_0=f_{\mid _{C_0}}$ and $f_s=f_{\mid _{C_s}}$. For $C_0$, the only
contribution to the normal bundle comes from $H^0(C_0,f_0^{*}T_{{\Bbb P}^n})$
which is a trivial vector bundle on ${\cal K}_{j_1,\ldots ,j_d}$. Moreover,
its equivariant Euler class is $\prod\limits_{j\neq i}(\lambda _i-\lambda _j)
$, since a basis is given by $f_0^{*}\left( \frac \partial {\partial
X_j}\right) $, $j\in \left\{ 0,\ldots ,\widehat{i},\ldots ,n\right\} $.

\noindent The only contribution to the normal bundle stemming from deformations of
each $C_s$ is given by vector fields in 
\begin{equation*}
{\cal V}:=\oplus_{s=1}^d\frac{H^0(C_s,f_s^{*}(T_{{\Bbb P}^n}-p_i))}{%
H^0(C_s,T_{C_s}(-q_s))}\text{.} 
\end{equation*}

\noindent After observing (see {\bf STEP 4)} of Lemma \ref{eleven}) that a basis for $%
{\cal V}$ is given by vector fields $z_1f_s^{*}\left( \frac \partial
{\partial X_\alpha }\right) $, $\alpha \in \left\{ 0,\ldots ,\widehat{i},%
\widehat{j}_s,\ldots ,n\right\} $, ${\cal E}({\cal V})=\prod_{s=1}^d\prod_{%
\underset{\alpha \neq i}{\alpha \neq j_s}}(\lambda _{j_s}-\lambda _{\alpha})$.

\noindent Apart from these deformations, there is only the contribution of
deformations smoothing the nodes. Thus, we can write 
\begin{equation*}
{\cal E}({\cal N}_{j_1,\ldots ,j_d})={\cal E}(\oplus_{s=1}^d(T_{q_s}C_0%
\otimes T_{q_s}C_s))\prod\limits_{j\neq i}(\lambda _i-\lambda
_j)\prod_{s=1}^d\prod_{\underset{\alpha \neq i}{\alpha \neq j_s}}(\lambda
_{j_s}-\lambda _{\alpha})\text{.} 
\end{equation*}

\noindent Setting ${\cal E}(T_{q_s}C_s)=c_s$, then 
\begin{equation*}
{\cal E}(\oplus_{s=1}^d(T_{q_s}C_0\otimes
T_{q_s}C_s))=\prod_{s=1}^d(c_s+\lambda _i-\lambda _{j_s})\text{,} 
\end{equation*}
(see {\bf STEP 4)} of Lemma \ref{eleven} for the computation ${\cal E}%
(T_{q_s}C_s)$. We can also see this class as the pull back under $\xi
_{\Gamma _{j_1,\ldots ,j_d}}$ of the Chern class $\widetilde{c}_s$, $1\leq
s\leq d$, of the {\sl universal tangent bundle} of the marked point
corresponding to $q_s$ on $\overline{{\cal M}}_{0,d+2}$.

To summarize, we obtain 
\begin{equation*}
\prod\limits_{j\neq i}(\lambda _i-\lambda _j)\sum\limits_{d\geq
0}Q^d\sum\limits_{(j_1,\ldots ,j_d)}\int_{{\cal K}_{j_1,\ldots ,j_d}}\frac{%
(-c)^{d-1}{\cal E}_{2,d}^{\prime }e_0^{*}(\phi _i)}{{\cal E}({\cal N}%
_{_{j_1,\ldots ,j_d}})}= 
\end{equation*}
\begin{equation*}
=\sum\limits_{d\geq 0}Q^d\sum\limits_{(j_1,\ldots ,j_d)}\frac{%
\prod\limits_{a=1}^r\prod\limits_{s=1}^d\prod_{m=1}^{l_a}[(l_a\lambda _i-\mu
_a)+m(\lambda _{j_s}-\lambda _i )]}{\prod_{s=1}^d\prod_{\underset{\alpha \neq
i}{\alpha \neq j_s}}(\lambda _{j_s}-\lambda _{\alpha})}\cdot
\end{equation*}
\begin{equation*}
\cdot
\int_{{\cal K}_{j_1,\ldots
,j_d}}\frac{(-c)^{d-1}}{\prod_{s=1}^d(c_s+\lambda _i-\lambda _{j_s})}\text{.}
\end{equation*}

\noindent If we develop the cohomology class in the integral, we further have 
\begin{equation*}
\sum\limits_{d\geq 0}Q^d\sum\limits_{(j_1,\ldots ,j_d)}\frac{%
\prod\limits_{a=1}^r\prod\limits_{s=1}^d\prod_{m=1}^{l_a}[(l_a\lambda _i-\mu
_a)+m(\lambda _{j_s}-\lambda _i)]}{\prod_{s=1}^d\prod_{\underset{\alpha \neq
i}{\alpha \neq j_s}}(\lambda _{j_s}-\lambda _{\alpha})}\cdot 
\end{equation*}
\begin{equation*}
\cdot \prod_{s=1}^d\sum_{k_s\geq 0}\frac 1{(\lambda _i-\lambda
_{j_s})^{k_s+1}}\int_{{\cal K}_{j_1,\ldots ,j_d}}(-c)^{d-1}(-c_s)^{k_s+1}= 
\end{equation*}
\begin{eqnarray*}
&=&\sum\limits_{d\geq 0}Q^d\sum\limits_{(j_1,\ldots ,j_d)}\frac{%
\prod\limits_{a=1}^r\prod\limits_{s=1}^d\prod_{m=1}^{l_a}[(l_a\lambda _i-\mu
_a)+m(\lambda _{j_s}-\lambda _i)]}{\prod_{s=1}^d\prod_{\underset{\alpha \neq
i}{\alpha \neq j_s}}(\lambda _{j_s}-\lambda _{\alpha})}\cdot \\
&&\cdot \prod_{s=1}^d\sum_{k_s\geq 0}\frac 1{(\lambda _i-\lambda
_{j_s})^{k_s+1}}\int_{\overline{{\cal M}}_{0,d+2}}(-\widetilde{c}_0)^{d-1}(-%
\widetilde{c}_s)^{k_s+1}\text{.}
\end{eqnarray*}

\noindent On the other hand, $(-\widetilde{c}_i)$, $0\leq i\leq d+2$, are Chern
classes of {\sl universal cotangent line bundles} of marked points on $%
\overline{{\cal M}}_{0,d+2}$. By dimension computations, it is clear that we
have to compute the integral 
\begin{equation*}
\int_{\overline{{\cal M}}_{0,d+2}}(-\widetilde{c}_0)^{d-1} 
\end{equation*}
which, by standard results on intersection theory on $\overline{{\cal M}}%
_{0,d+2}$ is $\left\langle \tau _{d-1},\tau _0,\ldots ,\tau _0\right\rangle
=\left\langle \tau _0,\tau _0,\tau _0\right\rangle =1$ because of the string
equation (see \cite{Wi}).

\noindent In other words, we showed that 
\begin{equation*}
\prod\limits_{j\neq i}(\lambda _i-\lambda _j)\sum\limits_{d\geq
0}Q^d\sum\limits_{(j_1,\ldots ,j_d)}\int_{{\cal K}_{j_1,\ldots ,j_d}}\frac{%
(-c)^{d-1}{\cal E}_{2,d}^{\prime }e_0^{*}(\phi _i)}{{\cal E}({\cal N}%
_{_{j_1,\ldots ,j_d}})}= 
\end{equation*}
\begin{equation*}
=\sum\limits_{d\geq 0}Q^d(-1)^d\frac 1{d!}\frac{\prod\limits_{a=1}^r\prod%
\limits_{s=1}^d\prod_{m=1}^{l_a}[(l_a\lambda _i-\mu _a)+m(\lambda
_{j_s}-\lambda _i)]}{\prod_{s=1}^d\prod_{\underset{\alpha \neq i}{\alpha
\neq j_s}}(\lambda _{j_s}-\lambda _{\alpha})\prod_{s=1}^d(\lambda _{j_s}-\lambda _i)%
}= 
\end{equation*}
\begin{equation*}
=\sum\limits_{d\geq 0}Q^d(-1)^d\frac 1{d!}\left( \sum_{j\neq i}\frac{%
\prod\limits_{a=1}^r\prod_{m=1}^{l_a}[(l_a \lambda _i - \mu _a)+m(p-\lambda _i)]}{%
\prod_{\alpha \neq j}(p-\lambda _\alpha )}_{\mid p=\lambda _j}\right) ^d= 
\end{equation*}
\begin{equation*}
=\exp \left\{ -Q\sum_{j\neq i}\frac{\prod\limits_{a=1}^r\
\prod_{m=1}^{l_a}[(l_a \lambda _i - \mu _a)+m(p-\lambda _i)]}{\prod_{\alpha \neq
j}(p-\lambda _\alpha )}_{\mid p=\lambda _j}\right\} 
\text{.}\end{equation*}

Let us now consider the function 
\begin{equation*}
\ \frac{\prod\limits_{a=1}^r\prod_{m=1}^{l_a}[(l_ap-\mu _a)+m(p-\lambda _i)]%
}{\prod_{\alpha \neq j}(p-\lambda _\alpha )} 
\end{equation*}
as a meromorphic function of $p$. Then by taking residues in $p=\lambda
_\alpha $, $\alpha \neq j$, $p=$ $\infty $ and $p=\lambda _i$, we have 
\begin{equation*}
\sum_{j\neq i}\frac{\prod\limits_{a=1}^r\ \prod_{m=1}^{l_a}[(l_ap-\mu
_a)+m(p-\lambda _i)]}{\prod_{\alpha \neq j}(p-\lambda _\alpha )}_{\mid
p=\lambda _j}=\frac{\prod\limits_{a=1}^r(l_a\lambda _i-\mu _a)^{l_a}}{%
\prod\limits_{\alpha \neq i}(\lambda _i-\lambda _\alpha )}+l_1!\ldots l_r!%
\text{.} 
\end{equation*}

\noindent This completely proves the Lemma.
\begin{flushright}
$\Box$
\end{flushright}

\begin{theorem}
\label{s7}The functions $z_i(Q,\hbar )$ satisfy the recursive relations 
\begin{equation}
z_i(Q,\hbar )=1+\sum\limits_{d>0}coeff_i(d)Q^d+\sum\limits_{j\neq
i}\sum\limits_{d^{\prime }>0}Q^{d^{\prime }}coeff_i^j(d^{\prime })z_j(Q,%
\frac{\lambda _j-\lambda _i}{d^{\prime }}){\it \text{,}}  \label{rel}
\end{equation}
\end{theorem}
where 
\begin{equation*}
coeff_i^j(d^{\prime })=\frac{\prod\limits_{a=1}^r\prod\limits_{m=1}^{l_ad^{%
\prime }}[\frac{l_a\lambda ^{\prime }-\mu _a}{\lambda _j-\lambda _i}%
d^{\prime }+m]}{d^{\prime }\left( 1+\frac{\lambda _i-\lambda j}{d^{\prime
}\hbar }\right) \prod\limits_{\underset{(\alpha ,m)\neq (j,d^{\prime })}{%
\alpha =0}}^n\prod\limits_{m=1}^{d^{\prime }}[\frac{\lambda _i-\lambda
_\alpha }{\lambda _j-\lambda _i}d^{\prime }+m]}{\it \text{,}} 
\end{equation*}
and

\begin{equation*}
coeff_i(d)=\sum\limits_{t=1}^d\frac 1{t!(d-t)!}\frac{\prod%
\limits_{a=1}^r(l_a\lambda _i-\mu _a)^{l_at}}{\prod\limits_{\alpha \neq
i}^{}(\lambda _i-\lambda _\alpha )^t}(-l_1!\ldots l_r!)^{d-t}{\it \text{.}} 
\end{equation*}

{\bf Proof. }By Proposition \ref{s6} and Proposition \ref{twelve}, we obtain
the recursive relations.
\begin{flushright}
$\Box$
\end{flushright}

\medskip\ 

\begin{corollary}
\label{s8} The functions $z_i(Q,1/\omega )$, $\omega =1/\hbar $, are power
series 

\noindent $\sum\limits_{d\geq 0}C_i(d,1/\omega )Q^d$ where the coefficients $%
C_i(d,1/\omega )$ are rational functions of $\omega $ of the form 
\[
C_i(d,1/\omega )=\frac{P_{nd}^{(i)}\omega ^{nd}+\ldots +P_0^{(i)}}{%
\prod\limits_{\alpha =0}^m\prod\limits_{m=1}^d[(\lambda _i-\lambda _j)\omega
+m]}\text{,}
\]
\[
C_i(0,1/\omega )=1\text{.}
\]
The functions $z_i(Q,\hbar )$ are uniquely determined by these properties,
the recursive relations of Theorem \ref{s7} and the initial conditions
\end{corollary}

\begin{equation*}
\sum\limits_{d>0}coeff_i(d)Q^d=\sum\limits_{d>0}\frac{P_{nd}^{(i)}}{%
d!\prod\limits_{\alpha \neq i}(\lambda _i-\lambda _\alpha )^d}Q^d{\it \text{.%
}} 
\end{equation*}

{\bf Proof.} By the expression for $coeff_i(d)$ and $coeff_i^j(d^{\prime })$%
, we find out that the coefficients $C_i(d,1/\omega )$ have the desired
form, since $d\sum_{a=1}^rl_a=dn$. Notice that the term $P_{nd}^{(i)}\omega
^{nd}$ in the numerator of $C_i(d,1/\omega )$ comes from $%
\sum\limits_{d>0}coeff_i(d)Q^d$. Conversely, the requirements of this
Corollary completely determine the functions $z_i(Q,\hbar )$, since by \ref
{rel}, we have 
\begin{equation*}
C_i(0,1/\omega )=1\text{,} 
\end{equation*}

\begin{equation*}
C_i(d,1/\omega )=coeff_i(d)+\sum\limits_{j\neq
i}\sum\limits_{m=1}^dcoeff_i^j(m)C_i(d-m,\frac{\lambda _j-\lambda _i}m)\text{%
.} 
\end{equation*}
\begin{flushright}
$\Box$
\end{flushright}

\medskip\ 

\begin{theorem}
\label{s9} The series 
\[
z_i(Q,1/\omega )=\sum\limits_{d\geq 0}Q^d\frac{\prod\limits_{a=1}^r\prod%
\limits_{m=1}^{l_ad}[l_a\lambda _i-\mu _a)\omega +m]}{d!\prod\limits_{\alpha
\neq i}^n\prod\limits_{m=1}^d[(\lambda _i-\lambda _\alpha )\omega +m]}\exp
\left\{ -Ql_1!\ldots l_r!\right\} 
\]
satisfy the requirements of Corollary \ref{s8} and the initial conditions 
\[
1+\sum\limits_{d>0}coeff_i(d)Q^d=\exp \left\{ Q\frac{\prod%
\limits_{a=1}^r(l_a\lambda _i-\mu _a)^{l_a}}{\prod\limits_{\alpha \neq
i}(\lambda _i-\lambda _\alpha )}\right\} \exp \left\{ -Ql_1!\ldots
l_r!\right\} {\it \text{.}}
\]
\end{theorem}

{\bf Proof. }By developing the series $z_i(Q,\hbar )$, the coefficients $%
C_i(d,1/\omega )$ have the form

\begin{equation*}
\sum\limits_{t=1}^d\frac 1{t!(d-t)!}\frac{\prod\limits_{a=1}^r\prod%
\limits_{m=1}^{l_at}[(l_a\lambda _i-\mu _a)\omega +m]}{\prod\limits_{\alpha
\neq i}^n\prod\limits_{m=1}^t[(\lambda _i-\lambda _\alpha )\omega +m]}%
(-l_1!\ldots l_r!)^{d-t}\text{,} 
\end{equation*}
($C_i(0,1/\omega )=1$). Obviously, these coefficients have the desired form
and the initial conditions are satisfied. In order to verify the recursive
relations, set $z_i^{\prime }(Q,\hbar ):=z_i(Q,\hbar )\exp \{-l_1!\ldots
l_r!Q\}$, and proceed similarly to Proposition \ref{frac}.
\begin{flushright}
$\Box$
\end{flushright}

\medskip\ 

{\bf Proof of Theorem \ref{s'2}} See Theorem \ref{four}.
\begin{flushright}
$\Box$
\end{flushright}

\bigskip\ 

\bigskip\ 

\section{The Calabi-Yau case}

In this section we are going to show how to relate solutions $s^{*}(t,\hbar
) $ of the Picard-Fuchs equation for the mirror symmetric family of a
Calabi-Yau projective complete intersection $X$ in ${\Bbb P}^n$ and
solutions of the equation $\nabla _\hbar s(T,\hbar )=0$.

\noindent Once again, the two generating functions are different, and again this is
due to the fact that Proposition \ref{type2} does not hold. Even in this
case we can find recursive relations, even though of a much more complicated
type, satisfied by coefficients of both functions, but with different
initial conditions.

Givental makes use of certain ``polynomial properties'' to characterize a
class ${\cal P}$ of solutions of these relations. These properties are
satisfied by solutions of $\nabla _\hbar s(T,\hbar )=0$ for {\sl geometric
reasons }(see Lemma \ref{mainlemma} and Corollary \ref{cormain}) , and by
solutions of the Picard-Fuchs for {\sl computational reasons} (see
Proposition \ref{polycomp}).

\noindent Once this is established, Givental introduces a set of transformations that
preserves the class ${\cal P}$, and describes precisely the effect of any
such transformation on initial conditions; so it becomes clear how to
recover functions $s^{*}(t,\hbar )$ from solutions $s_\beta (T,\hbar )$ of
the equation $\nabla _\hbar s(T,\hbar )=0$.

\smallskip\ 

\subsection{Properties and motivations arising from geometry}

We recall that 
\begin{equation*}
s_\beta (T,\hbar )=\displaystyle\int\limits_{{\Bbb P}^n}P^\beta S_{\left(
2\right) }(T,\hbar )\text{,} 
\end{equation*}
with $S_{\left( 2\right) }(T,\hbar )=e^{\frac{PT}\hbar }\displaystyle%
\sum\limits_{d\geq 0}e^{dT}(e_0)_{*}\left( \frac{{\cal E}_{2,d}}{\hbar +c}%
\right) \in H^{*}({\Bbb P}^n)$.

\noindent Take now the equivariant counterpart $S_{\left( 2\right) }^{\prime }(T,\hbar
)$ of $S_{\left( 2\right) }(T,\hbar )$, and consider the function 
\begin{equation*}
\Phi (t,\tau ):=\int_{{\Bbb P}^n}^{T^{n+1}\times T^r}\frac{S_{\left(
2\right) }^{\prime }(t,\hbar )S_{\left( 2\right) }^{\prime }(\tau ,-\hbar )}{%
Euler(\oplus_i{\cal O}(l_i))}\text{.} 
\end{equation*}

\noindent This function is exactly the $(0,0)$ entry of the matrix 
\begin{equation*}
\left\{ -\hbar ^2\frac{\partial ^2{\cal F}^1\left( x\right) }{\partial
t_i\partial \tau _j}\right\} 
\end{equation*}
described in section \ref{solnab}, with normal integrals replaced by
equivariant ones.

\noindent From Theorem \ref{decomp}, taking in account that our basis $\left\{ \phi
_i\right\} $ is orthogonal, but not orthonormal, we can prove that 
\begin{equation*}
\begin{array}{c}
\Phi (t,\tau )=\sum_{d_1,d_2\geq 0}e^{d_1t}e^{d_2\tau }\sum_{i=0}^n%
\displaystyle\frac{\prod_{a=1}^r(l_a\lambda _i-\mu _a)}{\prod_{j\neq
i}(\lambda _i-\lambda _j)}\displaystyle\int_{\overline{{\cal M}}_{0,2}({\Bbb %
P}^n,d_1)}\displaystyle\frac{e^{\left( \frac{pt}\hbar \right) }{\cal E}%
_{2,d_1}^{\prime }(e_0)^{*}(\phi _i)}{\hbar +c}\cdot \\ 
\cdot \displaystyle\int_{\overline{{\cal M}}_{0,2}({\Bbb P}^n,d_2)}%
\displaystyle\frac{e^{\left( -\frac{p\tau }\hbar \right) }{\cal E}%
_{2,d_2}^{\prime }(e_0)^{*}(\phi _i)}{-\hbar +c}\text{.}
\end{array}
\end{equation*}

\smallskip\ 

Let 
\begin{equation*}
L_d:=\left\{ 
\begin{array}{c}
\left[ C,x_0,x_1,(\psi _1,\psi _2)\right] \in \overline{{\cal M}}%
_{0,2}\left( {\Bbb P}^n\times {\Bbb P}^1,\left( d,1\right) \right) : \\ 
\psi _2\left( x_0\right) \in {\Bbb P}^n\times \left\{ 0\right\} ,\psi
_2\left( x_1\right) \in {\Bbb P}^n\times \left\{ \infty \right\}
\end{array}
\right\} \supset \left( \overline{{\cal M}}_{0,2}\left( {\Bbb P}^n\times 
{\Bbb P}^1,\left( d,1\right) \right) \right) ^{S^1} 
\end{equation*}
and $L_d^{\prime }:={\Bbb P}\left( {\Bbb C}^{n+1}\otimes H^0\left( {\Bbb P}^1 ,
{\cal O}\left( d\right) \right) \right) $; points in this space are $%
(n+1)$-tuples of homogeneous polynomials of degree $d$ in two variables.

\noindent The group $T=S^1\times \left( S^1\right) ^{n+1}$ acts on both spaces: the
action on $L_d$ has been already described, and the action on $L_d^{\prime }$
is induced by the diagonal action of $\left( S^1\right) ^{n+1}$ on ${\Bbb C}%
^{n+1}$, and by rotations on ${\Bbb P}^1$; more precisely we consider the following action:
\begin{equation}
\begin{array}{ccccc}
S^1 & \times & {\Bbb P}^1 & \rightarrow & {\Bbb P}^1\\ 
t & , & \left[ z_1 , z_2 \right] &  & \left[ t^{-1}  z_1 , z_2 \right]
\end{array}
\text{.}  
\end{equation}

\noindent Let $p=c_1^T\left( {\cal O}\left( 1\right) \right) \in H_T^{*}\left(
L_d^{\prime }\right) $.

\begin{lemma}
\label{mainlemma} {\bf ``The Main Lemma''.} There exists a regular $T$%
-equivariant map 
\[
\mu :L_d\rightarrow L_d^{\prime }\text{,}
\]
such that $\Phi \left( t,\tau \right) =\sum_{d\geq 0}\exp \left( d\tau
\right) \int_{L_d}\exp \left( \frac{\mu ^{*}\left( p\right) \left( t-\tau
\right) }\hbar \right) {\cal E}_{2,d}$.
\end{lemma}

First of all, let us give some motivation for proving this lemma; the proof
is postponed at the end of this section.

\noindent Consider the bundle $F_{2,d}$ on $L_d^{\prime }$ such that substitution of
the $n+1$ polynomials into $r$ homogeneous equations of degrees $l_1,...,l_r$
produces a section of this bundle. Introducing the addictional action of a
torus $T^{\prime }=\left( S^1\right) ^r$ on this bundle, we can see that 
\begin{equation*}
H_{T\times T^{\prime }}^{*}\left( L_d^{\prime }\right) \cong \frac{{\Bbb Q}%
\left[ p,\lambda _0,...,\lambda _n,\mu _1,...,\mu _r,\hbar \right] }{%
\prod_{j=0}^n\prod_{m=0}^d\left( p-\lambda _j-m\hbar \right) }\text{,} 
\end{equation*}
and that 
\begin{equation*}
{\cal F}_{2,d}:=Euler^{T\times T^{\prime }}\left( F_{2,d}\right)
=\prod_{a=1}^r\prod_{m=0}^{l_ad}\left( l_ap-\mu _a-m\hbar \right) \text{.} 
 \end{equation*}
 
\noindent Assume for a moment that $\mu ^{*}\left( {\cal F}_{2,d}\right) ={\cal E}%
_{2,d}$; then substituting in the result of the lemma we wolud have that 
\begin{eqnarray*}
\Phi \left( t,\tau \right) &=&\sum_{d\geq 0}\exp \left( d\tau \right)
\int_{L_d}\exp \left( \frac{\mu ^{*}\left( p\right) \left( t-\tau \right) }%
\hbar \right) \mu ^{*}\left( {\cal F}_{2,d}\right) = \\
&=&\sum_{d\geq 0}\exp \left( d\tau \right) \int_{L_d^{\prime }}\exp \left( 
\frac{p\left( t-\tau \right) }\hbar \right) {\cal F}_{2,d}= \\
&=&\frac 1{2\pi \sqrt{-1}} \cdot \\
&\cdot& \int \exp \left( \frac{p\left( t-\tau \right) }%
\hbar \right) \left( \sum_{d\geq 0}\exp \left( d\tau \right) \frac{%
\prod_{a=1}^r\prod_{m=0}^{l_ad}\left( l_ap-\mu _a-m\hbar \right) }{%
\prod_{j=0}^m\prod_{m=0}^d\left( p-\lambda _j-m\hbar \right) }\right) dp%
\text{,}
\end{eqnarray*}
and we would proceed exactly as in the simpler case $\sum_il_i<n$ .
Unfortunately, the assumption is false, and consequently the function $\Phi
\left( t,\tau \right) $ does not have exactly this form; anyway, the close
relation between the spaces $L_d$ and $L_d^{\prime }$ established in the
lemma allows us to write a similar expression for $\Phi \left( t,\tau
\right) $; namely, setting $q=\exp \left( \tau \right) $, $z=\frac{t-\tau }%
\hbar $, and $\Phi ^{\prime }(q,z)=\Phi (t,\tau )$, we can prove:

\begin{corollary}
\label{cormain} 
\[
\Phi ^{\prime }(q,z)=\frac 1{2\pi \sqrt{-1}}\int \exp \left( pz\right)
\left( \sum_{d\geq 0}\frac{q^dE_d\left( p,\lambda ,\mu ,\hbar \right) }{%
\prod_{j=0}^m\prod_{m=0}^d\left( p-\lambda _j-m\hbar \right) }\right) dp%
\text{,}
\]
where $E_d\left( p,\lambda ,\mu ,\hbar \right) $ $=\mu _{*}\left( {\cal E}%
_{2,d}\right) $ depends polynomially on all its variables.
\end{corollary}

{\bf Proof.} Apply integration formula: 
\begin{eqnarray*}
\Phi ^{\prime }(q,z) &=&\sum_{d\geq 0}q^d\int_{L_d}\exp \left( \mu
^{*}\left( p\right) z\right) {\cal E}_{2,d}= \\
\ &=&\sum_{d\geq 0}q^d\int_{L_d^{\prime }}\exp \left( pz\right) \mu
_{*}\left( {\cal E}_{2,d}\right) = \\
\ &=&\frac 1{2\pi \sqrt{-1}}\int \exp \left( pz\right) \left( \sum_{d\geq 0}%
\frac{q^dE_d\left( p,\lambda ,\mu ,\hbar \right) }{\prod_{j=0}^m%
\prod_{m=0}^d\left( p-\lambda _j-m\hbar \right) }\right) dp\text{.}
\end{eqnarray*}

\noindent Recall that ${\cal E}_{2,d}\in H_T^{*}(L_d)$ considered as a ${\Bbb Q}%
[\lambda _0,\ldots \lambda _n]$ module, hence a polynomial.
\begin{flushright}
$\Box$
\end{flushright}

\smallskip\ 

Take now the equivariant counterpart $S_{\left( 2\right) }^{\prime }(T,\hbar
)$ of $S(T,\hbar )$: if we compute the components of $S_{\left( 2\right)
}^{\prime }(T,\hbar )$ with respect to the basis $\phi _i$ of $%
H_{T^{n+1}\times T^r}({\Bbb P}^n)$ (see Lemma \ref{six}), the
functions 
\begin{equation*}
Z_i(q,\hbar ):=\sum_{d\geq 0}q^d\int_{\overline{{\cal M}}_{0,2}({\Bbb P}%
^n,d)}(e_0)^{*}(\phi _i)\frac{{\cal E}_{2,d}^{\prime }}{\hbar +c} 
\end{equation*}
are not yet written in an explicit form.

\medskip\ 

\subsection{The generating function for solutions of the Picard-Fuchs}

Let us now define the function 
\begin{equation*}
S_{(2)}^{*}(t,\hbar )=e^{\frac{Pt}\hbar }\sum_{d\geq 0}e^{dt}\frac{%
\prod_{a=1}^r\prod_{m=0}^{l_ad}(l_aP+m\hbar )}{\prod_{m=1}^d(P+m\hbar )^{n+1}%
}\text{.} 
\end{equation*}

\noindent On the basis of the analysis made in previous sections, we know that its
components $s_i^{*}(t,\hbar )$ satisfy the Picard-Fuchs equation (see
Corollary \ref{two}). We then consider the equivariant counterpart 
\begin{equation*}
S_{(2)}^{^{\prime }*}(t,\hbar )=e^{\frac{pt}\hbar }\sum_{d\geq 0}e^{dt}\frac{%
\prod_{a=1}^r\prod_{m=0}^{l_ad}(l_ap-\mu _a+m\hbar )}{\prod_{i=0}^n%
\prod_{m=1}^d(p-\lambda _i+m\hbar )^{n+1}} 
\end{equation*}
and introduce power series 
\begin{equation*}
Z_i^{*}(q,\hbar )=\sum_{d\geq 0}q^d\frac{\prod_{a=1}^r\prod_{m=1}^{l_ad}(l_a%
\lambda _i-\mu _a+m\hbar )}{\prod_{\alpha =0}^n\prod_{m=1}^d(\lambda
_i-\lambda _\alpha +m\hbar )}=\sum_{d\geq 0}q^dC_i^{*}(d,\hbar )\text{,} 
\end{equation*}
so as to obtain, as proved in Lemma \ref{six},

\begin{equation*}
s_i^{^{\prime }*}(t,\hbar )=\prod_{j\neq i}(\lambda _j-\lambda _i)e^{\frac{%
\lambda _it}\hbar }\prod_{a=1}^r(l_a\lambda _i-\mu _a)Z_i^{*}(e^t,\hbar )%
\text{.} 
\end{equation*}

\smallskip\ 

\subsection{Recursive ``Calabi-Yau'' relations}

First of all, we are going to prove that both $Z_i(q,\hbar )$ and $%
Z_i^{*}(q,\hbar )$ satisfy the same recursive relations. More precisely, as
in subsections 4.2 and 4.3, we will determine recursive relations for $%
Z_i(q,\hbar )$, $0\leq i\leq n,$ with the aid of the integration formula
over connected components of fixed points in $\overline{{\cal M}}_{0,2}(%
{\Bbb P}^n,d)$ and thereafter we will prove that $Z_i^{*}(q,\hbar )$ satisfy
the same recursive relations.
\smallskip

\begin{proposition}
\label{recursy} 
\[
Z_i(q,\hbar )=1+\sum_{d>0}\frac{q^d}{\hbar ^d}\frac{R_{i,d}}{d!}%
+\sum_{d^{\prime }>0}\sum_{j\neq i}\frac{q^{d^{\prime }}}{\hbar
^{d^{^{\prime }}}}\frac{coeff_i^j(d^{\prime })}{\lambda _i-\lambda
_j+d^{\prime }\hbar }Z_j\left( \frac q\hbar \frac{\lambda _j-\lambda _i}{%
d^{\prime }},\frac{\lambda _j-\lambda _i}{d^{\prime }}\right) ,
\]
with $R_{i,d}$ a polynomial of $(\hbar ,\lambda _\alpha ,\mu _a)$ and 
\[
coeff_i^j(d^{\prime })=\frac{\prod_{a=1}^r\prod_{m=1}^{l_ad^{\prime }}\left[
l_a\lambda _i-\mu _a+\frac m{d^{\prime }}(\lambda _j-\lambda _i)\right] }{%
d^{\prime }!\underset{(\alpha ,m)\neq (j,d^{\prime })}{\prod_{\alpha \neq
i}\prod_{m=1}^{d^{\prime }}}\left[ \lambda _i-\lambda _\alpha +\frac
m{d^{\prime }}(\lambda _j-\lambda _i)\right] }\text{.}
\]
\end{proposition}

{\bf Proof.} If we expand $\frac 1{\hbar +c}$ in power series, then 
\begin{equation*}
Z_i(q,\hbar )=1+\sum_{d>0}q^d\sum_{k\geq 0}\frac 1{\hbar ^{k+1}}\int_{%
\overline{{\cal M}}_{0,2}({\Bbb P}^n,d)}{\cal E}_{2,d}^{\prime }e_0^{*}(\phi
_i)(-c)^k\text{.} 
\end{equation*}

\noindent Since dim$\overline{{\cal M}}_{0,2}({\Bbb P}^n,d)=n+(n+1)d-1$, 
\begin{equation*}
Z_i(q,\hbar )=1+\sum_{d>0}q^d\left[ \sum_{k=0}^{d-1}\hbar ^{-k-1}\int_{%
\overline{{\cal M}}_{0,2}({\Bbb P}^n,d)}{\cal E}_{2,d}^{\prime }e_0^{*}(\phi
_i)(-c)^k\right] + 
\end{equation*}
\begin{equation*}
+\sum_{d>0}q^d\hbar ^{-d}\int_{\overline{{\cal M}}_{0,2}({\Bbb P}^n,d)}\frac{%
{\cal E}_{2,d}^{\prime }e_0^{*}(\phi _i)(-c)^d}{\hbar +c}\text{.} 
\end{equation*}

\noindent By Proposition \ref{type2} and Lemma \ref{s5}, the last sum is equal to 
\begin{equation*}
\sum_{d^{\prime }>0}q^{d^{\prime }}\hbar ^{-d^{\prime }}\frac{%
\prod_{a=1}^r\prod_{m=1}^{l_ad^{\prime }}\left[ l_a\lambda _i-\mu _a+\frac
m{d^{\prime }}(\lambda _j-\lambda _i)\right] }{d^{\prime }!\underset{(\alpha
,m)\neq (j,d^{\prime })}{\prod_{\alpha \neq i}\prod_{m=1}^{d^{\prime }}}%
\left[ \lambda _i-\lambda _\alpha +\frac m{d^{\prime }}(\lambda _j-\lambda
_i)\right] (\lambda _j-\lambda _i+d^{\prime }\hbar )}\cdot
\end{equation*}
\begin{equation*}
\cdot
Z_j\left( \frac q\hbar 
\frac{\lambda _j-\lambda _i}{d^{\prime }},\frac{\lambda _j-\lambda _i}{%
d^{\prime }}\right) \text{.} 
\end{equation*}

\noindent On the other hand, the first sum gives the contribution $\frac{R_{i,d}}{d!}.$ $\Box$

Obviously, we have the following

\begin{corollary}
\label{coij}The coefficients of the power series $Z_i(q,\hbar )=\sum_{d\geq
0}q^dC_i(d,\hbar )$ are rational functions 
\[
C_i(d,\hbar )=\frac{P_d^{\left( i\right) }}{d!\hbar ^d\prod_{j\neq
i}\prod_{m=1}^d(\lambda _i-\lambda _j+m\hbar )},
\]
where $P_d^{\left( i\right) }$ is a polynomial in $(\hbar ,\lambda ,\mu )\,$
of degree $(n+1)d$.
\end{corollary}
\begin{flushright}
$\Box$
\end{flushright}

\begin{proposition}
The functions $Z_i^{*}(q,\hbar )$ satisfy the recursive relations of
Proposition \ref{recursy}.
\end{proposition}

{\bf Proof. }See the proof of Proposition \ref{frac}. Notice that initial
conditions are obtained by dividing the numerator in $C_i^{*}(d,\hbar )$ by
the denominator.
\begin{flushright}
$\Box$
\end{flushright}

In the sequel, recursive relations as the ones satisfied by $Z_i(q,\hbar )$,
but possibly with $R_{i,d}$ polynomials of $(h,\lambda _\alpha ,\mu _a)$ of
degree at most $d,$ will be referred to as ``{\it Calabi -Yau relations''}.

\noindent We can prove a result of uniqueness.

\begin{proposition}
\label{coeff}For any given polynomial $R_{i,d}$, the ``{\it Calabi -Yau
relations''} are satisfied by a unique solution of the form $\sum_{d\geq
0}\left( \frac q\hbar \right) ^dK_i(d,\hbar )$, with $K_i(0,\hbar )=1$ and $%
K_i(d,\hbar )$ rational functions of $\hbar $ as in Corollary \ref{coij}
with numerator $F_d^{\left( i\right) }$.
\end{proposition}

{\bf Proof.} If we substitute functions $W_i(q,\hbar )$ in the relations,
then we see that coefficients $K_i(d,\hbar )$ are uniquely determined by the
relations 
\begin{equation*}
K_i(d,\hbar )=\frac 1{d!}R_{i,d}+\sum_{j\neq i}\sum_{m=1}^d\frac{Coeff_i^j(m)%
}{m!(\lambda _j-\lambda _i+m\hbar )}\left( \frac{\lambda _j-\lambda _i}%
m\right) ^{d-m}K_j\left( d-m,\frac{\lambda _j-\lambda _i}m\right) \text{.} 
\end{equation*}
\begin{flushright}
$\Box$
\end{flushright}

Both functions $Z_i(q,\hbar )$ and $Z_i^{*}(q,\hbar )$ enjoy another
remarkable property.

\subsection{The class of solutions with polynomial properties}

Fix a generic solution $W_i(q,\hbar )$ of the ``{\it Calabi-Yau relations'' }%
and define the correlator 
\begin{equation}
\Phi _{W_i}(q,z,\hbar ):=\sum_{i=0}^n\frac{\prod_{a=1}^r(l_a\lambda _i-\mu
_a)}{\prod_{j\neq i}(\lambda _i-\lambda _j)}e^{z\lambda _i}W_i(qe^{z\hbar
},\hbar )W_i(q,-\hbar )\text{.}  \label{fizeta}
\end{equation}

\noindent Let us compute the correlator for $Z_i(q,\hbar )$ and $Z_i^{*}(q,\hbar )$.
To this end, we recall that $L_d^{\prime }$ denotes the projective space
of $(n+1)$-tuples of polynomials in the complex variable $z$ of degree at
most $d$ and $p$ denotes the equivariant first Chern class of the Hopf
bundle over $L_d^{\prime }$.

\begin{proposition}
\label{polycomp} 
\[
\begin{array}{c}
\Phi _{Z_i^{*}}(q,z,\hbar )=\frac 1{2\pi \sqrt{-1}}\int e^{pz}\displaystyle%
\sum_{d\geq 0}q^d\frac{\prod_{a=1}^r\prod_{m=0}^{l_ad}(l_ap-\mu _a+m\hbar )}{%
\prod_{\alpha =0}^n\prod_{m=0}^d(p-\lambda _\alpha -m\hbar )}dp= \\ 
\\ 
=\int_{L_d^{\prime }}e^{pz}\sum_{d\geq
0}q^d\prod_{a=1}^r\prod_{m=1}^{l_ad}(l_ap-\mu _a+m\hbar )\text{.}
\end{array}
\]
\end{proposition}

{\bf Proof. }We have 
\begin{equation*}
\begin{array}{c}
\displaystyle\sum_{i=0}^n\displaystyle\frac{\prod_{a=1}^r(l_a\lambda _i-\mu
_a)}{\prod_{j\neq i}(\lambda _i-\lambda _j)}e^{z\lambda
_i}Z_i^{*}(qe^{z\hbar },\hbar )Z_i^{*}(q,-\hbar )= \\ 
\\ 
=\displaystyle\sum_{i=0}^n\displaystyle\frac{\prod_{a=1}^r(l_a\lambda _i-\mu
_a)}{\prod_{j\neq i}(\lambda _i-\lambda _j)}e^{z\lambda _i}\left( %
\displaystyle\sum_{d_1\geq 0}q^{d_1}e^{d_1z\hbar }C_i^{*}(d_1,\hbar )\right)
\left( \displaystyle\sum_{d_2\geq 0}q^{d_2} C_i^{*}(d_2,-\hbar
)\right) = \\ 
\\ 
=\displaystyle\sum_{i=0}^n\displaystyle\frac{\prod_{a=1}^r(l_a\lambda _i-\mu
_a)}{\prod_{j\neq i}(\lambda _i-\lambda _j)}e^{z\lambda _i}\displaystyle%
\sum_{d\geq 0}q^d\displaystyle\sum_{s=0}^de^{sz\hbar }\left( \displaystyle%
\frac{\prod_{a=1}^r\prod_{m=1}^{l_as}(l_a\lambda _i-\mu _a+m\hbar )}{%
\prod_{\alpha =0}^n\prod_{m=0}^s(\lambda _i-\lambda _\alpha +m\hbar )}\cdot
\right. \\ 
\cdot \left. \displaystyle\frac{\prod_{a=1}^r\prod_{m=1}^{l_a(d-s)}(l_a%
\lambda _i-\mu _a - m\hbar )}{\prod_{\alpha =0}^n\prod_{m=0}^{(d-s)}(\lambda
_i-\lambda _\alpha -m\hbar )}\right) = \\ 
\\ 
=\displaystyle\sum_{i=0}^n\displaystyle\sum_{d\geq 0}q^d\displaystyle%
\sum_{s=0}^de^{z(\lambda _i+s\hbar )}\displaystyle\frac{\prod_{a=1}^r%
\prod_{m=0}^{l_ad}(l_a\lambda _i-\mu _a-m\hbar +l_as\hbar )}{\underset{%
(\alpha ,m)\neq (i,s)}{\prod_{\alpha =0}^n\prod_{m=0}^d}(\lambda _i-\lambda
_\alpha -m\hbar +s\hbar )}= \\ 
\\ 
=\frac 1{2\pi \sqrt{-1}}\int e^{pz}\displaystyle\sum_{d\geq 0}q^d%
\displaystyle\frac{\prod_{a=1}^r\prod_{m=0}^{l_ad}(l_ap-\mu _a+m\hbar )}{%
\prod_{\alpha =0}^n\prod_{m=0}^d(p-\lambda _\alpha -m\hbar )}dp\text{.}
\end{array}
\end{equation*}
\begin{flushright}
$\Box$
\end{flushright}

On the other hand, by the definition of $Z_i(q,\hbar )$, we have 
\begin{equation*}
\Phi _{Z_i}(q,z,\hbar )=\sum_{i=0}^n\displaystyle\frac{\prod_{a=1}^r(l_a%
\lambda _i-\mu _a)}{\prod_{j\neq i}(\lambda _i-\lambda _j)}e^{z\lambda
_i}\cdot
\end{equation*}
\begin{equation*}
\cdot
\left( \displaystyle\sum_{d_1\geq 0}q^{d_1}e^{d_1z\hbar }\displaystyle%
\int_{\overline{{\cal M}}_{0,2}({\Bbb P}^n,d_1)}\displaystyle\frac{{\cal E}%
_{2,d_1}^{\prime }(e_0)^{*}(\phi _i)}{\hbar +c}\right) \cdot
\end{equation*}
\begin{equation*}
\cdot \left( \displaystyle\sum_{d_2\geq 0}q^{d_2}\displaystyle\int_{%
\overline{{\cal M}}_{0,2}({\Bbb P}^n,d_2)}\displaystyle\frac{{\cal E}%
_{2,d_2}^{\prime }(e_0)^{*}(\phi _i)}{-\hbar +c}\right) \text{.}
\end{equation*}

\noindent Notice that, if we set $q=e^\tau $, $z=\frac{t-\tau }\hbar $, then $\Phi
_{Z_i}(q,z,\hbar )=\Phi (t,\tau )$.

\noindent By the Main Lemma, we see that also for $Z_i(q,\hbar )$ the correlator $\Phi
_{Z_i}(q,z,\hbar )$ can be written as an equivariant integral over $%
L_d^{\prime }$, i.e. 
\begin{equation*}
\Phi _{Z_i}(q,z,\hbar )=\displaystyle\int_{L_d^{\prime }}e^{pz}\displaystyle%
\sum_{d\geq 0}q^dE_d(p,\lambda _\alpha ,\mu _a,\hbar )\text{,} 
\end{equation*}
with $E_d(p,\lambda _\alpha ,\mu _a,\hbar )$ polynomials in $p$ of degree at
most $(n+1)d+n$ with polynomial coefficients in $(\lambda _\alpha ,\mu
_a,\hbar )$. We have pointed out two properties enjoyed by both $Z_i(q,\hbar
)$ and $Z_i^{*}(q,\hbar ).$ This leads us to define a particular class of
solutions of the ``{\it Calabi-Yau relations''}.

\begin{definition}
{\bf \ }A solution $W_i(q,\hbar )$ of the ``{\it Calabi-Yau relations}''
belongs to ${\cal P}$ if $\Phi _{W_i}(q,z,\hbar )$ can be written as a sum
of residues of polynomials $P_d(p,\lambda _\alpha ,\mu _a,\hbar )$ in $p$ of
degree at most $(n+1)d+n$ with polynomial coefficients in $(\lambda _\alpha
,\mu _a,\hbar )$.
\end{definition}

\noindent Notice that if $P_d(p,\lambda _\alpha ,\mu _a,\hbar )$ exists, then it is
uniquely determined. Indeed, the following holds.

\begin{proposition}
\label{chissaseva}Let $W_i(q,\hbar )$ be as in Proposition \ref{coeff}. The
polynomial coefficients $P_d$ in $\Phi _{W_i}$ are uniquely determined by
their values 
\[
P_d(\lambda _i+k\hbar )=\prod_{a=1}^r(l_a\lambda _i-\mu _a)F_k^{\left(
i\right) }(\hbar )F_{d-k}^{\left( i\right) }\left( -\hbar \right) \text{, }%
0\leq i\leq n\text{, }0\leq k\leq d\text{.}
\]
\end{proposition}

{\bf Proof. }If we compute the integral in $\Phi _{W_i}$ and compare with 
\ref{fizeta}, then 
\begin{equation*}
\frac{P_d(\lambda _i+k\hbar )}{\underset{(\alpha ,s)\neq (i,k)}{%
\prod_{\alpha =0}^n\prod_{s=0}^d}(\lambda _i-\lambda _\alpha +k\hbar -s\hbar
)}=\frac{\prod_{a=1}^r(l_a\lambda _i-\mu _a)F_k^{\left( i\right) }(\hbar )F_{d-k}^{\left( i\right) }\left(
-\hbar \right) }{\prod_{i\neq \alpha }(\lambda _i-\lambda _\alpha )k!\hbar
^k\prod_{\alpha \neq i}\prod_{s=1}^k(\lambda _i-\lambda _\alpha +s\hbar )}%
\cdot 
\end{equation*}
\begin{equation*}
\cdot \frac 1{(d-k)!(-\hbar )^{d-k}\prod_{i\neq \alpha
}\prod_{s=1}^{d-k}(\lambda _i-\lambda _\alpha -s\hbar )}\text{.} 
\end{equation*}

\noindent Since 
\begin{equation*}
\underset{(\alpha ,s)\neq (i,k)}{\prod_{\alpha =0}^n\prod_{s=0}^d}(\lambda
_i-\lambda _\alpha +k\hbar -s\hbar )=\underset{}{\underset{}{\prod_{\alpha
\neq i}^n}\prod_{\underset{s\neq k}{s=0}}^d}(\lambda _i-\lambda _\alpha
+k\hbar -s\hbar )\prod_{\underset{s\neq k}{s=0}}^d(k-s)\hbar \prod_{\alpha
\neq i}(\lambda _i-\lambda _\alpha )= 
\end{equation*}
\begin{equation*}
=\prod_{i\neq \alpha }(\lambda _i-\lambda _\alpha )k!\hbar ^k\prod_{\alpha
\neq i}\prod_{s=1}^k(\lambda _i-\lambda _\alpha +s\hbar )(d-k)!(-\hbar
)^{d-k}\prod_{i\neq \alpha }\prod_{s=1}^{d-k}(\lambda _i-\lambda _\alpha
-s\hbar )\text{,} 
\end{equation*}
the Proposition follows.
\begin{flushright}
$\Box$
\end{flushright}

It is possible to prove a result of uniqueness for functions in the class $%
{\cal P}$.

\begin{proposition}
A solution from ${\cal P}$ is uniquely determined by the two first terms in
the expansion of $W_i(q,\hbar )=W_i^{\left( 0\right) }+W_i^{\left( 1\right)
}\frac 1\hbar +\ldots $ as power series in $\frac 1\hbar $.
\end{proposition}

{\bf Proof. }If we write $R_{i,d}=R_{i,d}^{\left( 0\right) }\hbar
^d+R_{i,d}^{\left( 1\right) }\hbar ^{d-1}+\ldots ,$ then 
\begin{equation*}
W_i^{\left( 0\right) }=1+\sum_dR_{i,d}^{\left( 0\right) }q^d\text{,} 
\end{equation*}
\begin{equation*}
W_i^{\left( 1\right) }=\sum_dR_{i,d}^{\left( 1\right) }q^d\text{.} 
\end{equation*}

\noindent Fix $d\geq 0$ and suppose that $W_i(q,\hbar )$ and $W_i^{\prime }(q,\hbar )$
are two solutions from ${\cal P}$ with $W_i^{\left( 0\right) }=W_i^{^{\prime
}\left( 0\right) }$ and $W_i^{\left( 1\right) }=W_i^{^{\prime }\left(
1\right) }$. By induction it is possible to prove that $P_k=P_k^{\prime }$
for each $k<d$ and thus $P_d$ vanishes at $p=\lambda _i+k\hbar $, $0\leq
i\leq n$, $1\leq k\leq d-1$. This implies that $\prod_j\prod_{m=1}^{d-1}(p-%
\lambda _j-m\hbar )$ divides $P_d-P_d^{\prime }$. In addition, we also have
(see Proposition \ref{chissaseva}) 
\begin{equation*}
P_d(\lambda _i+d\hbar )-P_d^{\prime }(\lambda _i+d\hbar
)=\prod_{a=1}^r(l_a\lambda _i-\mu _a)F_d^{(i)}(\hbar )= 
\end{equation*}
\begin{equation*}
=\prod_{a=1}^r(l_a\lambda _i-\mu _a)d!\hbar ^d\prod_{j\neq
i}\prod_{m=1}^d(\lambda _i-\lambda _j-m\hbar )\left\{ K_i(d,\hbar
)-K_i^{\prime }(d,\hbar )\right\} = 
\end{equation*}

\begin{equation*}
=\prod_{a=1}^r(l_a\lambda _i-\mu _a)\prod_{j\neq i}\prod_{m=1}^d(\lambda
_i-\lambda _j-m\hbar )\left\{ R_{i,d}-R_{i,d}^{^{\prime }}\right\} \text{.} 
\end{equation*}

\noindent So the polynomial $\left\{ R_{i,d}-R_{i,d}^{^{\prime }}\right\} $ should be
divisible by $\hbar ^{d-1}$ and, by our hypotheses, it is identically zero.
Thus if two solutions in ${\cal P}$ coincide in orders $\hbar ^0$ and $\hbar
^{-1}$, then the very solutions coincide.
\begin{flushright}
$\Box$
\end{flushright}

\subsection{Transformations in ${\cal P}$}

At this point we come to describe transformations among solutions in ${\cal P}
$ .

\noindent Consider the following operations on the class ${\cal P}$:

(a) simultaneous multiplication $W_i(q,\hbar )\mapsto \nu _1(q)W_i(q,\hbar )$
by a power series of $q$ with $\nu _1(0)=1$;

(b) changes $W_i(q,\hbar )\mapsto e^{\frac{\lambda _i\nu _2(q)}\hbar
}W_i(qe^{\nu _2(q)},\hbar )$ with $\nu _1(0)=0$;

(c) multiplication $W_i(q,\hbar )\mapsto \exp (C\nu _3(q)/\hbar )W_i(q,\hbar
)$, where $C$ is a linear function of $(\lambda ,\mu )$ and $\nu _3(0)=0$.

\noindent If we apply these operations to $W_i(q,\hbar )$, then we have to check that $%
\Phi _i$ transforms to an equivariant integral of power series of
polynomials. To this purpose, we need two lemmas.

\begin{lemma}
Suppose that a series 
\[
s=\sum_dq^dP_d(p,\lambda _\alpha ,\mu _a,\hbar )
\]
with coefficients which are polynomials of $p$ of degree $\leq $dim$%
L_d^{\prime }$, has the property that for every $k=0,1,\ldots ,$ the $q$
series $\int_{L_d^{\prime }}sp^k$ has polynomial coefficients in $(\lambda
_\alpha ,\mu _a,\hbar )$. Then the coefficients of all $P_d$ are polynomials
of $(\lambda _\alpha ,\mu _a,\hbar )$ and vice versa.
\end{lemma}

{\bf Proof.} Decomposing $P_d(p,\lambda _\alpha ,\mu _a,\hbar )=\sum
P_{dj}(\lambda _\alpha ,\mu _a,\hbar )p^j$ we get that 
\begin{equation*}
\int_{L_d^{\prime }}sp^k=\sum_dq^d\sum_jP_{dj}(\lambda _\alpha ,\mu _a,\hbar
)\frac 1{2\pi \sqrt{-1}}\int \frac{p^{j+k}}{\prod_{\alpha
=0}^n\prod_{m=0}^d(p-\lambda _\alpha -m\hbar )}dp\text{;} 
\end{equation*}
the integral vanishes  for $j+k<(d+1)(n+1)-1$, equals $1$ for $j+k=(d+1)(n+1)-1$
and is a polynomial in $\lambda _\alpha ,\hbar$ for $j+k>(d+1)(n+1)-1$.
The resulting triangular matrix is invertible, and still the entries of the
inverse matrix are polynomials in $\lambda _\alpha ,\hbar$.
Hence, on varying $k$,
we obtain that every $P_{dj}$ is a polynomial.
\begin{flushright}
$\Box$
\end{flushright}

\begin{lemma}
\label{jemofatta}The polynomiality of $P_d$ in $\Phi _{W_i}$ is invariant
with respect to the operations (a), (b), and (c).
\end{lemma}

{\bf Proof.} It is clear that multiplication by a series in $q$ does not
change the polynomiality property. If we perform a (b)-type operation and
set $Q=qe^{\nu _2(q)}$, then 
\begin{equation*}
\sum_{i=0}^n\frac{\prod_{a=1}^r(l_a\lambda _i-\mu _a)}{\prod_{j\neq
i}(\lambda _i-\lambda _j)}e^{z\lambda _i}W_i(Qe^{z\hbar -\nu _2(q)+\nu
_2(qe^{z\hbar })},\hbar )W_i(Q,-\hbar )= 
\end{equation*}
\begin{equation*}
=\int_{L_d^{\prime }}e^{\frac p\hbar \left\{ ^{z\hbar -\nu _2(q)+\nu
_2(qe^{z\hbar })}\right\} }\sum_dq^de^{d\nu _2(q)}P_d(p,\lambda _\alpha ,\mu
_a,\hbar )\text{.} 
\end{equation*}

\noindent By the previous Lemma, we have to prove that for all $k$ the $q$ series $%
\left( \frac \partial {\partial z}\right) _{\mid z=0}^k\Phi _{W_i}$ have
polynomial coefficients. Since the exponent $z\hbar -\nu _2(q)+\nu
_2(qe^{z\hbar })$ is divided by $\hbar $, the derivatives still have
polynomial coefficients. The case of operation (c) is analogous with $\exp
\left[ C\frac{-\nu _3(q)+\nu _3(qe^{z\hbar })}\hbar \right] $.
\begin{flushright}
$\Box$
\end{flushright}

\begin{theorem}
The class ${\cal P}$ is invariant with respect to the operations (a), (b),
and (c).
\end{theorem}

{\bf Proof.} Consider the term $\frac{q^d}{\hbar ^d}coeff_i^j(d)$ in the
recursive relations. Application of the operations (a), (b), and (c) to the
left- and right-hand side of the ``{\it Calabi-Yau relations''} causes,
respectively, the following modifications in this coefficient:

\begin{equation*}
q^d\mapsto q^d\frac{f(q)}{f\left( Q\right) }\text{,} 
\end{equation*}
\begin{equation*}
q^d\mapsto q^d\exp \left\{ \frac{\lambda _if(Q\hbar )}\hbar -df(Q\hbar )-%
\frac{\lambda _if(Q(\lambda _i-\lambda _j)/d)}{(\lambda _i-\lambda _j)/d}%
\right\} \text{,} 
\end{equation*}
\begin{equation*}
q^d\mapsto q^d\exp \left\{ C\frac{f(Q\hbar )}\hbar -C\frac{f(Q(\lambda
_j-\lambda _i)/d)}{(\lambda _j-\lambda _i)/d}\right\} \text{.} 
\end{equation*}

\noindent In the case of the change (b), additionally, the argument $Q:=(\lambda
_j-\lambda _i)/d$ in $Z_i$, on the right-hand side of the recursive relation
gets an extra-factor $$\exp \left[ f(Q\hbar )-f(Q(\lambda _j-\lambda
_i)/d\right] $$. All the modifying factors above take 1 at $\hbar =(\lambda
_j-\lambda _i)/d$. This means that the term responsible for the simple
fraction with the pole at $\hbar =(\lambda _j-\lambda _i)/d$ does not
change, and that the operations modify only the polynomial initial
conditions. Under our assumptions about $\nu _1$, $\nu _2$, and $\nu _3$,
the modifying factors depend on $\hbar $ only in the combination $Q\hbar $.
This implies that the degrees of the new initial conditions still do not
exceed $d.$

\noindent The Theorem is completely proved by Lemma \ref{jemofatta}.
\begin{flushright}
$\Box$
\end{flushright}

\subsection{Initial conditions}

Let us determine initial conditions for $Z_i(q,\hbar )$ and $Z_i^{*}(q,\hbar
)$.

\begin{proposition}
$Z_i^{\left( 0\right) }=1$, $Z_i^{\left( 1\right) }=0$.
\end{proposition}

{\bf Proof. }The statements follows from the recursive relations, since the
only term in degree zero is $1$ and the contribution in degree $1$ is given
by the series 
\begin{equation*}
\sum_{d>0}q^d\int_{\overline{{\cal M}}_{0,2}({\Bbb P}^n,d)}{\cal E}%
_{2,d}^{\prime }e_0^{*}(\phi _i)\text{,} 
\end{equation*}
which is identically zero for dimension computations.
\begin{flushright}
$\Box$
\end{flushright}

\begin{proposition}
{\bf \ } 
\[
\begin{array}{c}
Z_i^{*\left( 0\right) }=\sum_{d\geq 0}\displaystyle\frac{(l_1d)!\ldots
(l_rd)!}{(d!)^{n+1}}q^d:=f(q)\text{,} \\ 
Z_i^{*\left( 1\right) }=\lambda _i\sum_{a=1}^rl_a\left[
g_{l_a}(q)-g_1(q)\right] +\left( \sum_{\alpha =0}^n\lambda _\alpha \right)
g_1(q)-\sum_{a=1}^r\mu _ag_{l_a}(q)\text{,}
\end{array}
\]
\end{proposition}
where 
\begin{equation*}
g_l=\sum_{d\geq 1}q^d\frac{\prod_{a=1}^r(l_ad)!}{(d!)^{n+1}}\left(
\sum_{m=1}^{ld}\frac 1m\right) \text{.} 
\end{equation*}

{\bf Proof. }If we set $\omega =1/\hbar $, then 
\begin{equation*}
C_i^{*}(d,\omega )=\frac{\prod_{a=1}^r\prod_{m=1}^{l_ad}\left[ (l_a\lambda
_i-\mu _a)\omega +m\right] }{\prod_{\alpha =0}^n\prod_{m=1}^d\left[ (\lambda
_i-\lambda _\alpha )\omega +m\right] }=\frac{F\left( \omega \right) }{%
G\left( \omega \right) }, 
\end{equation*}
and so the statement follows easily.
\begin{flushright}
$\Box$
\end{flushright}

Finally, we can prove the following

\begin{theorem}
\label{cestfini}{\bf \ }Let $f(q)$ and $g_l(q)$ be as in Proposition 1.6. If
we perform the following operations with $Z_i(q,\hbar )$: 
\[
\]
(1) put 
\[
Q=q\exp \left\{ \sum_{a=1}^rl_a\left[ g_{l_a}(q)-g_1(q)\right] /f(q)\right\} 
\text{,}
\]
(2) multiply $Z_i(Q(q),\hbar )$ by 
\[
\exp \left\{ \frac 1{f(q)\hbar }\left[ \sum_{a=1}^r(l_a\lambda _i-\mu
_a)g_{l_a}(q)-(\sum_{\alpha =0}^n(\lambda _i-\lambda _\alpha ))g_1(q)\right]
\right\} ,
\]
(3) multiply all $Z_i(q,\hbar )$ simultaneously by $f(q)$. 
\[
\]
Then the resulting functions coincide with functions $Z_i^{*}(q,\hbar )$.
\end{theorem}

{\bf Proof. }The operations (1), (2) and (3) correspond to consecutive
applications of operations of type (a), (b) and (c). Indeed, (3) is exactly
(a) with $f(q)=\nu _1(q)$. In addition, if we set 
\begin{equation*}
\nu _2(q)=\left\{ \sum_{a=1}^rl_a\left[ g_{l_a}(q)-g_1(q)\right]
/f(q)\right\} \text{,} 
\end{equation*}

\noindent then the change of variable in (1) is one of those allowed by operations of
type (b). Finally, since $\sum_{a=1}^rl_a=n+1$, 
\begin{equation*}
\frac 1{f(q)\hbar }\left[ \sum_{a=1}^r(l_a\lambda _i-\mu
_a)g_{l_a}(q)-(\sum_{\alpha =0}^n(\lambda _i-\lambda _\alpha ))g_1(q)\right]
= 
\end{equation*}
\begin{equation*}
=\nu _2(q)-(\sum_{a=1}^r\mu _a)\frac{g_{l_a}(q)}{f(q)\hbar }+(\sum_{\alpha
=0}^n\lambda _\alpha )\frac{g_1(q)}{f(q)\hbar }\text{,} 
\end{equation*}
and thus (2) corresponds to operations of type (c). (Notice that $f(q)$ is
invertible).

\noindent It is straightforward to verify that the initial conditions of $Z_i(q,\hbar
) $ are transformed into those of $Z_i^{*}(q,\hbar )$. According to
Proposition 1.10. this transforms functions $Z_i(q,\hbar )$ to $%
Z_i^{*}(q,\hbar )$.
\begin{flushright}
$\Box$
\end{flushright}

\medskip\ 

\subsection{Proof of the ``Mirror Conjecture''}

The main and last theorem of this section is

\begin{theorem}
\label{calteo}Suppose $l_1+\ldots +l_r=n+1$. Consider the functions 
\[
f(e^t):=\sum_{d\geq 0}\frac{(l_1d)!\ldots (l_rd)!}{(d!)^{n+1}}e^{td}
\]
and 
\[
g_s(e^t):=\sum_{d\geq 1}e^{td}\left( \sum_{m=1}^{sd}\frac 1m\right) \frac{%
\prod_{a=1}^r(l_ad)!}{(d!)^{n+1}}.
\]
After performing the change of variable 
\[
T=t+\sum_{a=1}^r\frac{l_a\left[ g_{l_a}(e^t)-g_1(e^t)\right] }{f(e^t)}\text{,%
}
\]
and multiply the function $S_{\left( 2\right) }(T,\hbar )$ by $f(e^t)$, we
obtain the function
\end{theorem}

\begin{equation*}
S_{(2)}^{*}(t,\hbar )=e^{\frac{Pt}\hbar }\sum_{d\geq 0}e^{dt}\frac{%
\prod_{a=1}^r\prod_{m=0}^{l_ad}(l_aP+m\hbar )}{\prod_{m=1}^d(P+m\hbar )^{n+1}%
}\text{.} 
\end{equation*}

{\bf Proof. }By Theorem \ref{cestfini}, components of the equivariant
counterpart of $S_{\left( 2\right) }(T,\hbar )$, i.e. 
\begin{equation*}
s_i^{^{\prime }}(T,\hbar )=\prod_{j\neq i}(\lambda _j-\lambda _i)e^{\frac{%
\lambda _iT}\hbar }\prod_{a=1}^r(l_a\lambda _i-\mu _a)Z_i(e^T,\hbar )\text{,}
\end{equation*}
can be transformed to components 
\begin{equation*}
s_i^{^{\prime }*}(t,\hbar )=\prod_{j\neq i}(\lambda _j-\lambda _i)e^{\frac{%
\lambda _it}\hbar }\prod_{a=1}^r(l_a\lambda _i-\mu _a)Z_i^{*}(e^t,\hbar )%
\text{.} 
\end{equation*}

\noindent If we pass to equivariant to non-equivariant cohomolgy, the Theorem follows.
\begin{flushright}
$\Box$
\end{flushright} 

\begin{remark}
{\bf \ }Notice that the components $s_0^{*}(t,\hbar )$ and $s_1^{*}(t,\hbar )
$ in 
\[
S_{\left( 2\right) }^{*}(t,\hbar )=l_1\ldots l_r\left[ P^rs_0^{*}(t,\hbar
)+P^{r+1}s_1^{*}(t,\hbar )+\ldots +P^ns_n^{*}(t,\hbar )\right] 
\]
are exactly $f(e^t)$ and $tf\left( e^t\right) +\sum_{a=1}^rl_a\left[
g_{l_a}(e^t)-g_1(e^t)\right] $ respectively. Thus the inverse transformation
from $S_{\left( 2\right) }^{*}(t,\hbar )$ to $S_{\left( 2\right) }(T,\hbar )$
consists in division by $f(e^t)$ followed by the change $T=\frac{%
s_1^{*}(t,\hbar )}{s_0^{*}(t,\hbar )}$ \thinspace in accordance with the
method conjectured by physicists \cite{COGP}.
\end{remark}

\subsection{Proof{\bf \ }of the ``Main Lemma''}

Let 
\begin{equation*}
L_d^0:=\overline{{\cal M}}_{0,0}\left( {\Bbb P}^n\times {\Bbb P}^1,\left(
d,1\right) \right) \text{,} 
\end{equation*}
and let 
\begin{equation*}
\pi _0:L_d\rightarrow L_d^0 
\end{equation*}
be the map that forgets the two markings.

\noindent We shall construct a map 
\begin{equation*}
\mu _0:L_d^0\rightarrow L_d^{\prime }\text{,} 
\end{equation*}
and then compose to define the map $\mu :=\mu _0\pi _0$.

{\bf STEP I)} {\em Fixed points.\ }

\noindent We have already described fixed points in $L_d^0$ for the $S^1$ action (cfr.
Prop. \ref{fixs1}): 
\begin{eqnarray*}
\left( L_d^0\right) ^{S^1} &=&\bigcup_{k=0}^dM_{k,d-k}^0= \\
&=&\bigcup_{k=0}^d\left\{ 
\begin{array}{c}
\left[ C_1\cup C_0\cup C_\infty ,(\psi _1,\psi _2)\right] \in \left( \left(
Y\times {\Bbb P}^1\right) _{0,\left( d,1\right) }\right) : \\ 
\begin{array}{c}
\psi _{2\mid C_1}:C_1\overset{\symbol{126}}{\longrightarrow }{\Bbb P}^1,\psi
_2\left( C_0\right) =0,\psi _2\left( C_\infty \right) =\infty , \\ 
\deg \psi _{1\mid C_0}=k,\deg \psi _{1\mid C_\infty }=d-k,\psi
_1(C_1)=\left\{ pt.\right\}
\end{array}
\end{array}
\right\} \cong \\
&\cong &\bigcup_{k=0}^dY_{1,k}\times _YY_{1,d-k}\text{.}
\end{eqnarray*}

\noindent Among them, the fixed points for the $\left( S^1\right) ^{n+1}$ action are
given by 
\begin{equation*}
\bigcup_{i=0}^n\bigcup_{k=0}^dM_{k,d-k}^{0,i}\cong
\bigcup_{i=0}^n\bigcup_{k=0}^d\nu ^{-1}\left( p_i\right) \cap \left\{ \left(
Y_{1,k}\right) ^{\left( S^1\right) ^{n+1}}\times _Y\left( Y_{1,d-k}\right)
^{\left( S^1\right) ^{n+1}}\right\} \text{,} 
\end{equation*}
where $\nu $ is the evaluation map on the marked point, and $p_i$ is the
projectivization of the $i$-th coordinate line in ${\Bbb P}^n$. Analogously, 
\begin{equation*}
\left( L_d\right) ^{S^1}=\bigcup_{k=0}^dM_{k,d-k}\cong
\bigcup_{k=0}^dY_{2,k}\times _YY_{2,d-k}\text{,} 
\end{equation*}
and 
\begin{eqnarray*}
\left( L_d\right) ^{S^1\times \left( S^1\right) ^{n+1}}
&=&\bigcup_{i=0}^n\bigcup_{k=0}^dM_{k,d-k}^i\cong \\
&\cong &\bigcup_{i=0}^n\bigcup_{k=0}^d\nu ^{-1}\left( p_i\right) \cap
\left\{ \left( Y_{2,k}\right) ^{\left( S^1\right) ^{n+1}}\times _Y\left(
Y_{2,d-k}\right) ^{\left( S^1\right) ^{n+1}}\right\} \text{.}
\end{eqnarray*}

\noindent On the other hand, the fixed points in $L_d^{\prime }$ are $(n+1)$-tuples of
monomials of the following form: 
\begin{equation*}
\left[ 0,...,0,z_1^kz_2^{d-k},0,...,0\right] \text{, }0\leq k\leq d\text{.} 
\end{equation*}

\smallskip\ 

{\bf STEP 2)} {\em Definition and equivariancy of }$\mu _0$.

\noindent There is a natural way to define $\mu _0$ on the open set ${\cal M}%
_{0,0}\left( {\Bbb P}^n\times {\Bbb P}^1,\left( d,1\right) \right) \subset
L_d^0$. In fact, a point in this set is an equivalence class of maps ${\Bbb P%
}^1\rightarrow {\Bbb P}^n\times {\Bbb P}^1$ of bidegree $\left( d,1\right) $%
, and the image is the graph of a degree $d$ map ${\Bbb P}^1\rightarrow 
{\Bbb P}^n$, which is clearly given by a $(n+1)$-tuple of homogeneous
polynomials of degree $d$ in two variables, hence by a point in $L_d^{\prime
}$.

\noindent A generic point in $L_d^0$ is given by $\left[ C_0\cup C_1\cup ...\cup
C_r,(\psi _1,\psi _2)\right] $, with $\deg $ $\psi _{2\mid C_0}=1$, $\deg
\psi _{2\mid C_i}=0$, $i=1,...,r$, and $\deg $ $\psi _{1\mid C_j}=:d_j$, for
every $j$, with $\sum d_j=d$.

\noindent Let $q_i:=\psi _2\left( C_i\right) $, $i=1,...,r$, and pick a homogenous
polynomial $g$ of degree $d^{\prime }:=d-d_0$, in two variables, with a root
of order $d_i$ in each $q_i$ ; furthermore notice that the image of $\psi
_{2\mid C_0}$ is the graph of a degree $d_0$ map ${\Bbb P}^1\rightarrow 
{\Bbb P}^n$, hence is given by a $(n+1)$-tuple $\left[ f_0,...,f_n\right] $
of homogeneous polynomials of degree $d_0$.

Now we are ready to define $\mu _0\left( \left[ C_0\cup C_1\cup ...\cup
C_r,(\psi _1,\psi _2)\right] \right) :=\left[ gf_0,...,gf_n\right] $;
roughly speaking, the map $\mu _0$ forgets everything about the components $%
C_1,...,C_r$, except their image point in ${\Bbb P}^1$. Clearly, 
\begin{equation*}
\mu \left( \left[ C_0\cup C_1\cup ...\cup C_r,x_0,x_1,(\psi _1,\psi
_2)\right] \right) =\left[ gf_0,...,gf_n\right] . 
\end{equation*}

\noindent By this definition, the map is $S^1\times \left( S^1\right) ^{n+1}$
equivariant; notice that
\begin{equation*}
\mu \left( M_{k,d-k}^i\right) =\left[ 0,...,0,z_1^{d-k}z_2^{k},0,...,0\right] 
\text{.} 
\end{equation*}

{\bf STEP 2)} {\em Regularity of }$\mu _0$.

\noindent This is the most difficult step; we propose two different proofs of it; for
both, we need to recall that the moduli space of maps can be also
constructed (\cite{FP}) as the GIT quotient of a quasi-projective variety $%
{\cal J}$, such that there exist a semi-universal family ${\cal F}\overset{%
\pi }{\rightarrow }{\cal J}$, and the fiber of this map on each point is
isomorphic to the curve represented by the point itself.

\smallskip\ 

{\bf First proof.}

\noindent We consider the following incidence variety: 
\begin{equation*}
\Gamma =\left\{ \left( \left[ C,\psi _1\right] ,H,\left( p_1,...,p_d\right)
\right) \in {\cal J}\times \left( {\Bbb P}^n\right) ^{*}\times {\Bbb P}%
^{1^{(d)}}:\psi _1^{-1}\left( H\right) \cap C=\psi _2^{-1}\left(
p_1,...,p_d\right) \right\} \text{,} 
\end{equation*}
where ${\Bbb P}^{1^{(d)}}\cong {\Bbb P}^d$ is the $d$-fold symmetric product
of ${\Bbb P}^1$.

Let 
\begin{equation*}
{\cal J}\times \left( {\Bbb P}^n\right) ^{*}\times {\Bbb P}^{1^{(d)}}%
\overset{\gamma }{\rightarrow }{\cal J}\times \left( {\Bbb P}^n\right) ^{*}%
\overset{\delta }{\rightarrow }{\cal J} 
\end{equation*}
be the obvious projections. If we denote with $\Delta =\gamma (\Gamma )$,
the map $\delta $ restricted to $\Delta $ is still surjective, since every
degree $d$ curve in ${\Bbb P}^n$ intersects a generic hyperplane exactly in $%
d$ points; this allows us to say more, that is, the fiber $\delta
^{-1}\left( \left[ C,\psi _1\right] \right) $ in the following diagram 
\begin{equation*}
\Gamma \overset{\gamma }{\rightarrow }\Delta \overset{\delta }{\rightarrow }%
{\cal J} 
\end{equation*}
is $\left\{ \left[ C,\psi _1\right] \right\} \times U_{\left[ C,\psi
_1\right] }$, where $U_{\left[ C,\psi _1\right] }$ is the Zariski open
subset of hyperplanes which intersect transversally the curve $\psi _1\left(
C\right) $. Moreover, by definition of $\Gamma $, the fiber $\gamma
^{-1}\delta ^{-1}\left( \left[ C,\psi _1\right] \right) =\gamma ^{-1}\left(
\left\{ \left[ C,\psi _1\right] \right\} \times U_{\left[ C,\psi _1\right]
}\right) $ can be viewed as $\left\{ \left[ C,\psi _1\right] \right\} \times
G$, with $G$ the graph in $U_{\left[ C,\psi _1\right] }\times {\Bbb P}%
^{1^{(d)}}$of a regular map $g_{\left[ C,\psi _1\right] }:U_{\left[ C,\psi
_1\right] }\rightarrow {\Bbb P}^{1^{(d)}}$.

\noindent This map turns out to be linear, as it sends linear subspaces to linear
subspaces, and we can extend it to a linear map $g_{\left[ C,\psi _1\right]
}:{\Bbb P}^n\rightarrow {\Bbb P}^d$, thus we can close up $\Gamma $ to a
subvariety $\overline{\Gamma }$ whose fiber over each curve is the graph of
a linear map.

Consider the map 
\begin{equation*}
\begin{tabular}{ccc}
${\Bbb P}^n\times {\Bbb P}\left( Hom\left( {\Bbb C}^{n+1},{\Bbb C}%
^{d+1}\right) \right) $ & $\overset{\alpha }{\rightarrow }$ & ${\Bbb P}%
^n\times {\Bbb P}^d$ \\ 
$\left[ v\right] ,\left[ f\right] $ & $\rightarrow $ & $\left[ v\right]
,\left[ f\left( v\right) \right] $%
\end{tabular}
\end{equation*}

\noindent Our considerations imply that the subvariety $\widehat{\Gamma }:=\left(
\left( Id,\alpha \right) ^{-1}\left( \overline{\Gamma }\right) \right) \cap 
{\cal J}\times \left\{ pt\right\} \times {\Bbb P}\left( Hom\left( {\Bbb C}%
^{n+1},{\Bbb C}^{d+1}\right) \right) $ projects bijectively on the first
factor ${\cal J}$ ; since this variety is normal, $\widehat{\Gamma }$ is
isomorphic to the graph of a regular map 
\begin{equation*}
{\cal J}\rightarrow {\Bbb P}\left( Hom\left( {\Bbb C}^{n+1},{\Bbb C}%
^{d+1}\right) \right) \text{.} 
\end{equation*}

\noindent Now we can easily observe that this map coincides set-theoretically with our 
$\mu _0$, and we are done.

\smallskip\ 

{\bf Second proof.}

\noindent We are going to build a line bundle on ${\cal J}$ which naturally gives an
invariant map ${\cal J}\rightarrow L_d^{\prime }$, hence a map $%
L_d^0\rightarrow L_d^{\prime }$ , and then prove that coincides with our
set-theoretically defined $\mu _0.$ We also need the map $e:{\cal F}%
\rightarrow {\Bbb P}^n\times {\Bbb P}^1$ , the evaluation on the marked
point. Finally, let $H:=Hom\left( {\cal O}_{{\Bbb P}^n}\left( 1\right) ,%
{\cal O}_{{\Bbb P}^1}\left( d\right) \right) $ be a line bundle on ${\Bbb P}%
^n\times {\Bbb P}^1$ (for simplicity, we omit to write the pullbacks by the
projection maps). Consider the sheaf on ${\cal J}$ defined as follows: 
\begin{equation*}
{\cal H}^0\left( U\right) :=H^0\left( \pi ^{-1}\left( U\right) ,e^{*}\left(
H\right) \right) \text{.} 
\end{equation*}

\noindent We claim:

\begin{enumerate}
\item  \label{claim1}${\cal H}^0$ is a rank one locally free sheaf;

\item  \label{claim2}the fiber at $\left[ C=C_0\cup C_1\cup ...\cup
C_r,(\psi _1,\psi _2)\right] $ of the corresponding line bundle can be
identified with $H^0\left( C_0,\psi _{\mid C_o}^{*}\left( H\right)
\otimes {\cal O}\left( -q_1\right) ^{d_1}\otimes ...\otimes {\cal O}%
\left( -q_r\right) ^{d_r}\right) $.
\end{enumerate}

\noindent We postpone the proof of these facts.

\noindent A non zero vector $u_C$ in the above fiber gives a natural map 
\begin{equation*}
f_C:H^0\left( C,\psi _1^{*}\left( {\cal O}_{{\Bbb P}^n}\left( 1\right)
\right) \right) \rightarrow H^0\left( C,\psi _2^{*}\left( {\cal O}_{{\Bbb P}%
^1}\left( d\right) \right) \right) =H^0\left( {\Bbb P}^1,\left( {\cal O}_{%
{\Bbb P}^1}\left( d\right) \right) \right) \text{;} 
\end{equation*}
in fact, $u_{C\text{ }}$ is represented by a class $\left[ \sigma \right]
\in {\cal H}_{\left[ C,\psi \right] }^0$; the restriction of $\sigma $ to $%
\pi ^{-1}(\left[ C,\psi \right] )$ is well defined and does not depend on
the choice of the element of the class; hence it is well defined 
\begin{equation*}
\left[ \sigma \right] \rightarrow \sigma _{\mid \left[ C,\psi \right] }\in
H^0\left( \pi ^{-1}\left( \left[ C,\psi \right] \right) ,e^{*}\left(
H\right) \right) \text{.} 
\end{equation*}

\noindent Since $\pi ^{-1}\left( \left[ C,\psi \right] \right) \cong C$, we can see $%
\sigma _{\mid \pi ^{-1}(\left[ C,\psi \right] )}$ as a section $\tau _C\in
H^0\left( C,\psi ^{*}\left( H\right) \right) $ $=H^0\left( C,\psi ^{*}\left(
Hom\left( {\cal O}_{{\Bbb P}^n}\left( 1\right) ,{\cal O}_{{\Bbb P}^1}\left(
d\right) \right) \right) \right) $; the map $f_C$ is constructed via $\tau
_C $.

\noindent Moreover, it will be clear after the proof of the claim that the kernel of
the map $f_{C\text{ }}$ consists exactly of sections vanishing identically
on $C_0$; therefore we can proceed as follows: pick a basis $\left\{
X_0,...,X_n\right\} $ of $H^0\left( {\Bbb P}^n,{\cal O}_{{\Bbb P}^n}\left(
1\right) \right) $, pull it back to a set of elements in $H^0\left( C,\psi
_1^{*}\left( {\cal O}_{{\Bbb P}^n}\left( 1\right) \right) \right) $, which
cannot vanish simoultaneously on $C_0$; the map 
\begin{eqnarray*}
{\cal H}_{\left[ C,\psi \right] }^0\backslash \left\{ 0\right\} &\rightarrow
&\left( {\Bbb C}^{n+1}\otimes H^0\left( {\Bbb P}^1,{\cal O}\left(
d\right) \right) \right) \backslash \left\{ 0\right\} \\
u_C &\rightarrow &f_C\left( \psi _1^{*}\left( X_0\right) \right)
,...,f_C\left( \psi _1^{*}\left( X_n\right) \right)
\end{eqnarray*}
extends to a map 
\begin{equation*}
\left( {\cal H}^0\right) ^{*}\rightarrow {\Bbb C}^{n+1}\otimes H^0\left( 
{\Bbb P}^1,{\cal O}\left( d\right) \right) \text{,} 
\end{equation*}
which is equivariant with respect to translation on the fiber on the LHS,
and multiplication by a constant on the RHS. This map induces 
\begin{equation*}
{\Bbb P}\left( {\cal H}^0\right) \cong {\cal J}\rightarrow {\Bbb P}\left( 
{\Bbb C}^{n+1}\otimes H^0\left( {\Bbb P}^1,{\cal O}\left( d\right)
\right) \right) \cong L_d^{\prime }\text{, } 
\end{equation*}
which will be exactly our $\mu _0$.

\smallskip\ 

In order to justify claims \ref{claim1} and \ref{claim2}, we need to compute
the space of global sections of the sheaf $e^{*}\left( H\right) $ in a
formal neighborhood $\pi ^{-1}\left( U\right) $ of the fiber $\pi
^{-1}\left( \left[ C,\psi \right] \right) $ of the forgetful map $\pi :{\cal %
F}\rightarrow L_d^0$ .

{\bf A) }$C${\em \ is non singular. }In this case, we can find a formal
neighborhood 
\begin{equation*}
\pi ^{-1}\left( U\right) \cong C\times \Delta _\epsilon , 
\end{equation*}
with $\Delta _\epsilon = Spec{\Bbb C} \left[ \left[ \epsilon \right] \right] $.

We can prove that 
\begin{equation*}
{\cal H}^0\left( U\right) =H^0\left( \pi ^{-1}\left( U\right) ,e^{*}\left(
H\right) \right) \cong {\Bbb C}\left[ \left[ \epsilon \right] \right]
\otimes H^0\left( C,\psi ^{*}\left( H\right) \right) ; 
\end{equation*}
in fact, we can build an obvious covering of $\pi ^{-1}\left( U\right) $
trivializing the bundle, namely ${\cal W}=\left\{ W_1,W_2\right\} $, where 
\begin{equation*}
W_i\cong {\Bbb C}\times \Delta _\epsilon 
\end{equation*}
with coordinates $(x_i,\epsilon )$ related by $x_1=x_2^{-1}$. Let $%
g_{12}\left( x_2,\epsilon \right) $ be the transition function of the bundle 
$e^{*}\left( H\right) $ with respect to this covering; notice that 
\begin{equation*}
H^0(\pi ^{-1}\left( \left[ C,\psi \right] \right) ,e^{*}\left( H\right)
)\cong H^0(C,\psi ^{*}\left( H\right) )={\cal O}\text{,} 
\end{equation*}
hence $g_{12}\left( x_2,0\right) =k$ and $g_{12}\left( x_2,\epsilon \right)
=k\left( 1+\epsilon f_{12}\left( x_2,\epsilon \right) \right) $, $k\in {\Bbb %
C}$, is an invertible function .

\noindent A section in $H^0\left( \pi ^{-1}\left( U\right) ,e^{*}\left( H\right)
\right) $ is given by $\left\{ \sigma _1\left( x_1,\epsilon \right) ,\sigma
_2\left( x_2,\epsilon \right) \right\} $ with 
\begin{equation*}
\sigma _1\left( x_1,\epsilon \right) =g_{12}\left( x_2,\epsilon \right)
\sigma _2\left( x_2,\epsilon \right) =k\left( 1+\epsilon f_{12}\left(
x_2,\epsilon \right) \right) \sigma _2\left( x_2,\epsilon \right) \text{;} 
\end{equation*}
thus we can set $\tau _2\left( x_2,\epsilon \right) :=\left( 1+\epsilon
f_{12}\left( x_2,\epsilon \right) \right) \sigma _2\left( x_2,\epsilon
\right) $, and we can establish a bijection between couples $\left\{ \sigma
_1,\sigma _2\right\} $ and couples $\left\{ \sigma _1,\tau _2\right\} $.

\noindent By developing in power series with respect to $\epsilon $, we see that the
latter one is an element of ${\Bbb C}\left[ \left[ \epsilon \right] \right]
\otimes H^0\left( C,\psi ^{*}\left( H\right) \right) \cong {\Bbb C}\left[
\left[ \epsilon \right] \right] $. The map $f_{C\text{ }}$ is consequently a
multiplication by a constant. In particular $f_C\left( \psi ^{*}\left(
X_i\right) \right) =kf_i$, as we required.

{\bf B) }$C$ {\em is reducible.} For the sake of simplicity, suppose there
are two irreducible components, $C=C_0\cup C_1$. Once more $\pi ^{-1}\left(
\left[ C,\psi \right] \right) \cong C$; if $q$ is the node, we choose a
formal neighborhood $\pi ^{-1}\left( U\right) $ with the following
trivializing covering for $e^{*}\left( H\right) $: 
\begin{equation*}
W\cong \Delta _{x_0,x_1,\epsilon } 
\end{equation*}
is a neighborhood of $q$, and the curve has equation $x_0x_1=0,\epsilon =0$; 
\begin{equation*}
W_0\cong {\Bbb C}\times \Delta _{y_0,\eta ,\epsilon } 
\end{equation*}
is a neighborhood of $C_0\backslash \left\{ q\right\} $, and the curve has
equation $\eta =\epsilon =0$, and 
\begin{equation*}
W_1\cong {\Bbb C}\times \Delta _{y_1,\eta ,\epsilon } 
\end{equation*}
is a neighborhood of $C_1\backslash \left\{ q\right\} $, and the curve has
equation $\eta =\epsilon =0$; moreover, we have transition functions $%
y_0=x_0^{-1},\eta =x_0x_1$ from $W$ to $W_0$, and $y_1=x_1^{-1},\eta =x_0x_1$
from $W$ to $W_1$.

\noindent In order to see how $e^{*}\left( H\right) _{\mid \pi ^{-1}\left( \left[
C,\psi \right] \right) }$ looks like, it is useful to observe that 
\begin{equation*}
\psi _{\mid C_0}^{*}\left( Hom\left( {\cal O}_{{\Bbb P}^n}\left( 1\right) ,%
{\cal O}_{{\Bbb P}^1}\left( d\right) \right) \right) \cong {\cal O}\left(
d-d_0\right) \text{,} 
\end{equation*}
and 
\begin{equation*}
\psi _{\mid C_1}^{*}\left( Hom\left( {\cal O}_{{\Bbb P}^n}\left( 1\right) ,%
{\cal O}_{{\Bbb P}^1}\left( d\right) \right) \right) \cong {\cal O}\left(
-d_1\right) \text{.} 
\end{equation*}

\noindent Let $\left\{ \sigma \left( x_0,x_1,\epsilon \right) ,\sigma _0\left(
y_0,\eta ,\epsilon \right) ,\sigma _1\left( y_1,\eta ,\epsilon \right)
\right\} $ be a section in $H^0\left( \pi ^{-1}\left( U\right) ,e^{*}\left(
H\right) \right) $, with 
\begin{equation*}
\sigma \left( x_0,x_1,\epsilon \right) =g_1\left( y_1,\eta ,\epsilon \right)
\sigma _1\left( y_1,\eta ,\epsilon \right) \text{;} 
\end{equation*}
since $\psi ^{*}\left( H\right) _{\mid C_1}\cong {\cal O}\left( -d_1\right) $%
, $g_1\left( y_1,\eta ,\epsilon \right) =y_1^{d_1}f_1\left( y_1,\eta
,\epsilon \right) $, with $f_1\left( y_1,\eta ,\epsilon \right) $
invertible. We obtain $\sigma \left( x_0,x_1,\epsilon \right)
=y_1^{d_1}f_1\left( y_1,\eta ,\epsilon \right) \sigma _1\left( y_1,\eta
,\epsilon \right) $; observe that we should have $\sigma _1\left( y_1,\eta
,\epsilon \right) =\eta ^{d_1}\tau _1\left( y_1\eta ,\epsilon \right) $,
otherwise, since $y_1=x_1^{-1}$, $\sigma $ could not be a polynomial of
positive degree in $x_1$, hence $f_1^{-1}\sigma =y_1^{d_1}\eta ^{d_1}\tau
_1\left( y_1\eta ,\epsilon \right) =x_0^{d_1}\tau _1\left( x_0,\epsilon
\right) $; this means that from $f_1^{-1}\sigma $, which is a function on $W$
with a zero of order $\geq d_1$ in $q$, we recover $\sigma _1$.

Now we proceed as in {\bf A)} with the couple $\left\{ \sigma ,\sigma
_0\right\} $, related by a transition function $g_0\left( y_0,\eta ,\epsilon
\right) =y_0^{-d_0}f_0\left( y_0,\eta ,\epsilon \right) $, and prove that
there is a bijection between triples $\left\{ \sigma ,\sigma _0,\sigma
_1\right\} $ which give a section in $H^0\left( \pi ^{-1}\left( U\right)
,e^{*}\left( H\right) \right) $, and couples $\left\{ \sigma ,f_0\sigma
_0\right\} $ which give rise to an element of 
\begin{equation*}
{\Bbb C}\left[ \left[ \epsilon \right] \right] \otimes H^0\left( C_0,\psi
_{\mid C_0}^{*}\left( H\right) \otimes {\cal O}\left( -q\right)
^{d_1}\right) 
\end{equation*}
thus 
\begin{equation*}
{\cal H}^0\left( U\right) \cong {\Bbb C}\left[ \left[ \epsilon \right]
\right] \otimes H^0\left( C_0,\psi _{\mid C_0}^{*}\left( H\right)
\otimes {\cal O}\left( -q\right) ^{d_1}\right) \text{.} 
\end{equation*}
Since $\dim H^0\left( C,\psi _{\mid C_0}^{*}\left( H\right) \otimes {\cal %
O}\left( -q\right) ^{d_1}\right) =1$, this implies our claims in this case.

Moreover, the couple $\left\{ \sigma ,f_0\sigma _0\right\} $, evaluted in $%
\epsilon =0$ should give a homogeneous polynomial of degree $d_1=d-d_0$ on $%
C_0$, but we already proved that it has a zero of order $d_1$ in $q$, and
this determines it completely; if $g$ is this polynomial, then in this case
the map $f_C$ consists of multiplication by $g$, and $f_C\left( \psi
^{*}\left( X_i\right) \right) =gf_i$; an induction procedure on the number
of irreducible components completes the argument. Note that the map
coincides on points which correspond to isomorphic stable maps, hence
descends to a map on $L_d^0$, and that it is exactly $\mu _0$.

\smallskip\ 

{\bf STEP 3)} Let us look at the expression of $\sum_{d\geq 0}\exp \left(
d\tau \right) \int_{L_d}\exp \left( \frac{\mu ^{*}\left( p\right) \left(
t-\tau \right) }\hbar \right) {\cal E}_{2,d}$. The bijection established at
the end of {\bf STEP 1)} between connected components of fixed points in $%
L_d $ and $L_d^{\prime }$ tells us that the localization of $\mu ^{*}\left(
p\right) $ at $M_{k,d-k}^i$ equals $\lambda _i+k\hbar $, and thus the
pullback of $\mu ^{*}\left( p\right) $ to the fixed point set 
\begin{equation*}
M_{k,d-k}:=\bigcup_{i=0}^nM_{k,d-k}^i\cong \left\{ \left( Y_{2,k}\right)
^{\left( S^1\right) ^{n+1}}\times _Y\left( Y_{2,d-k}\right) ^{\left(
S^1\right) ^{n+1}}\right\} 
\end{equation*}
coincides with the $\left( S^1\right) ^{n+1}$equivariant class $\nu
^{*}\left( p+k\hbar \right) $. Hence, by computations of section \ref{separate},
taking in account that the action of $S^{1}$ is slightly different,
and recalling conventions for degenerate cases, the following holds:

\begin{equation*}
\sum_{d\geq 0}\exp \left( d\tau \right) \int_{L_d}\exp \left( \frac{\mu
^{*}\left( p\right) \left( t-\tau \right) }\hbar \right) {\cal E}_{2,d}= 
\end{equation*}
\begin{eqnarray*}
\ &=&\sum_{d\geq 0}\exp \left( d\tau \right) \sum_{k=0}^d\int_{M_{k,d-k}}%
\frac{\exp \left( \frac{\nu ^{*}\left( p+kh\right) \left( t-\tau \right) }%
\hbar \right) }{{\cal N}_{M_{k,d-k}/L_d}}{\cal E}_{2,d}= \\
\ &=&\sum_{d\geq 0}\exp \left( d\tau \right)
\sum_{k=0}^d\sum_{i=0}^n\frac{\prod_a\left( l_a\lambda _i-\mu _a\right) }{%
\prod_{j\neq i}\left( \lambda _j-\lambda _i\right) }\cdot
\end{eqnarray*}
\begin{equation*}
\cdot \left( \exp \left( kt\right) \int_{Y_{2,k}}\frac{e_0^{*}\left( \phi
_i\right) \exp \left( \frac{pt}\hbar \right) }{\hbar +c}{\cal E}%
_{2,d}^{\prime }\right) \cdot \left( \exp \left( -k\tau \right)
\int_{Y_{2,d-k}}\frac{e_0^{*}\left( \phi _i\right) \exp \left( -\frac{p\tau }%
\hbar \right) }{-\hbar +c}{\cal E}_{2,d-k}^{\prime }\right) = 
\end{equation*}
\begin{equation*}
=\sum_{i=0}^n\frac{\prod_a\left( l_a\lambda _i-\mu
_a\right) }{\prod_{j\neq i}\left( \lambda _j-\lambda _i\right) }\sum_{d\geq
0}\sum_{k=0}^d\left( \exp \left( kt\right) \int_{Y_{2,k}}\frac{e_0^{*}\left(
\phi _i\right) \exp \left( \frac{pt}\hbar \right) }{\hbar +c}{\cal E}%
_{2,d}^{\prime }\right) \cdot 
\end{equation*}
\begin{equation*}
\cdot \left( \exp \left( \left( d-k\right) \tau \right) \int_{Y_{2,d-k}}%
\frac{e_0^{*}\left( \phi _i\right) \exp \left( -\frac{p\tau }\hbar \right) }{%
-\hbar +c}{\cal E}_{2,d-k}^{\prime }\right) = 
\end{equation*}
\begin{equation*}
=\sum_{i=0}^n\frac{\prod_a\left( l_a\lambda _i-\mu
_a\right) }{\prod_{j\neq i}\left( \lambda _j-\lambda _i\right) }%
\sum_{d,d^{\prime }\geq 0}\exp \left( dt\right) \int_{Y_{2,k}}\frac{%
e_0^{*}\left( \phi _i\right) \exp \left( \frac{pt}\hbar \right) }{\hbar +c}%
{\cal E}_{2,d}^{\prime }\cdot 
\end{equation*}
\begin{equation*}
\cdot \exp \left( d^{\prime }\tau \right) \int_{Y_{2,d^{\prime }}}\frac{%
e_0^{*}\left( \phi _i\right) \exp \left( -\frac{p\tau }\hbar \right) }{%
-\hbar +c}{\cal E}_{2,d^{\prime }}^{\prime }=\Phi \left(
t,\tau \right) \text{.} 
\end{equation*}
\begin{flushright}
$\Box$
\end{flushright}

\newpage

\end{document}